\documentclass[12pt, reqno]{amsart}
\usepackage[english]{babel}
\usepackage[utf8]{inputenc}
\usepackage{amsmath, amsthm, amscd, amssymb, fullpage, fancyhdr, upgreek, amssymb, hyperref, enumerate, comment, siunitx, enumitem, graphicx, mathrsfs, mathtools}
\usepackage[dvipsnames]{xcolor}
\usepackage{microtype}
\usepackage[toc]{appendix}
\usepackage[parfill]{parskip}

\newtheorem{theorem}{Theorem}[section]
\newtheorem{proposition}[theorem]{Proposition}

\newtheorem{lemma}[theorem]{Lemma}
\newtheorem{conjecture}[theorem]{Conjecture}

\theoremstyle{remark}
\newtheorem*{remark}{Remark}

\DeclareMathOperator{\sym}{sym}
\DeclareMathOperator{\USp}{USp}
\DeclareMathOperator{\SpO}{SO}
\DeclareMathOperator{\U}{U}
\DeclareMathOperator{\Prob}{Prob}
\DeclareMathOperator{\ad}{ad}
\DeclareMathOperator{\GL}{GL}
\DeclareMathOperator{\res}{Res}
\DeclareMathOperator{\vol}{vol}
\DeclareMathOperator{\eff}{eff}

\renewcommand{\mod}{~\mathrm{mod}~}

\usepackage{tikz}
\usepackage{tkz-tab}
\usepackage{tkz-graph}
\usetikzlibrary{shapes.geometric,positioning}

\numberwithin{equation}{section}
\renewcommand{\Im}{\mathrm{Im}}
\renewcommand{\Re}{\mathrm{Re}}
\renewcommand{\mod}{~\mathrm{mod}~}

\newcommand{\hr}[1]{\href{#1}{\url{#1}}}

\newcommand{\paren}[1]{\bigg( #1 \bigg)}

\title{A random matrix model for a family of cusp forms
}

\author{Owen Barrett, Zo\"e X. Batterman, Aditya Jambhale, Steven J. Miller, Akash L. Narayanan, Kishan Sharma, Chris Yao}

\address{Department of Mathematics, University of California, Berkeley, Berkeley, CA 94720}
\email{\textcolor{blue}{\href{mailto:barrett@math.berkeley.edu}{barrett@math.berkeley.edu}}}

\address{Department of Mathematics and Statistics, Pomona College, Claremont, CA 91711}
\email{\textcolor{blue}{\href{mailto:zxba2020@mymail.pomona.edu}{zxba2020@mymail.pomona.edu}}}

\address{Department of Pure Mathematics and Mathematical Statistics, University of Cambridge, Cambridge CB3 0WA, UK}
\email{\textcolor{blue}{\href{mailto:aj644@cam.ac.uk}{aj644@cam.ac.uk}}}

\address{Department of Mathematics and Statistics, Williams College, Williamstown, MA 01267}
\email{\textcolor{blue}{\href{mailto:sjm1@williams.edu}{sjm1@williams.edu}},  \textcolor{blue}{\href{Steven.Miller.MC.96@aya.yale.edu}{Steven.Miller.MC.96@aya.yale.edu}}}

\address{Department of Mathematics, University of Michigan, Ann Arbor, MI 48104}
\email{\textcolor{blue}{\href{mailto:anaray@umich.edu}{anaray@umich.edu}}}

\address{Department of Pure Mathematics and Mathematical Statistics, University of Cambridge, Cambridge CB3 0WA, UK}
\email{\textcolor{blue}{\href{mailto:kds43@cam.ac.uk}{kds43@cam.ac.uk}}}

\address{Department of Mathematics, Yale University, New Haven, CT 06520}
\email{\textcolor{blue}{\href{mailto:chris.yao@yale.edu}{chris.yao@yale.edu}}}

\thanks{This work was supported in part by NSF Grants DMS1561945 and DMS1659037, the University of Michigan, and Williams College. We thank many of our colleagues from other summers for helpful conversations. In particular, we thank Nathan Ryan and his colleagues for sharing their data and pre-print.}
\date{\today}
\begin{document}
\begin{abstract}
The Katz-Sarnak philosophy states that statistics of zeros of $L$-function families near the central point as the conductors tend to infinity agree with those of eigenvalues of random matrix ensembles as the matrix size tends to infinity. While numerous results support this conjecture, S. J. Miller observed that for finite conductors, very different behavior can occur for zeros near the central point in elliptic curve $L$-function families. This led to the creation of the excised model of Due\~{n}ez, Huynh, Keating, Miller, and Snaith, whose predictions for quadratic twists of a given elliptic curve are well fit by the data. The key ingredients are relating the discretization of central values of the $L$-functions to excising matrices based on the value of the characteristic polynomials at 1 and using lower order terms (in statistics such as the one-level density and pair-correlation) to adjust the matrix size. We extended this model for a family of twists of an $L$-function associated to a given holomorphic cuspidal newform of odd prime level and arbitrary weight. We derive the corresponding ``effective'' matrix size for a given form by computing the one-level density and pair-correlation statistics for a chosen family of twists, and we show there is no repulsion for forms with weight greater than 2 and principal nebentype. We experimentally verify the accuracy of the model, and as expected, our model recovers the elliptic curve model.
\end{abstract}
\maketitle
\tableofcontents
\section{Introduction}
The Katz-Sarnak philosophy states that the statistics of zeros of $L$-function families near the central point as the conductors tend to infinity agree with those of eigenvalues of random matrix ensembles as the matrix size tends to infinity \cite{KS99a}. While the general philosophy yields remarkable predictive insights for both local and global statistics, classic matrix ensembles fail to reflect finer statistical properties of $L$-function zeros. 

In 2006, Steven J. Miller observed that the elliptic curve $L$-function zero statistics for finite conductors deviated significantly from the scaling limit of the expected model of orthogonal matrices, though the fit improved as the conductors increased \cite{Mil06}. Subsequently, Due\~{n}ez, Huynh, Keating, Miller, and Snaith \cite{DHKMS12} created the \emph{excised orthogonal model} for finite conductors to more accurately reflect the phenomena Miller observed in the elliptic curve case. They point to the Kohnen-Zagier formula for holomorphic cusp forms of even weight \cite[Theorem 1]{KZ81} as motivation for the discretization of the central values of $L$-functions observed. For ease of reading, we recall the statement in its most general form below \cite[Theorem 1.5]{Mao08}.
\begin{theorem}
Let $f$ be a normalized Hecke eigenform of weight $2k$ and odd level $N$, $g$ a Shimura correspondence of $f$, and $L(s,f\otimes\psi_d)$ the $L$-function of $f$ twisted by the quadratic character $\psi_d$ with fundamental discriminant $d$. The formula of Kohnen and Zagier is
\begin{align}\label{eqn:kohnen-zagier-formula}
L(k,f \otimes \psi_d) \ = \ \frac{c(|d|)^2}{|d|^{k-1/2}}\kappa_f, \quad \text{ where } \kappa_f \ = \ \frac{\pi^k}{(k-1)!}\frac{\langle f,f \rangle}{\langle g,g \rangle}\kappa
,\end{align}
where we have taken the critical strip of the $L$-function to be $0< \Re(s) < 2k$, the function $\langle \cdot,\cdot \rangle$ denotes the Petersson inner product, $c(|d|)$ are the $d$\textsuperscript{th} Fourier coefficients of $g$, $\kappa$ is some product given in Theorem 1.5 of \cite{Mao08}, and $k$ is the central point in the evaluation of $L(s,f\otimes\psi_d)$.
\end{theorem}

Since the Fourier coefficients $c(|d|)$ of the modular form $g$ are integers, the Kohnen-Zagier formula implies the value at the central point is discretized on the order of $|d|^{k-1/2}$. Baruch-Mao \cite{BM07} extend the result to square-free odd level, and Mao \cite{Mao08} extends to arbitrary odd level, with respective modifications to $\kappa$.

To model the behavior of low-lying zeros using matrices, the authors in \cite{DHKMS12} build a model with two parameters. First, they find an analogous discretization of the values of the characteristic polynomials at 1 by introducing a \emph{cutoff value}. The model sieves off or excises those matrices whose characteristic polynomial has value less than the cutoff. The other key ingredient the authors consider is a modification of the ensemble's size. They consider two matrix sizes: one related to the mean density of zeros called the \textit{standard matrix size}, and the other determined from lower-order terms of the one-level density, called the \textit{effective matrix size}.

We extend the family of Due\~{n}ez, Huynh, Keating, Miller, and Snaith in \cite{DHKMS12} to a family of quadratic twists with finite conductor of a prescribed cuspidal newform $f$ of level an odd prime. For our chosen family, we motivate the choice of conductor of the twists depending on the prescribed form's duality and principality of its nebentype by analyzing the sign of the form's twisted $L$-function. Thus, we construct a family of twists $\mathcal{F}_f^{+}(X)$ up to a chosen conductor $X$, generalizing the family of even twists of a given elliptic curve considered in \cite{DHKMS12}.

Given a newform, one of our aims is to pinpoint the associated classical compact group for its family of twists. We use expansions of the statistics to improve the random matrix model's size. The following table summarizes the properties of the form that determine the corresponding matrix group.
\begin{center}
\begin{tabular}{c|c}
    Case & Group \\
    \hline
    $\chi_f$ principal, even twists & $\SpO(2N)$\\
    $\chi_f$ principal, odd twists & $\SpO(2N+1)$\\
    $\chi_f$ non-principal and $f=\overline{f}$ (self-CM) & $\USp(2N)$\\ 
    $\chi_f$ non-principal and $f\neq\overline{f}$ (generic) & $\U(N)$\\ 
\end{tabular}
\end{center}
In Section \ref{sec:one-level-density-unscaled-and-scaled}, we recover arithmetic coefficients from the lower-order terms of the one-level density statistic for those newforms that are non-generic, and we use them to construct the effective matrix size following \cite{DHKMS12}. Since the scaled one-level density for unitary matrices has no lower-order terms, we cannot use the same procedure to obtain an effective matrix size for the generic case. Hence, we must use another statistic. Guided by the discussion in \cite{DHKMS12} which followed the argument by Bogolomny, Bohigas, Leboeuf, and Monastra in \cite{BBLM06}, we turn to the pair-correlation statistic. Again, we assume the Ratios Conjectures to find a series expansion for pair-correlation. With a minor restriction on the arithmetic terms, we suggest an effective matrix size for generic forms. As an aside, we recovered H. L. Montgomery's 1973 conjecture \cite{Mon73} while generalizing from the zeta function to the generic forms studied in this paper.

The key ingredient to the creation of the cutoff value is the Kohnen-Zagier formula in \cite[Theorem 1]{KZ81}, which applies to those forms with level an odd number. We heuristically show that the Kohnen-Zagier formula implies a repulsion from the central point for those forms of weight 2 with principal nebentype in Section \ref{sct:cutoff_value_for_forms_with_principal_nebentypus}. We must emphasize that the Kohnen-Zagier formula implies that weight is the only controller for discretization. Since we heuristically predict that no repulsion occurs at the origin for forms with principal nebentype of weight 4 and greater, we only need to introduce the cutoff when we restrict our family to the elliptic curve case; that is, we follow the recipe given in \cite[Section 5.2]{DHKMS12} only for those principal forms with weight 2. To support our predictions, the authors in \cite{CSLPRRV24} demonstrate that repulsion does not occur for twists of the forms (LMFDB label) \texttt{3.6.a.a}, \texttt{3.8.a.a}, and \texttt{3.10.a.a}. By fixing the level and varying the weight, they provide an example that verifies the heuristics given by the Kohnen-Zagier formula.

\begin{center}
\begin{figure}[htpb]	
	\begin{tabular}{c c c}
		Lowest zeros ($\Delta=+1$) & Lowest zeros ($\Delta=-1$)\\
	\includegraphics[scale=0.45]{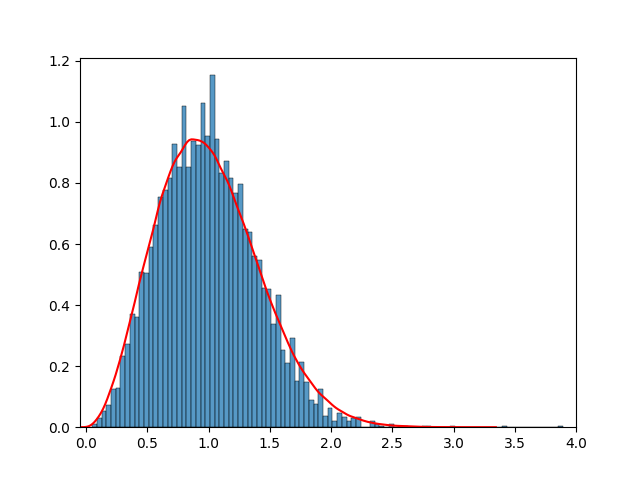} & \includegraphics[scale=0.45]{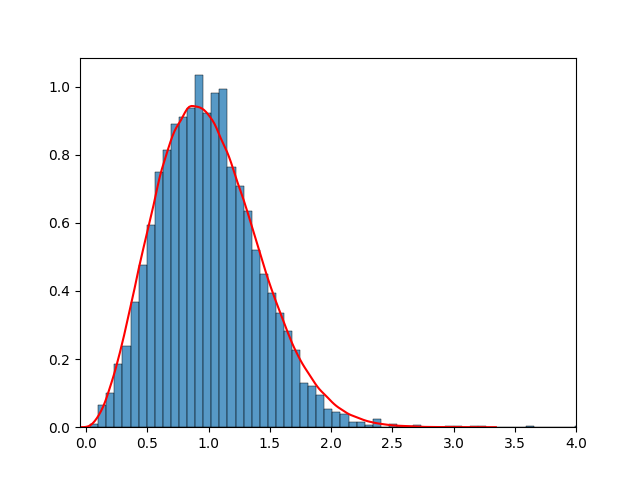} 
	\end{tabular}
	\caption{The left histogram shows the distribution of lowest-lying zeros for 5,458 twists of \texttt{3.7.b.a} with choice $\Delta =+1$ and discriminant up to 47,881, and the red line (left) shows the distribution of first eigenvalues of 1,000,000 randomly generated $\USp(20)$ matrices with characteristic polynomial evaluated at 1. The right histogram shows the distribution of lowest zeros for 5,726 twists of \texttt{3.7.b.a} with choice $\Delta =-1$ and discriminant up to 50,237, and the red line (right) shows the distribution of first eigenvalues of 1,000,000 randomly generated $\USp(20)$ matrices with characteristic polynomial evaluated at 1. The data have been normalized to have mean 1.}
 \label{tab:intro_zeros-mf3w7ba}
\end{figure}	
\end{center}

In Section \ref{sec:numerical-observations}, we gathered numerical data to verify the accuracy of our model. For each case, we select forms with low level and low weight to determine if the distribution of the lowest-lying zeros matches the eigenvalues of randomly generated matrices with characteristic polynomial evaluated at or near 1 from the predicted compact groups. For each form in the following table, we find agreement. For instance, in Figure \ref{tab:intro_zeros-mf3w7ba}, we superimpose the eigenvalues of random matrices from the ensemble $\USp(20)$ and the lowest-lying zeros of the self-CM form \texttt{7.3.b.a} and found agreement as predicted.
\begin{center}
\begin{tabular}{c|c|c}
    LMFDB Label & Type & Group \\
    \hline
    \texttt{11.2.a.a}, \texttt{5.4.a.a},\texttt{5.8.a.a}, \texttt{7.4.a.a} & $\chi_f$ principal, even twists & $\SpO(2N)$\\
    \hline
    \texttt{11.2.a.a}, \texttt{5.4.a.a},
    \texttt{5.8.a.a}, \texttt{7.4.a.a} & $\chi_f$ principal, odd twists & $\SpO(2N+1)$\\
    \hline
    \texttt{3.7.b.a}, \texttt{7.3.b.a} & self-CM & $\USp(2N)$\\
    \hline
    \texttt{13.2.e.a}, \texttt{11.7.b.b}, \texttt{7.4.c.a}, \texttt{17.2.d.a} & generic & $\U(N)$\\
\end{tabular}
\end{center}
We must highlight that we were not able to numerically find the effective matrix size for our family, except for the elliptic curve case as done in \cite{DHKMS12}. The computational difficulty unfortunately boils down to being unable to explicitly determine the associated $L$-function's Euler product used to calculate terms originating from the Ratios Conjectures \ref{conj:ratiosconjectures} and \ref{conj:ratios-lemma}.

We would like to remark that our numerical investigations in Section \ref{sec:numerical-observations} indicate that for most forms, the standard matrix size models exceptionally well the (non-vanishing) lowest-lying zeros of our family. For instance, in Figure \ref{tab:intro_zeros-mf3w7ba}, the random matrices from $\USp$ model extremely well twists of self-dual forms with non-principal nebentype. Using the model with the standard matrix size, we observe behavior previously unstudied. For instance, the non-vanishing lowest-lying zeros of a family of odd twists of a given form with principal nebentype with low weight (2 or 4) and even sign are ``attracted'' toward the origin (cf. Figure \ref{tab:intro_zeros-mf3w7ba}). This suggests a refinement of the random matrix model to incorporate this behavior. However, it is not clear what properties of a given form controls this attraction as the family of a principal form \texttt{5.4.a.a} with even sign does not show this attraction. We also find that the distribution of zeros of the family associated to certain generic forms (forms that are not self-dual) recover non-generic (self-CM) behavior. All of these surprising behaviors point to further theoretical investigations.
\begin{center}
\begin{figure}[htpb]	
	\begin{tabular}{c c}
		Lowest zeros (even twists) & Lowest zeros (odd twists) \\
	\includegraphics[scale=0.45]{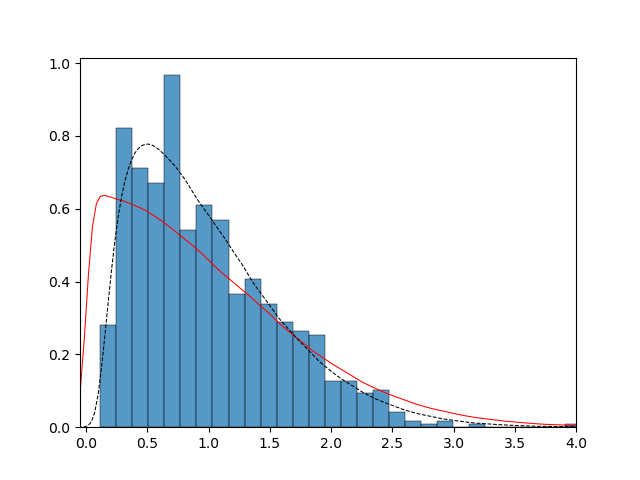} & \includegraphics[scale=0.45]{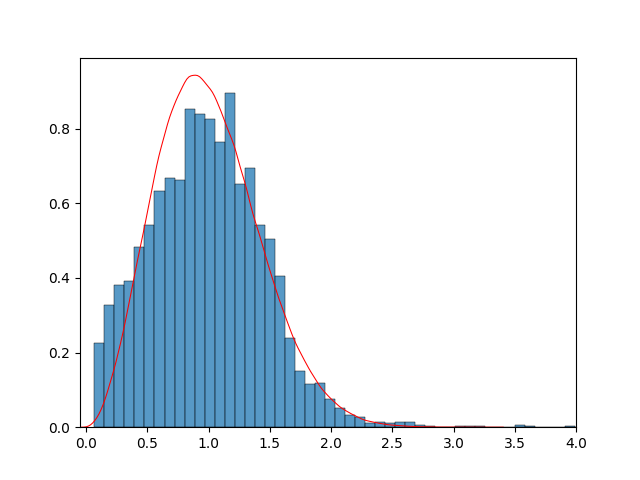}
	\end{tabular}
	\caption{The left histogram shows the distribution of lowest zeros for 1,380 even twists of \texttt{11.2.a.a} with discriminant up to 9,960; the red curve (left) shows the distribution of the first eigenvalues of 1,000,000 randomly generated $\text{SO}(20)$ matrices with characteristic polynomial evaluated at 1; the black dotted curve (left) is the same distribution but with excision. We varied the excision threshold numerically to obtain the optimal fit. The right histogram shows the distribution of lowest non-vanishing zeros for 4,563 odd twists of \texttt{11.2.a.a} with discriminant up to 32,897, and the red line (right) shows the distribution of the first eigenvalues of 1,000,000 randomly generated $\text{SO}(21)$ matrices with characteristic polynomial evaluated near 1. The data have been normalized to have mean 1.}
 \label{tab:intro_zeros-mf11w2aa}
\end{figure}	
\end{center}
\section{Background and notation}
We summarize some basic facts about random matrices and $L$-functions which will be vital for the paper's later sections.
\subsection{Matrix ensembles and one-level densities}
Denote by $G(N)$ any one of the compact matrix groups $\SpO(2N)$, $\SpO(2N+1)$, $\USp(2N)$, and $\U(N)$. Let $\varphi$ be an even Schwartz function. The unscaled one-level density $D(\varphi,G(N))$ formulas, found in \cite[Corollary AD.12.5.2]{KS99a}, are
\begin{align}
D(\varphi,\SpO(2N)) \ & \coloneqq \ \int_0^\pi \varphi(\theta)\left(\frac{2N-1}{2\pi} + \frac{\sin((2N-1)\theta)}{2\pi \sin\theta}\right)d\theta \\
D(\varphi,\SpO(2N+1)) \ & \coloneqq \ \int_0^{2\pi} \varphi(\theta)\left(\frac{N}{\pi} - \frac{\sin(2N\theta)}{2\pi\sin(\theta)} \right)d\theta \\
D(\varphi,\USp(2N)) \ & \coloneqq \ \int_0^{\pi} \varphi(\theta)\left(\frac{2N+1}{2\pi} -\frac{\sin((2N+1)\theta)}{2\pi \sin\theta} \right)d\theta \\
D(\varphi,\U(N)) \ & \coloneqq \ \int_0^{2\pi} \varphi(\theta)\left(\frac{N}{2\pi}\right)d\theta.
\end{align}
Since the mean spacing of the eigenangles of matrices in $G(N)$ depends on $N$, we may scale the eigenangles to have mean spacing one.
We obtain the asymptotic scaled one-level density for large $N$ by scaling the one-level density formulas and expanding them in powers of $1/N$ or $1/(2N+1)$:
\begin{align}
R(\varphi,\SpO(2N)) & \ \coloneqq \ \int_0^1 \varphi(\theta)\left(1-\frac{1}{2N} + \frac{\sin\left( (2N-1)\theta\pi /N \right) }{2N \sin(\theta\pi /N)} \right) d\theta \label{eq: scaledoneleveldensity_SO2N}\\
& \ = \ \int_0^1 \varphi(\theta)\left(1 + \frac{\sin(2\pi \theta)}{2\pi \theta} - \frac{1 + \cos(2\pi \theta)}{2N}- \frac{\pi \theta \sin(2\pi \theta)}{6N^2} + O(N^{-3})\right)d\theta \nonumber \\
R(\varphi,\SpO(2N+1)) & \ \coloneqq \ \int_0^1 \varphi(\theta)\left(1-\frac{1}{2N+1} - \frac{\sin(4\pi\theta N/(2N+1))}{(2N+1)\sin(2\pi\theta/(2N+1)} \right)d\theta \label{eq: scaledoneleveldensity_SO2N+1} \\
& \ = \ \int_0^1 \varphi(\theta) \left(1 - \frac{\sin(2\pi\theta)}{2\pi\theta} - \frac{1-\cos(2\pi\theta)}{2N+1} + \frac{2\pi\theta\sin(2\pi\theta)}{3(2N+1)^2} + O(N^{-3}) \right)d\theta \nonumber\\
R(\varphi,\USp(2N)) & \ \coloneqq \ \int_0^1 \varphi(\theta)\left(1 + \frac{1}{2N} - \frac{\sin\left( (2N+1)\theta \pi /N \right) }{2N \sin(\theta \pi /N)} \right)d\theta \\
& \ = \ \int_0^1 \varphi(\theta)\left(1- \frac{\sin (2\pi \theta)}{2\pi \theta} + \frac{1 -\cos(2\pi \theta)}{2N} + \frac{\pi \theta \sin (2\pi \theta)}{6N^2} + O(N^{-3}) \right)d\theta \nonumber\\
R(\varphi,\U(N)) & \ \coloneqq \ \int_0^1\varphi(\theta) d\theta
.\end{align}
\subsection{L-functions}
\label{sct:Lfunctions} 
We turn to cuspidal newforms and their associated $L$-functions. The following is adapted from \cite[Chapters 5.12, 5.13, and 14]{IK04}.

\subsubsection{Cuspidal newforms.}\label{sct:newforms-duals}
We consider the linear space $S_k(N,\chi_f)$ of cusp forms of level $N$, weight $k$, and nebentypus $\chi_f$ for the Hecke congruence subgroup $\Gamma_0(N)$. We focus on those cuspidal forms which are newforms.
In particular, if $f \in S_k^{\mathrm{new}}(M,\chi_f)$, then $f$ is an eigenform and has Fourier expansion $f(z)= \sum_{n = 1}^{\infty} a_f(n) e^{2\pi i nz}$
at the cusp $\infty$, where the $a_f(n)$'s are the Fourier coefficients. The normalized $L$-series $L(s,f)= \sum_{n\geq 1} \lambda_f(n)n^{-s}$ converges absolutely for $\mathrm{Re}(s)>1$ where we set $\lambda_f(n) \coloneqq a_f(n)n^{-s}$. Normalizing allows us to relate $s$ to $1-s$ through the functional equation of the completed $L$-series. We then analytically continue the $L$-series to the complex plane to construct the $L$-function associated to $f$ which has Euler product
\begin{align}
	L(s,f) \ = \ \prod_p \left( 1- \lambda_f(p)p^{-s} + \chi_f(p) p^{-2s} \right) ^{-1} \ = \ \prod_p \left( 1 - \alpha_f(p)p^{-s} \right) ^{-1} \left( 1 - \beta_f(p) p^{-s} \right) ^{-1},
\end{align}
where the Satake parameters $\alpha_f,\beta_f$ satisfy 
\begin{align}\label{eq:satakeidentity}
\alpha_f(p) + \beta_f(p) = \lambda_f(p) \quad \text{and} \quad \alpha_f(p) \beta_f(p) = \chi_f(p)
\end{align}
\cite[Proposition 1.3.6]{Bump97}. By comparing the coefficient of the $p^{-ms}$ term of the Euler product with the $L$-series, we obtain the relation
\begin{align} \label{eqn:lambda-coefficient-relation}
\lambda_f(p^{m}) \ = \ \sum_{\ell \geq 0} \alpha_f(p)^{\ell} \beta_f(p)^{m-\ell}
.\end{align}
The form $\overline{f}$ dual to $f$ has Fourier coefficients which satisfy the duality relation $\lambda_f(n) = \chi_f(n) \lambda_{\overline{f}}(n) = \chi_f(n) \overline{\lambda_{f}}(n)$ for $\mathrm{gcd}(n,M)=1$ by the adjointness formula for a cuspidal Hecke form $f$ \cite[Proposition 14.11]{IK04}.
We may relate $L(s,f)$ to $L(1-s,\overline{f})$ by the functional equation
\begin{align}
	L(s,f) \ = \ \epsilon_f \bigg( \frac{\sqrt{M} }{2\pi} \bigg) ^{1-2s} \frac{\Gamma\left( \frac{k+1}{2} - s \right) }{\Gamma\left( \frac{k-1}{2} + s \right) } L(1-s,\overline{f})
,\end{align}
where $\epsilon_f$ is the root number associated to $f$ and has absolute value 1.

\subsubsection{Rankin-Selberg convolution}\label{sct:rankin-selberg}
For $f\in S_k^{\mathrm{new}}(M_f,\chi_f)$ and $g\in S_\kappa^{\mathrm{new}}(M_g,\chi_g)$, the $L$-series of their Rankin-Selberg convolution is defined to be 
\begin{equation}
L(s,f\otimes g)  \ \coloneqq\ L(2s,\chi_f\chi_g) 
 \sum_{n\geq 1} \lambda_f(n) \lambda_g(n)n^{-s}, \end{equation} provided the least common multiple of $M_f$ and $M_g$ is square-free. It is well-known that the $L$-series $L(s,f\otimes g)$ has analytic continuation and admits an Euler product \cite[Theorem 1.6.2]{Bump97}.
The local factor at unramified primes (those primes not dividing the level of $f$ or $g$) is given by
\begin{align}
L_p(s,f\otimes g) & \ = \ 1 - \lambda_f(p)\lambda_g(p)p^{-s} + \left(\chi_f(p)\lambda_g(p)^2 + \chi_g(p)\lambda_f(p)^2 - 2 \chi_f(p) \chi_g(p) \right)p^{-2s} \nonumber \\
	& \qquad \ - \ \chi_f(p)\chi_g(p)\lambda_f(p)\lambda_g(p)p^{-3s} + \chi_f(p)^2\chi_g(p)^2p^{-4s} 
.\end{align}
\subsubsection{Quadratic twists}
An integer $d$ is a fundamental discriminant provided that $d$ is either square-free and congruent to 1 modulo 4 or is four times a square-free integer congruent to 2 or 3 modulo 4. Let $L(s,f_d) \coloneqq L(s,f\otimes \psi_d)$ denote the $L$-function obtained by twisting $L(s,f)$ by a quadratic character $\psi_d$ with fundamental discriminant $d$, that is, conductor $|d|$.
For $\mathrm{Re}(s)>1$, the twisted $L$-function $L(s,f_d) \ = \ \sum_{n\geq 1} \lambda_f(n) \psi_d(n)n^{-s}$ has Euler product \begin{align}\label{eqn:euler-twisted} 
L(s,f_d) \ = \ \prod_p \left( 1 - \lambda_f(p)\psi_d(p)p^{-s} + \chi_f(p)\psi_d(p)^2p^{-2s} \right) ^{-1}. 
\end{align}
Provided $\mathrm{gcd}(d,M) = 1$, the completed $L$-function satisfies the functional equation\footnote{\cite[Section 14.8]{IK04}.}
\begin{align}\label{eqn:twisted-functional-eq}
L(s,f_d) \ = \ \epsilon_{f\otimes \psi_d} \bigg( \frac{\sqrt{M} |d|}{2\pi} \bigg) ^{1-2s} \frac{\Gamma\left( \frac{k+1}{2}-s \right) }{\Gamma\left(\frac{k-1}{2} + s \right) } L(1-s,\overline{f}_d)
,\end{align}
with root number $\epsilon_{f\otimes \psi_d} = \chi_f(d) \psi_d(M) \tau(\psi_d)^2/d\epsilon_f =  \chi_f(d)\psi_d(-|D|) \epsilon_f$ where $\tau(\psi)$ denotes the Gauss sum and the second equality follows from \cite[Corollary 2.1.47]{Coh07} since $\psi_d=\left(\frac{d}{\cdot}\right)$ is a real character.
The approximate functional equation for a twist shifted by $\alpha$ is given in \cite[Section 5.2]{IK04} by
\begin{align}\label{eqn:approx-functional-eq-twisted}
L(1 /2+\alpha,f_d) & \ = \ \sum_{m < x} \frac{\lambda_f(m)\psi_d(m)}{m^{1 /2+\alpha}} \\
& \quad \ + \ \omega_f(d) \epsilon_f \bigg( \frac{\sqrt{M} |d|}{2\pi} \bigg) ^{-2\alpha} \frac{\Gamma\left(\frac{k}{2} - \alpha\right)}{\Gamma\left(\frac{k}{2} + \alpha\right)} \sum_{n < y} \frac{\overline{\lambda_f(n)}\psi_d(n)}{n^{1 /2-\alpha}} + \text{remainder}, \nonumber
\end{align}
where $xy=d^2 /2\pi$. We ignore the remainder term since the recipe given in~\cite{CFKRS05} instructs to discard the remainder in the
approximate functional equation, and to replace the root numbers and summands
with their expectations over the family.
\subsubsection{Complex multiplication and self-duality}\label{sct:CM} 
We say the newform $f$ has complex multiplication by $\eta$ if $\eta(p)\lambda_f(p) = \lambda_f(p)$ for all primes $p$ in a set of primes of density 1 \cite{Rib77}.  Those self-dual newforms with non-principal nebentype have complex multiplication by their nebentype \cite[Section 3, Remark 2]{Rib77}. Shimura constructed examples of integral-weight cusp forms with complex multiplication by their own nebentype. He constructed such forms from imaginary quadratic fields $K = \mathbb{Q}(\sqrt{-D} )$ with class number 1 and with unit group $\{\pm 1\} $. In particular, the associated modulus of $K$ is the entire ring of integers $\mathcal{O}_K$. We call such forms constructed by Shimura ``self-CM.'' Below, we document a property applicable to twists of a given self-CM form.

\begin{proposition}
	Consider the $L$-function $L(s,f)$ attached to a self-CM cuspidal newform $f$ with sign $\epsilon_f$. Then for each positive fundamental discriminant $d$ prime to $M$, the sign $\epsilon_{f \otimes \psi_d}$ of the twisted $L$-function $L(s,f_d)$ is equal to $\epsilon_f$.
\end{proposition}
\begin{proof} 
Assume $f$ is a self-CM form.
Then $f$ has level $M = |D|$ \cite[Theorem 12.5]{Iwa97}.
By \cite[Section 3, Remark 1]{Rib77}, we write the nebentype of the twisted form as $\chi_{f\otimes\psi_d}=\chi_{f}=\psi_D$. Thus, the root number of the twisted $L$-function in \ref{eqn:twisted-functional-eq} is
\begin{align}\label{eq:quadrecip} \epsilon_{f\otimes\psi_d}=\left(\frac{D}{d}\right)\left(\frac{d}{|D|} \right)\epsilon_f =(-1)^{(D'-1)(d'-1)/4}\epsilon_f,
\end{align}
where $D'$, $d'$ denote the odd parts (maximal odd divisors with the same sign) of $D$ and $d$, respectively. Since $D$ and $d$ are both fundamental discriminants and share no prime factors, one of $\left\{D,d\right\}$ must be odd and therefore congruent to $1$ mod $4$. We conclude $\epsilon_{f\otimes\psi_d}=\epsilon_f$.
\end{proof}

We introduce some notation. Let $S_{k}^{\text{new}}(M,\text{principal})$ denote the family of forms with $\chi_f$ principal; let $S_{k}^{\text{new}}(M,\text{self-CM})$ denote the family of self-dual forms with $\chi_f$ non-principal; and let $S_{k}^{\text{new}}(M,\text{generic})$ denote the family of generic forms with $\chi_f$ non-principal and $f \neq \overline{f}$.

\subsubsection{Symmetric and adjoint square L-functions.}\label{sct:sym-adj-sq}
We follow \cite[Section 5.12]{IK04} and define symmetric and adjoint square $L$-functions of newforms as factors of Rankin-Selberg convolutions. Let $\chi'_f$ denote the primitive character that induces the nebentypus $\chi_f$ of $f$. Denote $L(s,\sym^2f) \coloneqq L(s,f\otimes f)L(s,\chi'_f)^{-1}$ and $L(s,\ad^2f) \coloneqq L(s,f\otimes \overline{f})\zeta(s)^{-1}$. The symmetric square $L$-function $L(s,\sym^2f)$ has an Euler product with local factors at unramified primes given by
\begin{align}\label{eqn:local-factor-sym-sq}
L(s,\sym^2f) = (1-\alpha_f(p)^2p^{-s})^{-1}(1-\alpha_f(p)\beta_f(p)p^{-s})^{-1}(1-\beta_f(p)^2p^{-s})^{-1}
.\end{align}
We record analytic facts about these functions. First, note that $L(s,\ad^2f)$ is always entire since $L(s,f\otimes \overline{f})$ always has a simple pole at $s=1$ which cancels with the zero of $\zeta(s)$. For the upper half-plane $\mathcal{H}$, we have an arithmetically significant value at $s=1$ for $L(s,\ad^2f)$ by the equality 
\begin{align} \label{eqn:residue-at-1-f-bar-f}
	L(1,\ad^2f) \ = \ \frac{(4\pi)^k \langle f, f\rangle}{\Gamma(k)\vol (\Gamma_0(M) \backslash \mathcal{H})} \ = \ \res (L(s,f\otimes\overline{f}),1 )
.\end{align}
When $f \in S_{k}^{\text{new}}(M,\text{principal})$, then $\overline{f} = f$ and $L(s,\ad^2f) = L(s,\sym^2f)$ by $\chi'_f \equiv 1 \mod M$.
When $f$ has non-principal nebentypus, $L(s,\chi_f)$ is entire. The $L$-function $L(s,f\otimes f)$ may not be entire since $f$ may still be self-dual. If $f$ is self-CM, $L(s,\sym^2f)$ has a pole, and the symmetric and adjoint square $L$-functions do not coincide since $L(s,\chi'_f) \neq \zeta(s)$. In fact, the case of $f$ self-CM is the only case when $L(s,\sym^2f)$ inherits the pole from $L(s,f\otimes f)$ at $s = 1$ as $L(s,\chi_f')$ is entire.
\subsection{The set of ``good" fundamental discriminants}\label{sct:fundamental-discriminant-family}
Our goal is to determine a set of ``good'' fundamental discriminants for each form depending on its principality and self-duality so that we may apply the ratios recipe (see Section \ref{subsect:ratios_conjectures}). In particular, we seek fundamental discriminants which permit us to make useful replacements of $\omega_f(d)\epsilon_f$ and $\psi_d(-M)$ in the functional equation \eqref{eqn:approx-functional-eq-twisted}. Before we start, we insist that the level $M$ of the cuspidal newform be an odd prime.

We may choose even or odd twists of any given cuspidal newform of principal nebentype. Mao's generalization in \cite{Mao08} of the Kohnen-Zagier formula allows us to consider both positive and negative fundamental discriminants. For simplicity, we define our family using twists with positive discriminant.

We turn to self-dual forms with non-principal nebentype constructed by Shimura. Since we are unable to change the sign of such forms with quadratic twists by fundamental discriminants coprime to the level (see Section \ref{sct:CM}), we choose to restrict our family by twists with positive discriminants coprime to both the level and the sign. Recall that we only consider those self-CM forms that have positive sign.
Regardless, we may choose any parity for the twists.

In the generic case, there is no longer any notion of parity for the functional equation \eqref{eqn:approx-functional-eq-twisted}. The only constraint we impose is that our choice of fundamental discriminants is positive.

The above discussion gives rise to a notion of a set of ``good" fundamental discriminants associated to a form $f$. Let $\mathcal{D}^{+}$ denote the set of positive fundamental discriminants (one may choose to do the same for negative fundamental discriminants which should yield the same model), and fix a cuspidal newform $f \in S_k^{\text{new}} (M,\chi_f)$. Fix an integer $\Delta \in \{+1,-1\}$. The family from which we twist a given form by is
\begin{align}\label{D-eqn}
	\mathcal{D}_f^{+}(X) \ \coloneqq \ \begin{cases}
\{d \in \mathcal{D}^{+} \mid 0<d\leq X,\ \ \psi_d(-M) \epsilon_f \ = \ +1\} & \qquad \text{$\chi_f$ principal, even twists}, \\
\{d \in \mathcal{D}^{+} \mid 0<d\leq X,\ \ \psi_d(-M) \epsilon_f \ = \ -1\} & \qquad \text{$\chi_f$ principal, odd twists}, \\
\{d \in \mathcal{D}^{+} \mid 0<d\leq X,\ \ \psi_d(-M) \ = \ \Delta\} & \qquad \text{$f$ self-CM}, \\
\{d \in \mathcal{D}^{+} \mid 0 < d \leq X\} & \qquad \text{$f$ generic}
.\end{cases}
\end{align}
We state estimates on the cardinality $|\mathcal{D}_f^{+}(X)|$ of the family (c.f. Appendix \ref{appendix:counting_fundamental_discriminants}):
\begin{align}\label{eqn:cardinality-estimates} 
|\mathcal{D}_f^{+}(X)| \ = \ \begin{cases}
3MX(2\pi^2(M+1))^{-1} + O(X^{1 /2}) & \qquad f = \overline{f}, \\
3MX(\pi^2(M^2-1))^{-1} + O(X^{1 /2}) & \qquad f \neq \overline{f}.
\end{cases}
\end{align}

\subsection{The family of quadratic twists}
Denote the family of quadratic twists of a fixed holomorphic cuspidal newform $f$ with fundamental discriminants ranging over $\mathcal{D}_f^{+}(X)$ by $\mathcal{F}_f^{+}(X)$.
For $d \in \mathcal{D}_f^{+}(X)$, the quadratic character $\psi_d(M) = (d/M)$ assumes the value
\begin{align}
\mathcal{E}_f(M) \ \coloneqq \ \psi_d(M) \ = \ \begin{cases}
-\epsilon_f & \qquad \chi_f \text{ principal, even twists}, \\
\ \ \epsilon_f & \qquad \chi_f \text{ principal, odd twists}, \\
\ \ -\Delta & \qquad \chi_f \text{ non-principal, } f = \overline{f}, \\
(d/M) & \qquad \chi_f \text{ non-principal, } f \neq \overline{f}.
\end{cases}
\end{align}

\subsection{Ratios Conjectures}\label{subsect:ratios_conjectures}
The work of Conrey, Farmer, and Zirnbauer \cite{CFKRS05} formalized an observation by Nonnenmacher and Zirnbauer into the corresponding Ratios Conjectures for $L$-functions.\footnote{This observation happened at a workshop at MSRI in 1999 \cite[Introduction]{CFZ08}.} We must formulate our own form of the Ratios Conjectures for our family so that we may derive a formula for the one-level density of the zeros near the critical line $1/2$ of $L$-functions associated to $\mathcal{F}_f^{+}(X)$. We consider the average over the family of ``good" fundamental discriminants of a ratio of shifted $L$-functions:
\begin{align} \label{eqn:ratio-shifted}
R_f(\alpha,\gamma) \ \coloneqq \ \sum_{d \in \mathcal{D}_f^{+}(X)} \frac{L(1/2 + \alpha,f_d)}{L(1/2+ \gamma,f_d)}
.\end{align}
Using~\eqref{eqn:euler-twisted}, we have
\begin{align}\label{eq:recipseries}
\frac1{L(s,f_d)}\ =\ \sum_{n\geq1}\mu_f(n)\psi_d(n)n^{-s},
\end{align}
where $\mu_f(n)$ is a multiplicative function defined by
\begin{align} \mu_f(n)  \ = \ \begin{cases}
  -\lambda_f(p) & \text{if $n=p$,} \\
  \chi_f(p) & \text{if $n=p^2$,} \\
  0 & \text{if $n=p^j$, $j>2$.}
\end{cases} \label{eq:mudef}
\end{align}
We denote the first sum arising from the approximate
functional equation~\eqref{eqn:approx-functional-eq-twisted} by
\begin{align}\label{eq:R1def}R_f^1(\alpha,\gamma)\ \coloneqq \ \sum_{d \in \mathcal{D}_f^{+}(X)}
\sum_{h,m\geq 0}\frac{\lambda_f(m)\mu_f(h)\psi_d(mh)}{m^{1/2+\alpha}h^{1/2+\gamma}}.
\end{align}
Likewise for $R_f^2(\alpha,\gamma)$:
\begin{align}\label{eq:R2def}R_f^2(\alpha,\gamma)\ \coloneqq \
\omega_f(d)\epsilon_f\frac{\Gamma\left(\frac k2-\alpha\right)}{\Gamma\left(\frac k2+\alpha\right)}
\sum_{d \in \mathcal{D}_f^{+}(X)}\left(\frac{\sqrt M |d|}{2\pi}\right)^{-2\alpha}
\sum_{h,m\geq 0}\frac{\overline{\lambda_f(m)}\mu_f(h)\psi_d(mh)}
{m^{1/2-\alpha}h^{1/2+\gamma}}.
\end{align}
Our formulation of the Ratios Conjectures is the following.\footnote{For other examples of the Ratios Conjectures, see \cite{CS07}.}
\begin{conjecture}\label{conj:ratiosconjectures}
For the conditions $-1/4 < \mathrm{Re}(\alpha) < 1/4$, $1/\log x \ll \mathrm{Re}(\gamma) < 1/4$ and $\mathrm{Im}(\alpha),\Im(\gamma) \ll X^{1-\varepsilon}$, the average over the family of a ratio of shifted $L$-functions is
\begin{align}\label{eqn:R_f(alpha,gamma)-error}
	R_f(\alpha,\gamma) & \ = \ \sum_{d\in \mathcal{D}_f^{+}(X)} \left[ Y_f A_f(\alpha,\gamma) + \eta_f \bigg( \frac{\sqrt{M} |d|}{2\pi} \bigg)^{-2\alpha}  \frac{\Gamma\left(k/2 - \alpha\right)}{\Gamma\left(k/2 + \alpha\right)} \widetilde{Y}_f \widetilde{A}_f (-\alpha,\gamma) \right] \\
& \qquad \qquad \qquad \qquad \qquad \qquad \qquad \qquad \qquad \qquad \qquad + \ O(X^{1/2+\varepsilon}), \nonumber
\end{align}
where 
\begin{align}
Y_f(\alpha,\gamma) & \ = \ \frac{L(1 + 2\gamma,\chi'_f)L(1+2\alpha,\sym^2f)}{L(1+\alpha+\gamma,\chi'_f)L(1+\alpha+\gamma,\sym^2f)}, \label{eqn:Y_f(alpha,gamma)} \\ 
\widetilde{Y}_f(-\alpha,\gamma) & \ = \ \frac{L(1+2\gamma,\chi'_f) L(1-2\alpha,\sym^2\overline{f})}{\zeta(1-\alpha+\gamma)L(1-\alpha+\gamma,\ad^2f)}, \label{eqn:tilde-Y_f(alpha,gamma)}
\end{align}
\begin{align} \label{eqn:A_f(alpha,gamma)}
A_f(\alpha,\gamma) & \ = \ Y_f(\alpha,\gamma)^{-1} V_{\mid}(\alpha,\gamma) V_{\nmid}(\alpha,\gamma)
,\end{align}
\begin{align}\label{eqn:V_mid(alpha,gamma)}
V_{\mid}(\alpha,\gamma) \ = \ \sum_{m=0}^{\infty} \left( \frac{\lambda_f(M^{m}) \mathcal{E}_f^{m}(M)}{M^{m(1 /2+\alpha)}} - \frac{\lambda_f(M)}{M^{1 /2 + \gamma}} \frac{\lambda_f(M^{m}) \mathcal{E}_f^{m+1}(M)}{M^{m(1 /2+\alpha)}} \right)
,\end{align}
\begin{align} \label{eqn:V_nmid(alpha,gamma)}
V_{\nmid}(\alpha,\gamma) \ = \ \prod_{p\nmid M} \left( 1 + \frac{p}{p+1} \left( \sum_{m=1}^{\infty} \frac{\lambda_f(p^{2m})}{p^{m(1+2\alpha)}}- \frac{\lambda_f(p)}{p^{1+\alpha+\gamma}}\sum_{m=0}^{\infty} \frac{\lambda_f(p^{2m+1})}{p^{m(1+2\alpha)}} + \frac{\chi_f(p)}{p^{1+2\alpha}} \sum_{m=0}^{\infty} \frac{\lambda_f(p^{2m})}{p^{m(1+2\alpha)}} \right)  \right),
\end{align}
\begin{align} \label{eqn:tilde-A_f(alpha,gamma)}
\widetilde{A}_f(-\alpha,\gamma) \ & = \ \widetilde{Y}_f(-\alpha,\gamma)^{-1} \widetilde{V}_{\mid}(-\alpha,\gamma) \widetilde{V}_{\nmid}(-\alpha,\gamma),
\end{align}
\begin{align}\label{eqn:widetilde_V_mid(alpha,gamma)}
\widetilde{V}_{\mid}(-\alpha,\gamma) \ = \ \sum_{m=0}^{\infty} \left( \frac{\overline{\lambda_f(p^{m})} \mathcal{E}_f^{m}(p)}{p^{m(1 /2-\alpha)}} - \frac{\lambda_f(p)}{p^{1 /2 + \gamma}} \frac{\overline{\lambda_f(p^{m})} \mathcal{E}_f^{m+1}(p)}{p^{m(1 /2-\alpha)}} \right) ,
\end{align}
\begin{align}\label{eqn:widetilde_V_nmid(alpha,gamma)}
\widetilde{V}_{\nmid}(-\alpha,\gamma) \ = \ \prod_{p\nmid M} \left( 1 + \frac{p}{p+1} \left( \sum_{m=1}^{\infty} \frac{\overline{\lambda_f(p^{2m})}}{p^{m(1-2\alpha)}}- \frac{\lambda_f(p)}{p^{1-\alpha+\gamma}}\sum_{m=0}^{\infty} \frac{\overline{\lambda_f(p^{2m+1})}}{p^{m(1-2\alpha)}} + \frac{\chi_f(p)}{p^{1-2\alpha}} \sum_{m=0}^{\infty} \frac{\overline{\lambda_f(p^{2m})}}{p^{m(1-2\alpha)}} \right)  \right),
\end{align}
and the expectation of $\omega_f(d)\epsilon_f$ over $d$ is
\begin{align}\label{eq:rootexpect}
\eta_f \ \coloneqq \ \langle\omega_f(d)\epsilon_f\rangle \ = \ \begin{cases}
+1 & \qquad \text{$\chi_f$ principal, even twists}, \\
-1 & \qquad \text{$\chi_f$ principal, odd twists}, \\
+1 & \qquad \chi_f \text{ non-principal, $f = \overline{f}$},\\ 
\langle\omega_f(d)\epsilon_f\rangle & \qquad \text{$\chi_f$ non-principal, $f \neq \overline{f}$}.
\end{cases}
\end{align}
\end{conjecture}
The authors in \cite{HKS09} follow the recipe outlined in \cite{CFKRS05}, \cite{CFZ08} and the calculations in \cite{CS07} to derive a formula for \eqref{eqn:R_f(alpha,gamma)-error}. Since $\lambda_f(n)$ is a multiplicative function, we may express $R^1_f(\alpha,\gamma)$ as the Euler product
\begin{align}
R^1_f(\alpha,\gamma) \ \approx \ |\mathcal{D}_f^{+}(X)| V_{\mid}(\alpha,\gamma) V_{\nmid}(\alpha,\gamma)
\end{align}
where
\begin{align}
V_{\mid}(\alpha,\gamma) \ \coloneqq \ \prod_{p\mid N} \left( \sum_{h,m\geq 0} \frac{\lambda_f(p^m)\mu_f(p^h)\mathcal{E}_f^{m+h}(p)}{p^{m(1/2+\alpha)+h(1/2+\gamma)}} \right), \\
V_{\nmid}(\alpha,\gamma) \ \coloneqq \ \prod_{p\nmid N} \left( 1 + \frac{p}{p+1} \sum_{\substack{m,h\geq 0\\m+h\text{ even}}} \frac{\lambda_f(p^m)\mu_f(p^h)}{p^{m(1/2+\alpha)+h(1/2+\gamma)}} \right)
.\end{align}
The definition \eqref{eq:mudef} requires that we need only consider $h=0,1$ for~\eqref{eqn:V_mid(alpha,gamma)} and $h=0,1,2$ for~\eqref{eqn:V_nmid(alpha,gamma)}. We are left with
\begin{align}
V_{\mid} \ & = \ \prod_{p\mid M} \left( \sum_{m=0}^{\infty} \left( \frac{\lambda_f(p^m)\mathcal{E}_f^m(p)}{p^{m(1/2+\alpha)}} - \frac{\lambda_f(p)\lambda_f(p^m)\mathcal{E}_f^{m+1}(p)}{p^{m(1/2+\alpha)+(1/2+\gamma)}}
\right)\right), \\
V_{\nmid} \ & = \ \prod_{p\nmid M} \left(1+\frac{p}{p+1}\left(\sum_{m=1}^{\infty} \frac{\lambda_f(p^{2m})}{p^{m(1+2\alpha)}} - \frac{\lambda_f(p)}{p^{1+\alpha+\gamma)}}\sum_{m=0}^{\infty} \frac{\lambda_f(p^{2m+1})}{p^{m(1+2\alpha)}} + \frac{\chi_f(p)}{p^{1+2\gamma}} \sum_{m=0}^{\infty} \frac{\lambda_f(p^{2m})}{p^{m(1+2\alpha)}} \right)\right) \label{eqn:euler-prod-V-nmid}
.\end{align}
Recall the local Euler factors at unramified primes of the symmetric square $L$-function $L(s,\sym^2f)$ given in~\eqref{eqn:local-factor-sym-sq}. We wish to rewrite the Euler product for $L(s,\sym^2f)$ in terms of the Fourier coefficients of the given form. 
Using~\eqref{eq:satakeidentity} and~\eqref{eqn:lambda-coefficient-relation}, we have that the local factor of $L(s,\sym^2f)$ at an unramified prime, is given by
\begin{align}\label{eq:symlocalfactor}
  \left(1-\frac{\lambda_f(p)^2-\alpha_f(p)\beta_f(p)}{p^s}
  + \frac{\alpha_f(p)\beta_f(p)(\lambda_f(p)^2-\alpha_f(p)\beta_f(p))}{p^{2s}}
  - \frac{\alpha_f(p)^3\beta_f(p)^3}{p^{3s}}\right)^{-1}.
\end{align}
We impose the restrictions $-1/4 < \mathrm{Re}(s) <1/4$ and $\log X \ll \mathrm{Re}(s) \ll 1/4$ to control the convergence of the Euler product in \ref{eqn:euler-prod-V-nmid}. With these restrictions, we rewrite \ref{eqn:euler-prod-V-nmid} as follows:
\begin{align}
V_{\nmid} \ & = \ \prod_{p\nmid M} \left( 1+\frac{\lambda_f(p^2)}{p^{1+2\alpha}}-\frac{\lambda_f(p^2)+\chi_f(p)}{p^{1+\alpha+\gamma}}+\frac{\chi_f(p)}{p^{1+2\gamma}}+\cdots \right),
\end{align}
where the $\cdots$ indicate terms that converge like $1/p^2$ for restricted $\alpha$ and $\gamma$. We use the following approximation identities to factor out the divergent or slowly convergent terms. By Equation \eqref{eq:symlocalfactor}, the local factor of $L(1+2\alpha,\sym^2f)$ at an unramified prime is
\begin{align}
\left( 1-\frac{\lambda_f(p^2)}{p^{1+2\alpha}} +\frac{\lambda_f(p^2)\chi_f(p)}{p^{2(1+2\alpha)}}-\frac{\chi_f(p)^3}{p^{3(1+2\alpha)}}\right)^{-1} = \left( 1 + \frac{\lambda_f(p^2)}{p^{1+2\alpha}} + \cdots \right),
\end{align}
where $\cdots$ indicate those terms which converge under the imposed conditions. We recall the local factor at an unramified prime of $L(1+\alpha+\gamma,\sym^2f)^{-1}L(1+\alpha+\gamma,\chi_f)^{-1}$ coincides with that of $L(1+\alpha+\gamma,f\otimes f)^{-1}$ by definition: this factor is
\begin{align}
\left(1 - \frac{\lambda(p^2)+\chi_f(p)}{p^{1+\alpha+\gamma}} + \cdots \right)
.\end{align}
We account for the last factor $\chi_f(p)/p^{1+2\gamma}$. If $\chi_f$ is non-principal, this term is convergent since $L(s,\chi_f)$ is entire. If $\chi_f$ is principal, we note both $L(s,\chi_f)$ and $L(s,\chi_f')$ have the same local factors at unramified primes. Thus, we let a factor of $L(1+2\gamma,\chi_f')$ account for the divergence of the term $\chi_f(p)/p^{1+2\gamma}$. 

We obtain the other sum $R_f^2(\alpha,\gamma)$ by using the second term of the approximate functional equation and carrying out the same steps as for $R_f^1(\alpha,\gamma)$. We find analogous expressions for $V_{\mid}(\alpha,\gamma)$ and $V_{\nmid}(\alpha,\gamma)$, which we denote by $\tilde{V}_{\mid}(\alpha,\gamma)$ and $\tilde{V}_{\nmid}(\alpha,\gamma)$, respectively:
\begin{align}
	\tilde{V}_{\mid}(\alpha,\gamma) \ & =\ \prod_{p\mid M} \left( \sum_{m=0}^{\infty} \left( \frac{\overline{\lambda_f(p^{m})\mathcal{E}_f^{m}(p)}}{p^{m(1 /2-\alpha)}} - \frac{\lambda_f(p)\overline{\lambda(p^{m})}\mathcal{E}_f^{m+1}(p)}{p^{m(1 /2-\alpha)+1 /2+\gamma}} \right)  \right)  \\
	\tilde{V}_{\nmid}(\alpha,\gamma) \ & =\ \prod_{p\nmid M} \left( 1 + \frac{p}{p+1} \left( \sum_{m=1}^{\infty} \frac{\overline{\lambda_f(p^{2m})}}{p^{m(1-2\alpha)}} - \frac{\lambda_f(p)}{p^{1-\alpha+\gamma}} \sum_{m=0}^{\infty} \frac{\overline{\lambda_f(p^{2m+1})}}{p^{m(1-2\alpha)}} + \frac{\chi_f(p)}{p^{1+2\gamma}}\sum_{m=0}^{\infty} \frac{\overline{\lambda_f(p^{2m})}}{p^{m(1-2\alpha)}}\right)  \right) \\
					 & =\ \prod_{p\nmid M}\left( 1 + \frac{\overline{\lambda_f(p^2)}}{p^{1-2\alpha}} - \frac{|\lambda_f(p)|^2}{p^{1-\alpha+\gamma}} + \frac{\chi_f(p)}{p^{1+2\gamma}} + \cdots \right) 
.\end{align}
Note the expression of $L_{\overline{f}}(1-2\alpha,\sym^2f)$ with local factor at unramified primes is given by
\begin{align}
	\prod_{p \nmid M} \left( 1 + \frac{\overline{\lambda_f(p^2)}}{p^{1-2\alpha}} + \cdots \right) 
,\end{align}
and this expression allows us to handle the $\overline{\lambda_f(p^2)}/p^{1-2\alpha}$ term. We also note that $\chi_f(p) /p^{1 + 2\gamma}$ is handled identically to the previous case. To handle the $p^{1-\alpha+\gamma}$ term, note the local factor of $L(1+\alpha+\gamma,f\otimes \overline{f})^{-1}$ is
\begin{align}
	\left( 1- \frac{|\lambda_f(p)|^2}{p^{s}} + \ldots \right) 
.\end{align}
Upon substituting our approximations, we conclude with our desired formula.

\section{One-level density: Unscaled and Scaled}
\label{sec:one-level-density-unscaled-and-scaled}
\subsection{Averaging the logarithmic derivative}
To calculate the unscaled and scaled one-level density, we need the average of the logarithmic derivative of $L$-functions defined by
\begin{align}
\sum_{d \in \mathcal{D}_f^{+}(X)} \frac{L'}{L}(1 /2 + r, f_d) \ = \ \frac{\partial}{\partial \alpha} \bigg|_{\alpha=\gamma=r} R_f (\alpha,\gamma).
\end{align}
\begin{proposition}\label{prop:avg-logarithmic-derivative} Assume the Ratios Conjectures and that $1 /\log X \ll \mathrm{Re}(r) < 1 /4$ and $\mathrm{Im}(r) \ll X^{1- \varepsilon}$.
Then the average of the logarithmic derivative over the family $\mathcal{F}_f^{+}(X)$ is
\begin{align}
\sum_{d\in \mathcal{D}_f^{+}(X)} \frac{L'}{L} (1 /2 + r, f_d) & \ = \ \sum_{d\in \mathcal{D}_f^{+}(X)} \bigg( - \frac{L'}{L}(1+2r,\chi'_f) + \frac{L'}{L}(1+2r,\sym^2f) + A^{1}_f(r,r) \nonumber \\		
& \quad \ - \ \eta_f \bigg( \frac{\sqrt{M} |d|}{2\pi} \bigg)^{-2r} \frac{\Gamma\left( \frac{k}{2} - r\right)}{\Gamma\left( \frac{k}{2} + r\right)} \frac{L(1+2r,\chi'_f)L(1-2r,\sym^2\overline{f})}{L(1,\ad^2f)} \nonumber \\
& \qquad \ \times \ \widetilde{A}_f (-r,r) \bigg) + O(X^{1/2 + \varepsilon}),
\end{align}
where $\widetilde{A}_f$, $A^{1}_f(r,r)$ are given by \eqref{eqn:tilde-A_f(alpha,gamma)} and $
A^{1}_f(r,r)= \frac{\partial}{\partial \alpha} \bigg|_{\alpha=\gamma=r} A_f(\alpha,\gamma)$, respectively.
\end{proposition}
\begin{proof}
The proof is the same as that of Theorem 2.2 in \cite{HKS09}, mutatis mutandis. 
\end{proof}
\subsection{Unscaled one-level density}
With a formula for the average of logarithmic derivatives, we turn to finding the lower-order terms of the scaled one-level density functions for our family. For some Schwartz function $\varphi$ and ordinates $\gamma_d$ of the zeros on the critical line, the unscaled one-level density function is defined by
\begin{align}
D_1(\varphi,f) \ \coloneqq \ \sum_{d \in \mathcal{D}_f^{+}(X)} \sum_{\gamma_d} \varphi(\gamma_d)
.\end{align}
By the argument principle, we rewrite $D_1(\varphi;f)$ as
\begin{align}
\sum_{d \in \mathcal{D}_f^{+}(X)} \frac{1}{2\pi i} \left( \int_{(c)} - \int_{(1-c)} \right) \frac{L'}{L}(s,f_d) \varphi(i(1 /2-s ))\,ds
,\end{align}
where $1 /2 + 1 /\log X < c < 3 /4$ is fixed and $(c)$ denotes the path from $c-i\infty$ to $c+i\infty$.
Turning to the integral on the $(c)$-line, we have
\begin{align}\label{eqn:unscaled-integrand}
\frac{1}{2\pi} \int_{\mathbb{R}} \varphi(t-i(c-1 /2)) \sum_{d \in \mathcal{D}_f^{+}(X)} \frac{L'}{L} (c + it,f_d) \, dt
,\end{align}
and the sum over $d$ can be replaced by Proposition \ref{prop:avg-logarithmic-derivative}.
The bounds on the size $t$ coming from the Ratios Conjectures do not pose a problem; see \cite{HKS09,CS07}.

\begin{lemma}\label{lem:unscaled-regular-integrand}
	The integrand in~\eqref{eqn:unscaled-integrand} is regular at $t = 0$.
\end{lemma}
\begin{proof} Since we assumed the Generalized Riemann Hypothesis for Dirichlet $L$-functions, then $L(s,\chi_f)$ does not vanish on the line $\mathrm{Re}(s)=1$. In particular, the poles of $L(s,\chi_f')$ and $L(s,\sym^2f)$ at $s=1$ may obstruct regularity of the integrand in~\eqref{eqn:unscaled-integrand} at $t = 0$. The $L$-function $L(s,\chi_f)$ has a pole at $s = 1$ if and only if $\chi_f$ is principal, and $L(s,\sym^2f)$ has a pole only when $f$ is self-CM.

For the case $f$ has principal nebentype, $\widetilde{A}_f(-r,r) = A_f(-r,r)$ which implies $A_f$ is analytic.
As $r \to 0$, $\zeta(1+2r) = (2r)^{-1} + O(1)$ and $-(\zeta' /\zeta)(1+2r) = (2r)^{-1} + O(1)$.
After substituting the expansions into the expression in Proposition \ref{prop:avg-logarithmic-derivative}, the integrand in~\eqref{eqn:unscaled-integrand} is regular at $t=0$ if and only if $\widetilde{A}_f(0,0)=1$ which happens since $A_f(r,r)=1$.
		
For the case $f$ is self-CM, the symmetric square $L$-function of $f$ inherits a pole at $s=1$ from $L(s,f\otimes f)$, and $A_f$ is analytic. Using the relation $\lambda(p^{2m+1})\lambda(p) = \lambda(p^{2m+2}) + \lambda(p^{2m})$ for $p \nmid M$ and the multiplicativity of $\lambda(p)$ for $p = M$ (see Section \ref{sct:newforms-duals}), we get that $\widetilde{A}_f(-r,r) = A_f(-r,r)=1$. We use the constructed expectation $\eta_f=+1$ to obtain
\begin{align}
\lim_{r \to 0} \bigg[ & \frac{L'}{L} (1+2r,\sym^2f) - \left( \frac{\sqrt{M} |d|}{2\pi} \right) ^{-2r} \frac{\Gamma\left(\frac{k}{2}  - r \right)}{\Gamma\left(\frac{k}{2} + r \right)} \frac{L(1+2r,\chi'_f) L(1 -2r, \sym^2\overline{f})}{L(1,\ad^2f)} \widetilde{A}_f (-r,r) \bigg] \nonumber\\
& = \ \lim_{r \to 0} \left( -(2r)^{-1} + O(1) \right) - \frac{L(1+2r,\chi'_f) L(1-2r,\chi'_f)^{-1} L(1-2r,f\otimes f)}{\res(L(s,f\otimes f), 1)} \nonumber\\   
& = \ \lim_{r\to 0}\left( - (2r)^{-1} + O(1) \right) - \left( -(2r)^{-1} + O(1) \right) \  = \ O(1), 
\end{align}
which completes the proof.
\end{proof}	
\begin{proposition}\label{prop:unscaled-density}
Assume the Generalized Riemann Hypothesis and the Ratios Conjectures. The unscaled one-level density for the zeros of the family $\mathcal{F}_f^{+}(X)$ is
\begin{align}\label{eqn:unscaled-density}
D_1(\varphi,f) \ = \ \frac{1}{2\pi} & \int_{\mathbb{R}} \varphi(t) \sum_{d \in \mathcal{D}_f^{+}(X)} \bigg(\log \bigg( \frac{\sqrt{M} |d|}{2\pi} \bigg)^2 + \frac{\Gamma'}{\Gamma}\bigg( \frac{k}{2} + it \bigg) + \frac{\Gamma'}{\Gamma} \bigg( \frac{k}{2} - it \bigg) \\
& \qquad + \ \frac{L'}{L} \bigg( \frac{1}{2} + it, f_d \bigg) + \frac{L'}{L} \bigg( \frac{1}{2}+ it, \overline{f}_d \bigg) \bigg) \, dt + O(X^{1/2 + \varepsilon}) \nonumber
.\end{align}
\end{proposition}	
\begin{proof}
By Lemma~\ref{lem:unscaled-regular-integrand} we move the path of integration to $c = 1 /2$ and obtain 
\begin{align}
\frac{1}{2\pi} \int_{\mathbb{R}} \varphi(t) & \sum_{d \in \mathcal{D}_f^{+}(X)} \bigg(- \frac{L'}{L}(1 + 2it, \chi'_f) + \frac{L'_f}{L_f} (1+2it,\sym^2) + A_f^{1} (it,it) \nonumber \\
& - \ \eta_f \bigg( \frac{\sqrt{M} |d|}{2\pi} \bigg) ^{-2it} \frac{\Gamma\left( \frac{k}{2}-it \right)}{\Gamma\left( \frac{k}{2} +it \right)} \frac{L(1+2it,\chi_f)L_{\overline{f}}(1-2it,\sym^2)}{L_f(1,\ad^2)} \widetilde{A}_f(-it,it) \bigg) \, dt \nonumber \\
& \quad + \ O(X^{1 /2 + \varepsilon})
.\end{align}
The rest of the proof is the same as that of Theorem 2.3 in \cite{HKS09}, mutatis mutandis. 
\end{proof}

\subsection{Scaled one-level density} 
We calculate the one-level density for scaled zeros and recover the limit and next-to-leading-order term from the unscaled one-level density statistic. We rescale the zeros to have unit mean spacing by rescaling the variable $t$ by $\tau = tR/\pi$ in Equation \eqref{eqn:unscaled-density}. Motivated by \cite[Theorem 5.8]{IK04}, we equate the mean densities of eigenvalues with the mean densities of zeros by setting, as in \cite{DHKMS12},
\begin{align}
R \ \coloneqq \ 
\begin{dcases} \log \bigg( \frac{\sqrt{M} X}{2\pi} \bigg) - \frac{1}{2} & \qquad \chi_f \text{ principal}, \text{odd twists}, \\
\log \bigg( \frac{\sqrt{M} X}{2\pi} \bigg) & \qquad \chi_f \text{ principal, even twists or $f$ self-CM},\\
2 \log \bigg( \frac{\sqrt{M} X}{2\pi} \bigg) & \qquad \text{$f$ generic}.
\end{dcases}
\end{align}
We define the even test function $\varphi(t) \coloneqq g(\tau)$, and we define the scaled one-level density function to be $S_1(g,f) \coloneqq D_1(\varphi,f)/|\mathcal{D}_f^{+}(X)|$. Summation by parts and the cardinality estimate \eqref{eqn:cardinality-estimates} yield the following approximations (see Appendix \ref{appendix:counting_fundamental_discriminants}):
\begin{align}
\sum_{d\in \mathcal{D}_f^{+}(X)} \log \bigg( \frac{\sqrt{M} |d|}{2\pi} \bigg) \ & = \ |\mathcal{D}_f^{+}(X)| \bigg( \log\bigg( \frac{\sqrt{M} X}{2\pi }\bigg)-1 \bigg)   + O(X^{1 /2}), \\
\sum_{d\in \mathcal{D}_f^{+}(X)}\bigg( \frac{\sqrt{M} |d|}{2\pi} \bigg) ^{-2i\pi\tau /R} \ & = \ |\mathcal{D}_f^{+}(X)| \bigg( 1 -  \frac{2i \pi\tau}{R} \bigg)^{-1} e^{-2i\pi\tau} + O(X^{1 /2})
.\end{align}
Recall our aim is to obtain a series expansion of the scaled one-level density in terms of $R$. The normalized scaled one-level density is
\begin{align}\label{eq:principalnormalizeddensity}
&S_1(g,f) \ \coloneqq \ \frac{1}{|\mathcal{D}_f^{+}(X)|}\sum_{d\in\mathcal{D}_f^{+}(X)} \sum_{\gamma_d}g\left(\frac{\gamma_d R}\pi\right) \nonumber \\
  &=\frac{1}{2R}\int_{-\infty}^\infty g(\tau)\sum_{d\in\mathcal{D}_f^{+}(X)}\Bigg(2\log\left(\frac{\sqrt M|d|}{2\pi}\right)
    +\frac{\Gamma'}{\Gamma}\left(\frac{k}{2}+\frac{i\pi\tau}{R}\right)
    +\frac{\Gamma'}{\Gamma}\left(\frac{k}{2}-\frac{i\pi\tau}{R}\right) \nonumber \\
    &\qquad \qquad \qquad +\Bigg[-\frac{L'}
    {L}\left(1+\frac{2i\pi\,\tau}{R},\chi'_f\right)
    +\frac{L'}{L} \left(1+\frac{2i\pi\tau}{R},\sym^2f\right)
    +A_f^1\left(\frac{i\pi\tau}{R},\frac{i\pi\,\tau}{R}\right) \nonumber \\
    &\qquad \qquad \qquad-\left(\frac{\sqrt M |d|}{2\pi}\right)^{-2i\pi\tau/R}
    \frac{\Gamma\left(\frac k2-\frac{i\pi\tau}{R}\right)}
    {\Gamma\left(\frac{k}{2}+\frac{i\pi\tau}{R}\right)}
    \frac{L\left(1+\frac{2i\pi\tau}{R},\chi'_f\right)
      L\left(1-\frac{2i\pi\tau}{R},\sym^2\overline{f}\right)}{L(1,\ad^2f)} \nonumber \\
    &\qquad \qquad \qquad \times{A}_f\left(-\frac{i\pi\tau}{R},\frac{i\pi\tau}{R}\right)
    \Bigg] \nonumber\\
    &\qquad \qquad \qquad +\Bigg[-\frac{L'}
    {L}\left(1+\frac{2i\pi\,\tau}{R},\chi'_f\right)
    +\frac{L'}{L} \left(1+\frac{2i\pi\tau}{R},\sym^2f\right)
    +A_f^1\left(\frac{i\pi\tau}{R},\frac{i\pi\,\tau}{R}\right) \nonumber \\
    &\qquad \qquad \qquad-\left(\frac{\sqrt M |d|}{2\pi}\right)^{-2i\pi\tau/R}
    \frac{\Gamma\left(\frac k2-\frac{i\pi\tau}{R}\right)}
    {\Gamma\left(\frac{k}{2}+\frac{i\pi\tau}{R}\right)}
    \frac{L\left(1+\frac{2i\pi\tau}{R},\chi'_f\right)
      L\left(1-\frac{2i\pi\tau}{R},\sym^2\overline{f}\right)}{L(1,\ad^2f)} \nonumber \\
    &\qquad \qquad \qquad \times{A}_f\left(-\frac{i\pi\tau}{R},\frac{i\pi\tau}{R}\right)
    \Bigg]
    \Bigg)\, d t +O(X^{1/2+\epsilon}).
\end{align}

\begin{proposition}
Assume the Ratios Conjectures. 
Then the scaled one-level density for the zeros of the family $\mathcal{F}_f(X)$ is given by
\begin{align}
S_1(g,f) \ = \ \int_{\mathbb{R}} g(\tau) \left( 1 + Q(\tau) + O(R^{-3}) \right) d\tau,
\end{align}
where the lower order terms not in the error term are
\begin{align*}
Q(\tau)\ =\  \begin{dcases}
\frac{\sin(2\pi \tau)}{2\pi\tau} - a_1 \frac{1 + \cos(2\pi \tau)}{R} - a_2 \frac{\pi \tau \sin(2\pi \tau)}{R^{2}}& \quad \chi_f \ \text{\emph{principal}}, \text{\emph{ even twists}}, \\
- \frac{\sin(2\pi\tau)}{2i\pi\tau} - a_3 \frac{1-\cos(2\pi\tau)}{2R+1} + a_4 \frac{2\pi\tau\sin(2\pi\tau)}{(2R+1)^2} & \quad \chi_f \ \text{\emph{principal}}, \text{\emph{ odd twists}}, \\
- \frac{\sin(2\pi \tau)}{2\pi\tau}  + b_1 \frac{1 - \cos(2\pi \tau)}{R} + b_2 \frac{\pi \tau \sin(2\pi \tau)}{R^2} & \quad \chi_f \ \text{\emph{non-principal}},\ f = \overline{f}, \\ 
\frac{c_1 + c_2 \cos(2\pi \tau)}{R} + d_1 \frac{\pi \tau \sin(2\pi \tau)}{R^2}  & \quad \chi_f \ \text{\emph{non-principal}},\ f \neq \overline{f}
\end{dcases}
\end{align*}
with coefficients 
\begin{align}
a_1 \ & = \ 1 - \psi(k/2) - A_f^{1}(0,0) + \gamma - \frac{L'}{L}(1,\sym^2f) \label{eqn:coeff-a_1}, \\
a_2 \ & = \ -2\psi(k/2) - 2\psi(k/2)\gamma + 2\gamma - 2\gamma_1 + (2\psi(k/2)-2-2\gamma-B'(0)) \frac{L'}{L}(1,\sym^2f) \label{eqn:coeff-a_2} \nonumber \\
& \quad + \ (\gamma + 1 - \psi(k/2)) B'(0) + \frac{1}{4}B''(0) + 2\frac{L''}{L}(1,\sym^2f) ,\\
a_3 \ & = \ 2 - 2\psi(k/2) + 2\gamma_1 - 2 \frac{L'}{L}(1,\sym^2f) - 2 A_f^1(0,0), \label{eqn:coeff-a_3}\\
a_4 \ & = \ 4\psi(k /2) - 4i\pi\tau \gamma + 4\psi(k /2) \gamma + 4\gamma_1 +(2\psi(k /2) -2 - 2\gamma) B'(0) \label{eqn:coeff-a_4} \nonumber \\
& + ( 4 + 4\gamma + 2B'(0) - 4\psi(k /2))  \frac{L'}{L}(1,\sym^2f) - \frac{B''(0)}{2} - \frac{L''}{L}(1,\sym^2f), \\
b_1 \ & = \ 1-\psi(k/2) - \xi_0 \frac{L(1,\chi'_f)}{L(1,\ad^2f)} - A^{1}_f(0,0) + \frac{L'}{L}(1,\chi'_f) \label{eqn:coeff-b_1} ,\\
b_2 \ & = \ -2\psi(k/2) +B'(0)- \psi(k/2) B'(0)+ \frac{B''(0)}{4} + 2 \frac{L''}{L}(1,\chi'_f) \label{eqn:coeff-b_2}\\
&+\frac{L'}{L}(1,\chi'_f)\Big(-2\xi_0 +B'(0) +2-2\psi(\textstyle{\frac{k}{2}})\Big) + \frac{L(1,\chi'_f)}{L(1,\ad^2f)}\Big(2\psi(\textstyle{\frac{k}{2}}) \xi_0 -2\xi_0 +2\xi_1 -\xi_0 B'(0)\Big) ,\nonumber
\end{align}
\begin{align}
c_1 \ & = \ \psi(k/2) + \frac{1}{2} ( (A^{1}_f + A_{\overline{f}}^{1})(0,0)  - \frac{L'}{L}(1,\chi'_f) - \frac{L'}{L}(1,\chi'_{\overline{f}}) + \frac{L'}{L}(1,\sym^2f) + \frac{L'}{L}(1,\sym^2\overline{f}) ) \label{eqn:coeff-c_1} ,\\
c_2 \ & = \ -\frac{1}{2} ( \eta_f \widetilde{A}_f(0,0)L(1,\chi'_f) \frac{L (1,\sym^2\overline{f})}{L(1,\ad^2f)} + \eta_{\overline{f}} \widetilde{A}_{\overline{f}}(0,0)L(1,\chi'_{\overline{f}}) \frac{L(1,\sym^2f)}{L(1,\ad^2\overline{f})} )  \label{eqn:coeff-c_2},\\
d_1 \ & = \ \eta_f \frac{L(1,\sym^2\overline{f})}{L(1,\ad^2f)}\bigg( -\frac{1}{2}\widetilde{B}'_f(0)L(1,\chi'_f) + \psi(k /2) \widetilde{A}_f(0,0)L(1,\chi'_f) - \widetilde{A}_f(0,0)L(1,\chi'_f) \bigg) \label{eqn:coeff-d_1} \nonumber \\
& \quad + \ \eta_{\overline{f}} \frac{L(1,\sym^2f)}{L(1,\ad^2\overline{f})} \bigg( -\frac{1}{2}\widetilde{B}'_{\overline{f}}(0)L(1,\chi'_{\overline{f}}) + \psi(k /2) \widetilde{A}_{\overline{f}}(0,0)L(1,\chi'_{\overline{f}}) - \widetilde{A}_{\overline{f}}(0,0)L(1,\chi'_{\overline{f}}) \bigg) \nonumber\\
& \quad  + \ \eta_f \frac{L'(1,\sym^2\overline{f})}{L(1,\ad^2f)} \widetilde{A}_f(0,0) L(1,\chi'_f) + \eta_f \frac{L'(1,\sym^2f)}{L(1,\ad^2\overline{f})} \widetilde{A}_{\overline{f}}(0,0) L(1,\chi'_{\overline{f}}),
\end{align}
 where $\psi \coloneqq \Gamma'/{\Gamma}$ is the digamma function $B_f(s) = A_f(-s,s)$, $-B'(0)/2 = A^{1}(0,0)$, and $B_f^{(n)}(s) = \frac{d^n}{d r^n} \Big |_{r=s} A_f(-r,r)$.
\end{proposition}
\begin{remark}
We note that specializing \eqref{eqn:coeff-a_1} to $k=2$ gives a slightly different expression than Equation (3.18) of \cite{HKS09}. This is because they scale by $\log\left( \frac{\sqrt{M} X}{2\pi} \right) $ while we scale instead by $\log\left( \frac{\sqrt{M} X}{2\pi e} \right) $. Otherwise, the two match.
\end{remark}
\begin{proof}
Fix a cuspidal newform $f$ with principal nebentype. We note the simplifications $\widetilde{A}_f = A_f$ and $L(1+r,\chi'_f) = \zeta(1+r)$. For even twists, the expectation $\eta_f$ is equal to $+1$. For odd twists, $\eta_f$ equals $-1$. The proof is the same as that of (3.19) in \cite{HKS09} with appropriate substitutions of the following series expansions: 
\begin{align}
\frac{\Gamma'}{\Gamma}\left(\frac{k}{2}+\frac{i\pi\tau}{R}\right) & \ = \ 2\psi(k/2) + \frac{i\pi\tau}{R}\psi^{(1)}(k/2)+O\left(R^{-2}\right),\\
\frac{\Gamma'}{\Gamma}\left(\frac{k}{2}-\frac{i\pi\tau}{R}\right) & \ = \ 2\psi(k/2)- \frac{i\pi\tau}{R}\psi^{(1)}(k/2)+O\left(R^{-2}\right),\\
\frac{\Gamma\left(\frac{k}{2}+\frac{i\pi\tau}{R}\right)}{\Gamma\left(\frac{k}{2}-\frac{i\pi\tau}{R}\right) } & \ = \ 1-\frac{2i\pi\tau}{R}\psi(k/2) +O\left(R^{-2}\right),\\
A^1\left(\frac{i\pi\tau}{R},\frac{i\pi\tau}{R}\right) & \ = \ A^1(0,0)+\frac{i\pi\tau}{R}A^2(0,0)+O\left(R^{-2}\right),\\
A\left(-\frac{i\pi\tau}{R},\frac{i\pi\tau}{R}\right) & \ = \ A(0,0) + \frac{i\pi\tau}{R}A^1(0,0) - \frac{\pi^2\tau^2}{2R^2}A^2(0,0) +O\left(R^{-3}\right),\\
L\left(1-\frac{2i\pi\tau}{R},\sym^2f\right) & \ = \ -\frac{R}{2i\pi\tau} \frac{L(1,\ad^2f)}{L(1,\chi'_f)}+\xi_0-\frac{2i\pi\tau}{R}\xi_1 +O\left(R^{-2}\right),\\
\frac{L'}{L}\left(1+\frac{2i\pi\tau}{R},\sym^2f\right) & \ =
\frac{L'}{L}\left(1,\sym^2f\right)\\
& \quad + \frac{2i\pi\tau}{R}\frac{L\left(1,\sym^2f\right)L''\left(1,\sym^2f\right)-L'\left(1,\sym^2f\right)^2}{L\left(1,\sym^2f\right)^2}+O\left(R^{-2}\right)\nonumber
,\end{align}
where $\psi=\Gamma'/\Gamma$ is the digamma function.

Fix a cuspidal newform $f$ with self-CM. Since we assumed the sign is $+1$, then the expectation $\eta_f$ over $d$ equals 1.The symmetric square $L$-function $L(s,\sym^2f)$ has a simple pole at $s=1$ with residue
\begin{align}
r_{\sym^2} \ \coloneqq \ \res(L(s, \sym^2f), 1) \ =\ \frac{\res(L(s, f \otimes f), 1)}{L(1,\chi')} \ =\ \frac{L(1,\ad^2 f)}{L(1,\chi'_f)}.
\end{align}
The Laurent expansion at $s=0$ of $L(1+s,\sym^2f)$ is
\begin{align}
L(1+s,\sym^2f) \ = \ r_{\sym^2}s ^{-1} + \xi_0 + \xi_1 s + O(s^2),
\end{align}
where
\begin{align}
\xi_n \ = \ \oint_{\mathcal{C}} \frac{L(s,\sym^2f)}{s^{n+1}} \,ds
\end{align}
and the contour $\mathcal{C}$ encloses the point $s=1$.
The Laurent expansion at $s=0$ of the logarithmic derivative evaluated at $1+s$ is
\begin{align}
	\frac{L'}{L} \left( 1+s,\sym^2f \right) \ = \ - s ^{-1} + \frac{\xi_0}{r_{\sym^2}} + \left( \frac{2\xi_1}{r_{\sym^2}} - \frac{\xi_0^2}{r_{\sym^2}^2} \right)s + O(s^2) 
.\end{align}
Since $\chi_f$ is non-principal, $L(s,\chi_f')$ has no pole and we may expand the following:
\begin{align}
L(1+s,\chi'_f) & \ = \ L(1,\chi'_f) + L'(1,\chi'_f)s + L''(1,\chi'_f)s^2 ,\\
\frac{L'}{L}(1+s,\chi_f') & \ = \ \frac{L'}{L}(1,\chi') + \frac{L''}{L}(1,\chi'_f)s - \bigg( \frac{L'}{L}(1,\chi'_f)\bigg)^2 s
.\end{align}
We substitute the expansions in the rescaled version of \eqref{eqn:unscaled-density} and clear out odd terms to obtain our desired result.

Fix $f$ a cuspidal newform with non-principal nebentype that is not self-dual, i.e., it is generic. We return to \eqref{eqn:unscaled-density} for $\chi_f$ non-principal and $f \neq \overline{f}$. In this case, $L(s,\chi'_f)$ and $L(s,\sym^2f)$ are entire (see Section \ref{sct:sym-adj-sq} for details on the latter), and the values of $L(s,\chi'_f)$ and its derivatives at $s=1$ are not well known except in particular cases, such as when $\chi_f'$ is a quadratic character associated to a fundamental discriminant. In this case, $L(1,\chi'_f)$ is given by Dirichlet's class number formula. Note that $\widetilde{A}_f(0,0)$ can no longer be simplified in an obvious way.
\end{proof}

The one-level density of the cuspidal newforms converges to that of the one-level density of eigenvalues near 1 in certain compact groups. Following the definition proposed by \cite{DHKMS12}, we define the effective matrix size
\begin{align}
	N_{\eff} \ \coloneqq \ \begin{dcases}
		\frac{1}{2a_1}\log\left( \frac{\sqrt{M} X}{2\pi} \right), & \qquad \text{$\chi_f$ principal, even twists},\\
		\frac{1}{a_3}\left( \log\left( \frac{\sqrt{M} X}{2\pi} \right) -\frac{1}{2} \right)-\frac{1}{2}, & \qquad \text{$\chi_f$ principal, odd twists},\\
		\frac{1}{b_1}\log \left( \frac{\sqrt{M} X}{2\pi} \right), & \qquad \text{$f$ self-CM},
	\end{dcases}
\end{align}
where
\begin{align}
a_1 \ & = \ 1 - \psi(k/2) - A_f^{1}(0,0) + \gamma - \frac{L'}{L}(1,\sym^2f), \\
a_3 \ & = \ 2 - 2\psi(k/2) + 2\gamma_1 - 2 \frac{L'}{L}(1,\sym^2f) - 2 A_f^1(0,0),\\
b_1 \ & = \ 1-\psi(k/2) - \xi_0 \frac{L(1,\chi'_f)}{L(1,\ad^2f)} - A^{1}_f(0,0) + \frac{L'}{L}(1,\chi'_f),
\end{align}
 where $\psi \coloneqq \Gamma'/{\Gamma}$ is the digamma function $B_f(s) = A_f(-s,s)$, $-B'(0)/2 = A^{1}(0,0)$, and $B_f^{(n)}(s) = \frac{d^n}{d r^n} \Big |_{r=s} A_f(-r,r)$.

\section{Standard matrix size}\label{sct:standard-matrix-size}
The standard approach to choose the matrix size is to equate mean densities of eigenvalues with the mean density of zeros; following the work of \cite{DHKMS12}, this implies we should choose
\begin{align}
N_{\text{std}} \ = \ \log\bigg( \frac{\sqrt{M} |d|}{2\pi} \bigg) 
,\end{align}
where $M$ is the level of our cusp form $f$. The approach \cite{BBLM06} and \cite{DHKMS12} take is to multiply $N_{\text{std}}$ obtained in this way by a constant that originates from the arithmetic of the form so as to see agreement in lower-order terms of the one-level density or the pair-correlation statistic, respectively.

\section{Pair-correlation}\label{sct:pair_correlation}
The series expansion of the scaled one-level density for the unitary group has no lower-order terms. This means we cannot match the lower-order terms of the one-level density of a generic form we found above to that coming from the matrix side, and so we cannot refine the fit of the random matrix model for the generic case. Since we cannot extract any arithmetic in the one-level density of the generic case, we turn to pair-correlation. We chose this particular statistic through the discussions in \cite{DHKMS12,BBLM06}.

The pair-correlation statistic is insensitive to any finite set of zeros and, hence, is computed for one $L$-function only. We obtain a series expansion for the pair-correlation statistics of a single $L$-function in large $T$ to obtain the lower-order terms of arithmetic origin necessary for computing the effective matrix size. Then, under a mild restriction, we obtain the effective matrix size in our random matrix ensemble simply by averaging over all choices of quadratic twists ranging over our family. This fact simplifies the argument as we will not need to average over an infinite family of $L$-functions. As a byproduct, we also show that Montgomery's pair-correlation conjecture holds in the general setting of generic forms of level an odd prime under the Ratios Conjectures. The overall method, with necessary modifications, closely follows the work of \cite[Section 4]{CS07}, which computed the pair-correlation statistic for the Riemann zeta function using the Ratios Conjectures.

Let $\gamma,\gamma'$ denote the imaginary coordinates of non-trivial zeros of $L(s,f_d)$, and suppose $\varphi(s)$ is a holomorphic function throughout the strip $|\mathrm{Im}(s)| < 2$, real-valued on the real line, even, and satisfies the bound $\varphi(x) \ll 1 /(1+x^2)$ as $x \to \infty$. Since we assume the Ratios Conjectures, we are able to circumvent working directly with the Fourier transform of $\varphi$ as in \cite{RS96}, and thus we do not require any condition on the Fourier transform's support. We wish to evaluate the pair-correlation statistic given by
\begin{align}\label{eqn:pair_correlation_definition}
	P(f_d;\varphi) \ \coloneqq \ \sum_{0 < \gamma,\gamma' < T} \varphi(\gamma - \gamma')
.\end{align}
Roughly speaking, $P(f_d;\varphi)$ measures the spacing between pairs of zeros on the critical line. Throughout the argument, we assume GRH. To compute \eqref{eqn:pair_correlation_definition}, we apply another version of the Ratios Conjectures to
\begin{align}\label{eqn:avg-log-derivative-shifted-pair}
	\int_0^{T} \frac{L'}{L}(s+\alpha,f_d) \frac{L'}{L}(1-s+\beta,\overline{f}_d)\,dt,
\end{align}
where $s = 1/2 + it$.
We apply Conjecture 5.1 in \cite{CFZ08} to
\begin{align}
	\mathcal{T}_{f_d}(\alpha,\beta,\gamma,\delta) \ \coloneqq \ \int_{0}^{T} \frac{L(s+\alpha,f_d) L(1-s+\beta,\overline{f}_d)}{L(s+\gamma,f_d)L(1-s+\delta,\overline{f}_d)} \,dt,
\end{align}
since a generic cuspidal newform has unitary symmetry.
Hence, we substitute $K=L=1$ and the group $\Xi_{1,1} = \{(1),(12)\} $ which consists of the identity permutation and the transposition $(12)$, and we identify $\alpha_1 = \alpha$, $\alpha_2 = -\beta$, $\gamma_1=\gamma$, and $\delta_1 = \delta$ into the Ratios Conjectures.
We now state the Ratios Conjectures' lemma for our family.
\begin{conjecture}\label{conj:ratios-lemma}
For $-1 /4 < \mathrm{Re}(\alpha), \mathrm{Re}(\beta) <1 /4$, $1/\log(T) \ll \mathrm{Re}(\delta) < 1/4$, and $\mathrm{Im}(\alpha), \mathrm{Im}(\beta) \ll_{\varepsilon} T^{1-\varepsilon}$ for all $\varepsilon > 0$, we have
\begin{align}
\mathcal{T}_{f_d}(\alpha,\beta,\gamma,\delta) \ & = \ \int_0^{T} \bigg( Y_U(\alpha,\beta,\gamma,\delta) A_L(\alpha,\beta,\gamma,\delta) \nonumber \\
& \quad + \ \bigg( \frac{\sqrt{M} |d| t}{2\pi} \bigg)^{-2(\alpha+\beta)}  Y_U(-\beta,-\alpha,\gamma,\delta) A_L(-\beta,-\alpha,\gamma,\delta) \bigg) dt + O(T^{1 /2+\varepsilon}),
\end{align}
where
\begin{align}
Y_U(\alpha,\beta,\gamma,\delta) \ = \ \frac{L(1+\alpha+\beta,f_d\otimes \overline{f}_d)L(1+\gamma+\delta,f_d\otimes \overline{f}_d)}{L(1+\alpha+\delta,f_d\otimes \overline{f}_d)L(1+\beta+\gamma,f_d\otimes \overline{f}_d)}
,	\end{align}
and
\begin{align}
A_{L}(\alpha,\beta,\gamma,\delta) \ & = \ \prod_p \frac{(1-p^{-(1+\alpha+\beta)})(1-p^{-(1+\gamma+\delta)})}{(1-p^{-(1+\alpha+\delta)})(1-p^{-(1+\beta+\gamma)})} \nonumber \\
& \quad \times \ \sum_{m+h=n+k}\frac{\mu\lambda}{p^{(1 /2+\alpha)m+(1 /2+\beta)n+(1 /2+\gamma)h+(1 /2+\delta)k}}
.\end{align}
Here, $\mu = \mu_d$ is the coefficient on $n^{-s}$ of the reciprocal series $L(s, f_d)^{-1}$ for $\mathrm{Re}(s) > 1$; explicitly,
\begin{equation}
    \mu(n) \ = \ 
    \begin{cases}
        \lambda(n) & n = p ,\\
        \chi(n) & n = p^2 ,\\
        0 & n = p^j \text{ for } j > 2.
    \end{cases}
\end{equation}
\end{conjecture}
We also remark that
\begin{equation}
A_L(\alpha,\beta,\gamma,\delta) \ = \ Y_U(\alpha,\beta,\gamma,\delta)^{-1} W_{\mid}(\alpha,\beta,\gamma,\delta) W_{\nmid}(\alpha,\beta,\gamma,\delta),
\end{equation}
where
\begin{align}
W_\mid(\alpha,\beta,\gamma,\delta) \ & = \ \prod_{p\mid M}
    \Bigg[
      \sum_{m=0}^{\infty}\frac{\left|\lambda(p^m)\right|^2}{p^{(1+\alpha+\beta)m}}
      -\sum_{m=0}^{\infty}\frac{\lambda(p^{m+1})\overline\lambda(p^m)\overline\lambda(p)}
      {p^{(1+\alpha+\beta)m+(1+\alpha+\delta)}} \nonumber \\
      & \qquad -\sum_{m=0}^{\infty}\frac{\lambda(p^m)\overline\lambda(p^{m+1})\lambda(p)}
      {p^{(1+\alpha+\beta)m+(1+\beta+\gamma)}}
      +\sum_{m=0}^{\infty}\frac{\left|\lambda(p^m)\right|^2\left|\lambda(p)\right|^2}
      {p^{(1+\alpha+\beta)m+(1+\gamma+\delta)}}
      \Bigg],
    \end{align}
    and
    \begin{align}
  W_\nmid(\alpha,\beta,\gamma,\delta) \ & = \ \prod_{p\nmid N}
    \Bigg[
      \sum_{m=0}^{\infty}\frac{\left|\lambda(p^m)\right|^2}{p^{(1+\alpha+\beta)m}}-\sum_{m=0}^{\infty}\frac{\lambda(p^{m+1})\overline\lambda(p^m)\overline\lambda(p)}
      {p^{(1+\alpha+\beta)m+(1+\alpha+\delta)}}
      -\sum_{m=0}^{\infty}\frac{\lambda(p^m)\overline\lambda(p^{m+1})\lambda(p)}
      {p^{(1+\alpha+\beta)m+(1+\beta+\gamma)}} \nonumber \\
      &\qquad +\sum_{m=0}^{\infty}\frac{\overline\lambda(p^m)\lambda(p^{m+2})\overline\chi(p)}
      {p^{(1+\alpha+\beta)m+(2+2\alpha+2\delta)}}
      +\sum_{m=0}^{\infty}\frac{\lambda(p^m)\overline\lambda(p^{m+2})\chi(p)}
      {p^{(1+\alpha+\beta)m+(2+2\beta+2\gamma)}} \nonumber \\
      &\qquad +\sum_{m=0}^{\infty}\frac{\left|\lambda(p^m)\right|^2\left|\lambda(p)\right|^2}
      {p^{(1+\alpha+\beta)m+(1+\gamma+\delta)}}
      -\sum_{m=0}^{\infty}\frac{\lambda(p^{m+1})\overline\lambda(p^m)\lambda(p)\overline\chi(p)}
      {p^{(1+\alpha+\beta)m+(2+\alpha+\gamma+2\delta)}} \nonumber \\
      &\qquad -\sum_{m=0}^{\infty}\frac{\lambda(p^m)\overline\lambda(p^{m+1})\overline\lambda(p)\chi(p)}
      {p^{(1+\alpha+\beta)m+(2+\beta+2\gamma+\delta)}}
      +\sum_{m=0}^{\infty}\frac{\left|\lambda(p^m)\right|^2\left|\chi(p)\right|^2}
      {p^{(1+\alpha+\beta)m+(2+2\gamma+2\delta)}}
      \Bigg].
\end{align}
We used the definition of $\mu$ to state that, if $p \mid N = M |d|^2$, we are free to discard all terms except $h,k \in \{0,1\}$. Otherwise, we may discard all terms except $h,k \in \{0,1,2\}$.

\subsection{Averaging the logarithmic derivative} 
We obtain the formula for \eqref{eqn:avg-log-derivative-shifted-pair} by differentiating the result of Conjecture \ref{conj:ratios-lemma} which allows us to compute $P(f_d; \varphi)$ using contour integration. By expanding in series the formula for $P(f_d; \varphi)$ in large $T$, we obtain the desired lower-order terms.
\begin{proposition}\label{prop:int-unramified}
	Assume Conjecture \ref{conj:ratios-lemma}, and let $\alpha$, $\beta$, $\gamma$, and $\delta$ be as above. Then
\begin{align}
	& \int_0^{T} \frac{L'}{L}(s+\alpha,f_d) \frac{L'}{L}(1-s+\beta,\overline{f}_d) \, dt \ = \ \int_0^{T} \bigg( \bigg( \frac{L'_{\star}}{L_{\star}} \bigg)' (1+\alpha+\beta,f_d\otimes \overline{f}_d) \nonumber \\
	& \qquad + \ \frac{1}{c_{f_d}^2}\bigg( \frac{\sqrt{M} |d|t}{2\pi} \bigg) ^{-2(\alpha+\beta)} L(1+\alpha+\beta,f_d\otimes \overline{f}_d)L(1-\alpha-\beta,f_d\otimes \overline{f}_d) \nonumber \\
	& \qquad \times \ A_L(-\beta,-\alpha,\alpha,\beta) + \mathscr{C}(1+\alpha+\beta) \bigg)\,dt + O(T^{1 /2+\varepsilon})
,\end{align}
where $N \coloneqq M |d|^2$, $c_{f_d} = \emph{Res}(L(s, f_d \otimes \overline{f}_d), 1)$, $\mathscr{C}(1+\alpha+\beta)$ is given in \eqref{script_C}, and $L_{\star}(s, f_d \otimes \overline{f}_d)$ is the unramified part of the Rankin-Selberg convolution as defined in \eqref{rankin_selberg}.
\end{proposition}
\begin{proof}
For notation, put 
\begin{align}
\Lambda \ & \coloneqq \ Y_U(\alpha,\beta,\gamma,\delta)
A_L(\alpha,\beta,\gamma,\delta), \\
\Omega \ & \coloneqq \ \bigg( \frac{\sqrt{M}|d|t}{2\pi} \bigg)^{-2(\alpha+\beta)} Y_U(-\beta,-\alpha,\gamma,\delta) A_L(-\beta,-\alpha,\gamma,\delta).
\end{align}
We differentiate the integrand in Conjecture \ref{conj:ratios-lemma} with respect to $\alpha$ and $\beta$ and then map $\gamma \mapsto \alpha$ and $\delta \mapsto \beta$ which yields
\begin{align}
\int_{0}^{T} \frac{\partial^2\Lambda}{\partial\beta\,\partial\alpha}\bigg|_{(\gamma,\delta)\ =\ (\alpha,\beta)} + \frac{\partial^2\Omega}{\partial\beta\,\partial\alpha}\bigg|_{(\gamma,\delta)\ =\ (\alpha,\beta)}\,dt + O(T^{1/2+\varepsilon}),
\end{align}
where the third equality comes from Conjecture \ref{conj:ratios-lemma}.
We first compute the derivative of $\Lambda$ evaluated at $(\gamma,\delta) = (\alpha,\beta)$, noting that $A_L(\alpha,\beta,\alpha,\beta)=1$: 
\begin{align}\label{eqn:lambda_derivative_step_1}
    \frac{\partial^2\Lambda}{\partial\beta\,\partial\alpha} \bigg|_{(\gamma,\delta) \ = \ (\alpha,\beta)} \ = \ & \frac{L''(1+\alpha+\beta, f_d \otimes \overline{f}_d)L(1+\alpha+\beta, f_d \otimes \overline{f}_d) - L'(1+\alpha+\beta, f_d \otimes \overline{f}_d)^2}{L(1+\alpha+\beta, f_d \otimes \overline{f}_d)^2} \nonumber \\
    & \quad + \ \frac{\partial^2}{\partial\beta\,\partial\alpha}\bigg|_{(\gamma,\delta) \ = \ (\alpha,\beta)} A_L(\alpha,\beta,\gamma,\delta).
\end{align}
Since the derivative of $A_L(\alpha,\beta,\gamma,\delta)$ is difficult to compute, we discuss the main ideas at the end of the proof. 

We turn to the derivative of $\Omega$. By applying the usual product rule and setting $\gamma \mapsto \alpha$ and $\beta \mapsto \delta$ by taking the limits appropriately, we are left with
\begin{align}
\frac{\partial^2\Omega}{\partial\beta\,\partial\alpha} \bigg|_{(\gamma,\delta) \ = \ (\alpha,\beta)} \ = \ &\bigg( \frac{\sqrt{M}|d| t}{2\pi} \bigg)^{-2(\alpha+\beta)} \bigg(\frac{L'}{L^2} \bigg)^2(1, f_d\otimes \overline{f}_d) \\
&\times \ L(1-\alpha-\beta, f_d\otimes \overline{f}_d)L(1+\alpha+\beta, f_d\otimes \overline{f}_d)A_L(-\beta,\alpha,\alpha,\beta),\nonumber
\end{align}
upon substituting $(\alpha, \beta) = (\gamma,\delta)$. 

We now wish to evaluate $(L'/L^2)^2(1, f_d \otimes \overline{f}_d)$. Since $L(1, f_d \otimes \overline{f}_d)$ has a simple pole at 1 and $L'(1,f_d\otimes \overline{f}_d)$ has a double pole at 1, the logarithmic derivative $(L'/L)(s, f_d \otimes \overline{f}_d)$ has a simple pole at 1. 
To compute its residue, take a small circular contour $\mathcal{C}$ oriented counterclockwise with center at $s = 1$. By the residue theorem and the argument principle, we have
\begin{equation}
    \text{Res}\bigg( \frac{L'}{L}(s, f_d \otimes \overline{f}_d), 1 \bigg) \ = \ \frac{1}{2\pi i} \oint_{\mathcal{C}} \frac{L'}{L}(s, f_d \otimes \overline{f}_d)\,ds \ = \ \#\{\text{zeros of } L\} - \#\{\text{poles of } L\}.
\end{equation}
The contour may be chosen sufficiently small so that $L(1, f_d \otimes \overline{f}_d)$ has no zeros within the contour which gives the residue $\res( (L'/L)(s, f_d \otimes \overline{f}_d), 1 ) =  -1$. 
The residue of $L(s, f_d \otimes \overline{f}_d)$ at $s = 1$ is $c_{f_d}$ (see equation \eqref{eqn:residue-at-1-f-bar-f}), and so $(L'/L^2)(s, f_d \otimes \overline{f}_d)$ is entire with value $-1/c_{f_d}$ at $s=1$.
We get
\begin{align}
    \frac{\partial^2\Omega}{\partial\beta\,\partial\alpha} \bigg|_{(\gamma,\delta) \ = \ (\alpha,\beta)} \ &= \ \frac{1}{c_{f_d}^2} \bigg( \frac{\sqrt{M}|d| t}{2\pi} \bigg)^{-2(\alpha+\beta)} L(1-\alpha-\beta, f_d \otimes \overline{f}_d) \nonumber\\
    & \quad \times \ L(1+\alpha+\beta, f_d \otimes \overline{f}_d)A_L(-\beta,\alpha,\alpha,\beta).
\end{align}
It remains to compute
\begin{equation}
    \frac{\partial^2}{\partial\beta\,\partial\alpha}\bigg|_{(\gamma,\delta) \ =  \ (\alpha,\beta)} A_L(\alpha,\beta,\gamma,\delta).
\end{equation}
After applying the product rule and the fact that $Y_U(\alpha,\beta,\alpha,\beta) = 1$, we obtain the expression
\begin{align}
    & \frac{\partial^2}{\partial\beta\,\partial\alpha}\bigg|_{(\gamma,\delta) \ = \ (\alpha,\beta)} A_L (\alpha,\beta,\gamma,\delta) \\
    & \quad = \ \frac{\partial^2}{\partial\beta\,\partial\alpha}\bigg|_{(\gamma,\delta) \ = \ (\alpha,\beta)} \bigg( Y_U(\alpha,\beta,\gamma,\delta)^{-1} + W_{\mid}(\alpha,\beta,\gamma,\delta) + W_{\nmid}(\alpha,\beta,\gamma,\delta) \bigg) .\nonumber
\end{align}
By repeated application of the product and quotient rule, we obtain
\begin{equation}\label{eqn:partial-Y-expanded}
    \frac{\partial^2}{\partial\beta\,\partial\alpha}\bigg|_{(\gamma,\delta) \ = \ (\alpha,\beta)} Y_U(\alpha,\beta,\gamma,\delta)^{-1} \ = \ -\left( \frac{L'}{L} \right)'(1+\alpha+\beta,f_d \otimes \overline{f}_d).
\end{equation}
We focus our attention to the derivative of $W_{\mid}(\alpha,\beta,\gamma,\delta)$, and fix a prime $p$ dividing $M$ (such a prime, by assumption, must be $M$).
Noting the Fourier coefficients $\lambda_{f_d}(p)$ are completely multiplicative at primes dividing $M$, we use the formula for the sum of a geometric series to obtain
\begin{equation}\label{eqn:second-partial-W-mid-log-sq}
    \frac{\partial^2}{\partial\beta\,\partial\alpha}\bigg|_{(\gamma,\delta) \ =  \ (\alpha,\beta)} W_{\mid}(\alpha,\beta,\gamma,\delta) \ = \ \frac{\log(M)^2}{ |\lambda_{f_d}(M)|^{-2} \cdot M^{1+\alpha+\beta} - 1}.
\end{equation}
Finally, it remains to calculate the derivative of $W_{\nmid}(\alpha,\beta,\gamma,\delta)$. The computations are considerably more involved as there is no guarantee that the Fourier coefficients $\lambda_{f_d}(p)$ are multiplicative. The strategy involves carefully re-indexing the sums appearing in $W_{\nmid}(\alpha,\beta,\gamma,\delta)$ and repeatedly applying the relation \eqref{eqn:lambda-coefficient-relation}.
We then expand $\lambda_{f_d}(p^{m+2}) - \chi(p) \lambda_{f_d}(p^m)$ in terms of the Satake parameters attached to $f_d$ (see Section \ref{sct:newforms-duals}).
Letting $L_{\star}$ denote the unramified part of $L$, given by the Euler product
\begin{equation}\label{rankin_selberg}
L_{\star}(s) \ = \ \prod_{p\nmid N} L_p(s) \ = \ L(s) \prod_{p \mid N} L_p(s)^{-1}
,\end{equation}
in the half-plane of convergence and its analytic continuation elsewhere, we obtain that
\begin{align}\label{script_B}
& \frac{\partial^2}{\partial\beta\,\partial\alpha}\bigg|_{(\gamma,\delta) \ =  \ (\alpha,\beta)} W_{\nmid}(\alpha,\beta,\gamma,\delta) \ = \ \left( \frac{L_{\star}'}{L_{\star}} \right)'(1+\alpha+\beta, f_d \otimes \overline{f}_d) \nonumber \\
& \quad - \ \sum_{p \nmid N} \bigg[ \bigg( \frac{\alpha_{f_d}(p) \overline{\alpha}_{f_d}(p) \log(p)}{p^{1+\alpha+\beta} - \alpha_{f_d}(p) \overline{\alpha}_{f_d}(p)} \bigg)^2
+ \bigg( \frac{\alpha_{f_d}(p) \overline{\beta}_{f_d}(p) \log(p)}{p^{1+\alpha+\beta} - \alpha_{f_d}(p) \overline{\beta}_{f_d}(p)} \bigg)^2 \nonumber \\
& \quad + \ \bigg( \frac{\overline{\alpha}_{f_d}(p) \beta_{f_d}(p) \log(p)}{p^{1+\alpha+\beta} - \overline{\alpha}_{f_d}(p) \beta_{f_d}(p)} \bigg)^2 + \bigg( \frac{\beta_{f_d}(p) \overline{\beta}_{f_d}(p) \log(p)}{p^{1+\alpha+\beta} - \beta_{f_d}(p) \overline{\beta}_{f_d}(p)} \bigg)^2
    \bigg].
\end{align}
Let $\mathscr{B}(1+\alpha+\beta)$ denote the sum over $p\nmid N$ on the RHS, which allows us to write
\begin{align}\label{eqn:second-partial-W-nmid-unramified-B}
\frac{\partial^2}{\partial\beta\,\partial\alpha}\bigg|_{(\gamma,\delta) \ =  \ (\alpha,\beta)} W_{\nmid}(\alpha,\beta,\gamma,\delta) \ = \ &\left( \frac{L_{\star}'}{L_{\star}} \right)'(1+\alpha+\beta, f_d \otimes \overline{f}_d) - \mathscr{B}(1+\alpha+\beta).
\end{align}
Note that the Ramanujan-Petersson conjecture is known for holomorphic cusp forms of even weight. In particular, we get that $|\alpha_{f_d}| = |\beta_{f_d}| = 1$ for primes $p$ not dividing the level $N$; see \cite{IK04} and \cite{Sar05}. Therefore, $\mathscr{B}(v)$ exists and is analytic in a neighborhood of $v = 1$. Finally, put
\begin{equation}\label{script_C}
    \mathscr{C}(1 + \alpha + \beta) \ \coloneqq \ -\mathscr{B}(1 + \alpha + \beta) + \frac{\log(M)^2}{ |\lambda_{f_d}(M)|^{-2} \cdot M^{1+\alpha+\beta} - 1}.
\end{equation}
Putting everything together yields the result.
\end{proof}

\begin{remark}
The Ramanujan-Petersson conjecture also states that $\alpha_{f_d}(p) = \overline{\beta}_{f_d}(p)$; see \cite[Section 3]{ILS00}. Therefore, $\mathscr{B}(v)$ simplifies to
\begin{align}
    \mathscr{B}(v) \ = \ \sum_{p \nmid N} \bigg[ 2 \cdot \bigg( \frac{\log(p)}{p^{1+v} - 1} \bigg)^2
    + \bigg( \frac{\alpha_{f_d}(p)^2 \log(p)}{p^{1+v} - \alpha_{f_d}(p)^2} \bigg)^2 + \bigg( \frac{ \beta_{f_d}(p)^2 \log(p)}{p^{1+v} - \beta_{f_d}(p)^2} \bigg)^2
    \bigg].
\end{align}
\end{remark}

\subsection{Contour integration to compute pair-correlation}
The formula for the average of the logarithmic derivative for shifted $L$-functions allows us to obtain a formula for $P(f_d; \varphi)$, which will later be used to perform a series expansion to obtain lower-order terms of arithmetic origin with which to calibrate our effective matrix size. Much of the work in this section is analogous to that presented in \cite{CS07}.
\begin{proposition}\label{contour_theorem}
Set
\begin{align}
    I_r \ \coloneqq \ &\int_{-T}^{T} \varphi(r) \bigg( 2\log^2\bigg( \frac{\sqrt{M} |d| t}{2\pi} \bigg) + \bigg( \frac{L'_{\star}}{L_{\star}} \bigg)' (1 + ir, f_d \otimes \overline{f}_d) \\
    &+ \ \frac{1}{c_{f_d}^2} \bigg( \frac{\sqrt{M} |d| t}{2\pi} \bigg)^{-2ir} L(1+ir, f_d \otimes \overline{f}_d)L(1-ir, f_d \otimes \overline{f}_d) \mathscr{A}(ir) + \mathscr{C}(1+ir) \bigg)\,dr \nonumber.
\end{align}
Assuming the Ratios Conjectures and with $\varphi$ as above, we have
\begin{equation}
    P(f_d; \varphi) \ = \ \sum_{0 < \gamma,\gamma' < T} \varphi(\gamma-\gamma') \ = \ \frac{1}{2\pi^2} \int_{0}^{T} \bigg[ 2\pi\varphi(0) \log\bigg( \frac{\sqrt{M}|d| t}{2\pi} \bigg) + I_r \bigg]\,dt + O(T^{1/2 + \varepsilon}).
\end{equation}
Here, $I_r$ should be regarded as a principal-value integral near $r=0$.
\end{proposition}
The proof mainly relies on contour integration and some asymptotic analysis; the formula for \eqref{eqn:avg-log-derivative-shifted-pair} will help with the integrals. Choose real $a$ and $b$ such that
\begin{equation}
    \frac{1}{2} + \frac{1}{\log T} \ <\ a \ <\ b \ <\ \frac{3}{4}.
\end{equation}
Then, let $\mathcal{C}_1$ be the contour with vertices at $a$, $a + iT$, $1-a+iT$, and $1-a$, oriented counterclockwise. Let $\mathcal{C}_2$ be the contour with vertices at $b$, $b+iT$, $1-b+iT$, and $1-b$, also oriented counterclockwise.

By GRH, the poles of $(L'/L)(s)$ within the contours occur at the zeros $z = 1/2 + i\gamma$ and $w = 1/2 + i\gamma'$ of the $L$-function; moreover, all of these poles are simple. Therefore, by an adaptation of the residue theorem applied to both contours, we obtain 
\begin{equation}
    P(f_d; \varphi) \ = \ \frac{1}{(2\pi i)^2} \int_{\mathcal{C}_1} \int_{\mathcal{C}_2} \frac{L'}{L}(z, f_d) \, \frac{L'}{L}(w, f_d) \, \varphi(-i(z-w))\,dw\,dz.
\end{equation}
As in \cite{CS07}, a standard convexity estimate for $\GL(2)$ allows us to absorb the integrals along the horizontal segments into the error term of \eqref{eqn:avg-log-derivative-shifted-pair}.
We end up with four double integrals $I_1, \hdots, I_4$. In particular,
\begin{enumerate}[label=(\roman*)]
    \item $I_1$ is the (double) integral whose vertical segments have real parts $a$ and $b$;
    \item $I_2$ is the (double) integral whose vertical segments have real parts $1-a$ and $1-b$;
    \item $I_3$ is the (double) integral whose vertical segments have real parts $a$ and $1-b$;
    \item $I_4$ is the (double) integral whose vertical segments have real parts $1-a$ and $b$.
\end{enumerate}

\begin{lemma}
The integral $I_1$ is equal to $O(T^{\varepsilon+1/2})$.
\end{lemma}
\begin{proof}
We are free to use GRH to move the contours to the right of 1 without picking up additional residues. Integrating term-by-term, we obtain $I_1 \ll T^{\varepsilon} \ll T^{\varepsilon+1/2}$.
\end{proof}

\begin{lemma} The integral $I_2$ is given by
\begin{equation}
    I_2 \ = \ \frac{4}{(2\pi)^2} \int_{-T}^{T} \varphi(\rho) \int_{0}^{T} \log\bigg(\frac{\sqrt{M} |d|v}{2\pi}\bigg)^2 \,dv\,d\rho + O(T^{\varepsilon}).
\end{equation}
\end{lemma}
\begin{proof}
Working with $I_2$, we parametrize both paths and use the functional equation
\begin{equation}
    \frac{L'}{L}(1/2+ix, f_d) \ = \ \frac{\Phi_d'}{\Phi_d}(1/2+ix, f_d) - \frac{L'}{L}(1/2+ix, \overline{f}_d).
\end{equation}
Parametrization causes an $i \cdot i = i^2$ to fall to the front, cancelling with the copy in the denominator. 
This yields
\begin{equation}
    I_2 \ = \ \frac{1}{(2\pi)^2} \int_{0}^{T} \int_{0}^{T} \frac{\Phi_d'}{\Phi_d}\left(\frac{1}{2}+iu, f_d\right) \frac{\Phi_d'}{\Phi_d}\left(\frac{1}{2}+iv, f_d\right) \varphi(u-v)\,du\,dv + O(T^{\varepsilon}),
\end{equation}
where the error term comes from the second term in the functional equation. 
The integrand is symmetric in $u$ and $v$ since $\varphi$ is even, and so we may write
\begin{equation}
    I_2 \ = \ \frac{2}{(2\pi)^2} \int_{0}^{T} \int_{v}^{T} \frac{\Phi_d'}{\Phi_d}\left(\frac{1}{2}+iu, f_d\right) \frac{\Phi_d'}{\Phi_d}\left(\frac{1}{2}+iv, f_d\right) \varphi(u-v)\,du\,dv + O(T^{\varepsilon}).
\end{equation}
As in \cite{CS07}, we use the asymptotic estimate
\begin{equation}
    \frac{\Phi_d'}{\Phi_d}(1 /2+ix, f_d) \ = \ -2\log\bigg(\frac{\sqrt{M} |d|}{2\pi}x\bigg)\big( 1 + O(x^{-1}) \big)
    ,\end{equation}
and perform the change of variable $u \mapsto v+\rho$ while pulling out $\varphi(\rho)$ to get
\begin{equation}
    I_2 \ = \ \frac{8}{(2\pi)^2} \int_{0}^{T} \varphi(\rho) \int_{0}^{T-\rho} \log\bigg(\frac{\sqrt{M} |d|}{2\pi}(v+\rho)\bigg) \log\bigg(\frac{\sqrt{M}|d|}{2\pi}v\bigg) \,dv\,d\rho + O(T^{\varepsilon}).
\end{equation}
Once more, we perform the change of variable $v \mapsto vT$ to get
\begin{equation}
    I_2 \ = \ \frac{8T}{(2\pi)^2} \int_{0}^{T} \varphi(\rho) \int_{0}^{1-\rho/T} \log\bigg(\frac{\sqrt{M} |d|}{2\pi}(vT+\rho)\bigg) \log\bigg(\frac{\sqrt{M} |d|}{2\pi}vT\bigg) \,dv\,d\rho + O(T^{\varepsilon}).
\end{equation}
Since $\varphi(x) \ll 1/(1+x^2)$ for real $x$, we may extend the upper limit of integration; the error that results is of order $\log(T)^3$, which gets absorbed into the existing error term.
Due to the same asymptotic on the test function, we may also replace $\log(\sqrt{M}|d|(vT+\rho)/(2\pi))$ with $\log(\sqrt{M}|d|vT/(2\pi))$, also with an error of $\log(T)^3$. This error gets absorbed, resulting in
\begin{equation}
    I_2 \ = \ \frac{8T}{(2\pi)^2} \int_{0}^{T} \varphi(\rho) \int_{0}^{1} \log\bigg(\frac{\sqrt{M} |d|vT}{2\pi}\bigg)^2 \,dv\,d\rho + O(T^{\varepsilon}).
\end{equation}
Finally, the evenness of $\varphi$ and the change of variable $vT \mapsto T$ yields our desired form.
\end{proof}

\begin{lemma}
The integral $I_3$ is given by
\begin{equation}
    I_3 \ = \ \frac{1}{(2\pi)^2 i} \int_{0}^{T} \int_{-\delta-iT}^{-\delta+iT} \varphi(-i\rho) \mathscr{I}(-\rho,t)\,d\rho\,dt + O(T^{1/2+\varepsilon}),
\end{equation}
where $\mathscr{I}(\rho,t)$ is defined in \eqref{script_I}.
\end{lemma}
\begin{proof}
We consider $I_3$ and parametrize to get
\begin{equation}
    I_3 \ = \ \frac{1}{(2\pi i)^2} \int_{1-b+iT}^{1-b} \int_{a}^{a+iT} \frac{L'}{L}(w, f_d) \, \frac{L'}{L}(z, f_d) \, \varphi(-i(z-w))\,dw\,dz,
\end{equation}
where here the paths of integration are line segments as indicated by the bounds. Note that the integral along the segment with real part $1-b$ is traversed from top to bottom due to the counterclockwise orientation of $\mathcal{C}_2$. Switching the bounds on the outer integral yields
\begin{equation}
    I_3 \ = \ -\frac{1}{(2\pi i)^2} \int_{1-b}^{1-b+iT} \int_{a}^{a+iT} \frac{L'}{L}(w, f_d) \, \frac{L'}{L}(z, f_d) \, \varphi(-i(z-w))\,dw\,dz.
\end{equation}
We perform the change of variable $z = w + \rho$ and obtain
\begin{equation}
    I_3 \ = \ -\frac{1}{(2\pi)^2 i} \int_{1-a-b-iT}^{1-a-b+iT} \varphi(-i\rho) \int_{T_1}^{T_2} \frac{L'}{L}(a+it, f_d) \frac{L'}{L}(a+it+\rho)\,dt\,d\rho,
\end{equation}
where
\begin{equation}
    T_1 \ = \ \max\{0, -\text{Im}(\rho)\} \quad \text{and} \quad T_2 \ = \ \min\{T, T - \text{Im}(\rho)\}.
\end{equation}
We again use the functional equation, this time in the form
\begin{equation}
    \frac{L'}{L}(a+it+\rho, f_d) \ = \ \frac{\Phi_d'}{\Phi_d}(a+it+\rho) - \frac{L'}{L}(1-a-it-\rho, \overline{f}_d).
\end{equation}
Note that the term with $\Phi_d'/\Phi_d$ is small, which may be seen by moving the contour to the right. Letting $s = 1/2 + it$, which runs along the critical line, we get
\begin{align}
	I_3 & \ = \ \frac{1}{(2\pi)^2 i} \int_{1-a-b-iT}^{1-a-b+iT} \varphi(-i\rho) \\
	    & \qquad \qquad \qquad \times  \int_{T_1}^{T_2} \frac{L'}{L}\bigg(s + \left(a - \frac{1}{2}\right), f_d\bigg) \frac{L'}{L}\bigg(1-s+\left(\frac{1}{2}-a-\rho\right), \overline{f}_d\bigg)dt d\rho + O(T^{\varepsilon}) \nonumber 
.\end{align}
The point of rewriting the integrand in this fashion is to obtain an average of two shifted $L$-functions. Using the formula for \eqref{eqn:avg-log-derivative-shifted-pair}, the inner integral becomes
\begin{align}
    &\int_{T_1}^{T_2} \bigg(\frac{L'_{\star}}{L_{\star}}\bigg)'(1-\rho, f_d \otimes \overline{f}_d) + \frac{1}{c^2_{f_d}} \bigg( \frac{\sqrt{M}|d| t}{2\pi} \bigg)^{2\rho} \nonumber\\
    & \quad \times \ L(1-\rho, f_d \otimes \overline{f}_d) L(1+\rho, f_d \otimes \overline{f}_d) A_L\left(-\frac{1}{2}+a+\rho, \frac{1}{2}-a, a-\frac{1}{2}, \frac{1}{2}-a-\rho\right) \nonumber \\
    & \quad - \ \mathscr{B}(1-\rho) + \frac{\log(M)^2}{|\lambda_{f_d}(M)|^{-2} \cdot M^{1-\rho} - 1}\,dt + O(T^{1/2 + \varepsilon}). 
\end{align}
Then, the definition of $\mathscr{C}(1-\rho)$ yields
\begin{align}
    &\int_{T_1}^{T_2} \bigg(\frac{L'_{\star}}{L_{\star}}\bigg)'(1-\rho, f_d \otimes \overline{f}_d) + \frac{1}{c^2_{f_d}} \bigg( \frac{\sqrt{M}|d| t}{2\pi} \bigg)^{2\rho} \nonumber \\
    & \quad \times \ L(1-\rho, f_d \otimes \overline{f}_d) L(1+\rho, f_d \otimes \overline{f}_d) A_L\left(-\frac{1}{2}+a+\rho, \frac{1}{2}-a, a-\frac{1}{2}, \frac{1}{2}-a-\rho\right) \nonumber \\
    & \quad + \ \mathscr{C}(1-\rho) \,dt + O(T^{1/2 + \varepsilon}) .
\end{align}
For the sake of notation, put $\delta\coloneqq a+b-1$, and let $\mathscr{I}(-\rho,t)$ denote the inner integral from $T_1$ to $T_2$. Moreover, the integral over $t$ may be extended to the interval $[0,T]$, rather than merely $[T_1, T_2]$, with an error of at most
\begin{equation}
    T^{\varepsilon} \int_{\rho} |\rho| \cdot |f(\rho)| \,d\rho \, \ll \, T^{\varepsilon}.
\end{equation}
Thus, after several substitutions, we obtain our desired expression.
\end{proof}

\begin{lemma}
The integral $I_4$ is given by
\begin{equation}
    I_4 \ = \ \frac{1}{(2\pi)^2 i} \int_{0}^{T} \int_{\delta-iT}^{\delta+iT} \varphi(-i\rho) \mathscr{I}(\rho,t)\,d\rho\,dt + O(T^{1/2+\varepsilon}),
\end{equation}
where $\mathscr{I}(\rho,t)$ is defined in \eqref{script_I}.
\end{lemma}
\begin{proof}
Recall that $I_4$ has vertical paths with real parts $1-a$ and $b$. The computations are similar to those for $I_3$, and we again appeal to the formula for \eqref{eqn:avg-log-derivative-shifted-pair}.
Parametrizing and interchanging the order of integration, we see that
\begin{equation}
    I_4 \ = \ \frac{1}{(2\pi i)^2} \int_{1-a+iT}^{1-a} \int_{b}^{b+iT} \frac{L'}{L}(w, f_d) \, \frac{L'}{L}(z, f_d) \, \varphi(-i(z-w))\,dw\,dz,
\end{equation}
where here the paths of integration are line segments as indicated by the bounds. Note that the integral along the segment with real part $1-a$ is traversed from top to bottom due to the counterclockwise orientation of $\mathcal{C}_1$. Switching the bounds on the outer integral therefore yields
\begin{equation}
    I_4 \ = \ -\frac{1}{(2\pi i)^2} \int_{1-a}^{1-a+iT} \int_{b}^{b+iT} \frac{L'}{L}(w, f_d) \, \frac{L'}{L}(z, f_d) \, \varphi(-i(z-w))\,dw\,dz.
\end{equation}
Now, we perform the change of variable $z = w + \rho$ and interchange the order of integration to get
\begin{equation}
    I_4 \ = \ -\frac{1}{(2\pi)^2 i} \int_{a+b-1-iT}^{a+b-1+iT} \varphi(-i\rho) \int_{T_1}^{T_2} \frac{L'}{L}(1-a+it, f_d) \frac{L'}{L}(1-a+it+\rho)\,dt\,d\rho,
\end{equation}
where
\begin{equation}
    T_1 \ = \ \max\{0, -\text{Im}(\rho)\} \quad \text{and} \quad T_2 \ = \ \min\{T, T - \text{Im}(\rho)\}.
\end{equation}
Now, we again use the functional equation, this time in the form
\begin{equation}
    \frac{L'}{L}(1-a+it, f_d) \ = \ \frac{\Phi_d'}{\Phi_d}(1-a+it) - \frac{L'}{L}(a-it, \overline{f}_d).
\end{equation}
Moreover, the term with $\Phi_d'/\Phi_d$ is small, which may be seen by moving the contour to the right. Letting $s = 1/2 + it$, which runs along the critical line, we get
\begin{align}
	I_4 \ & = \ \frac{1}{(2\pi)^2 i} \int_{a+b-1-iT}^{a+b-1+iT} \varphi(-i\rho) \int_{T_1}^{T_2} \frac{L'}{L}\bigg(s+\left(\frac{1}{2}-a+\rho\right), f_d\bigg) \nonumber \\
    & \qquad \times \ \frac{L'}{L}\bigg(1-s + \left(a - \frac{1}{2}\right), \overline{f}_d\bigg)dt\,d\rho + O(T^{\varepsilon}).
\end{align}
Using the formula for the average of the logarithmic derivatives of shifted $L$-functions and arguing akin to $I_3$, we get
\begin{align}
    I_4 \ = \ &\frac{1}{(2\pi)^2 i} \int_{a+b-1-iT}^{a+b-1+iT} \varphi(-i\rho) \int_{T_1}^{T_2} \bigg(\frac{L'_{\star}}{L_{\star}}\bigg)'(1+\rho, f_d \otimes \overline{f}_d) \nonumber \\
    &+ \ \frac{1}{c^2_{f_d}} \bigg( \frac{\sqrt{M}|d| t}{2\pi} \bigg)^{-2\rho} L(1+\rho, f_d \otimes \overline{f}_d) L(1-\rho, f_d \otimes \overline{f}_d) \nonumber \\
    & \times \ A_L\left(\frac{1}{2}-a, -\frac{1}{2}+a-\rho, \frac{1}{2}-a+\rho, a-\frac{1}{2}\right) + \mathscr{C}(1+\rho) \,dt\,d\rho + O(T^{1/2 + \varepsilon}). 
\end{align}
It is readily verified that
\begin{equation*}\label{equivalence}
    A_L\left(-\frac{1}{2}+a+\rho, \frac{1}{2}-a, a-\frac{1}{2}, \frac{1}{2}-a-\rho\right) =  A_L\left(\frac{1}{2}-a, -\frac{1}{2}+a-\rho, \frac{1}{2}-a+\rho, a-\frac{1}{2}\right).
\end{equation*}
Let $\mathscr{A}(\rho)$ denote the above. Reusing the same notation as in the computations for $I_3$ and performing an analogous extension in the interval of integration, we get our desired equation, which is allowed by \eqref{equivalence}.
\end{proof}
\begin{proof}
We now prove Proposition \ref{contour_theorem}.
Adding up the results of $I_1, I_2, I_3, I_4$, we get
\begin{align}
    P(f_d; \varphi) \ &= \ \frac{4}{(2\pi)^2} \int_{-T}^{T} \varphi(\rho) \int_{0}^{T} \log(Av)^2 \,dv\,d\rho \nonumber \\
    &\quad + \ \frac{1}{(2\pi)^2 i} \int_{0}^{T} \int_{-\delta-iT}^{-\delta+iT} \varphi(-i\rho) \mathscr{I}(-\rho,t)\,d\rho\,dt \nonumber \\
    &\quad + \ \frac{1}{(2\pi)^2 i} \int_{0}^{T} \int_{\delta-iT}^{\delta+iT} \varphi(-i\rho) \mathscr{I}(\rho,t)\,d\rho\,dt + O(T^{1/2+\varepsilon}). 
\end{align}
We now focus our attention on $I_3 + I_4$, which are the final two integrals in the equation above. 
We perform the change of variable $-\rho \mapsto \rho$ in $I_3$.
Since $\varphi$ is even, we see that
\begin{equation}
    I_3 + I_4 \ = \ \frac{2}{(2\pi)^2 i} \int_{0}^{T} \int_{\delta-iT}^{\delta+iT} \varphi(i\rho) \mathscr{I}(\rho,t)\,d\rho\,dt.
\end{equation}
We have
\begin{align}\label{almost_there}
    P(f_d; \varphi) \ &= \ \frac{1}{2\pi^2} \bigg(\frac{1}{i} \int_{0}^{T} \int_{\delta-iT}^{\delta+iT} \varphi(i\rho) \mathscr{I}(\rho,t)\,d\rho\,dt \nonumber \\
		      & \qquad \qquad + \ \int_{0}^{T} \int_{-T}^{T} \varphi(\rho) 2\log(At)^2 \,d\rho\,dt \bigg) + O(T^{1/2+\varepsilon})
		      .\end{align}
By the computations for $I_3$ and $I_4$, we have that $\mathscr{I}(\rho,t)$ takes the form
\begin{align}\label{script_I}
    \mathscr{I}(\rho, t) \ = \ &\bigg( \frac{L_{\star}'}{L_{\star}} \bigg)'(1+\rho, f_d \otimes \overline{f}_d) + \mathscr{C}(1+\rho) \nonumber \\
    &+ \ \frac{1}{c_{f_d}^2} \bigg( \frac{\sqrt{M} |d| t}{2\pi} \bigg)^{-2\rho} L(1+\rho, f_d\otimes \overline{f}_d) L(1-\rho, f_d \otimes \overline{f}_d) \mathscr{A}(\rho).
\end{align}
Analogous to \cite{CS07}, we have that $\mathscr{A}'(0) = 0$. Moreover, recall $\mathscr{A}(0) = 1$.
By the factorization of $L(s, f_d \otimes \overline{f}_d)$ given earlier, along with the fact that $(L_{\star}'/L_{\star})'$ differs from $(L'/L)'$ by a sum of a finite number of primes constituting an entire function, we may state the Laurent series about $\rho=0$ of $\mathscr{I}(\rho,t)$ is given by
\begin{equation}\label{laurent}
    \mathscr{I}(\rho,t) \ = \ \frac{2 \log(\sqrt{M} |d| \cdot t /2\pi)}{\rho} + O(1).
\end{equation}
Finally, we move the path of integration in $\rho$ in the second term of \eqref{almost_there} to the imaginary axis. To this end, we let $\delta \rightarrow 0$. We use the principal value while passing through zero as it is clear from \eqref{laurent} that $\mathscr{I}(\rho,t)$ has a simple pole at $\rho=0$.

We cannot simply, however, shift the path of integration to the imaginary axis due to the aforementioned pole at the origin. Instead, upon shifting to the imaginary axis, we must additionally consider a small counterclockwise semicircular path around the origin of radius $\varepsilon$ to avoid the pole. Taking $\varepsilon \rightarrow 0$ recovers the original segment. However, we now have an additional contribution from the semicircular piece as a result of the fractional residue theorem (see \cite{G01}). In particular, this contribution yields
\begin{equation}
    \int_{\delta-iT}^{\delta+iT} \varphi(i\rho) \mathscr{I}(\rho,t)\,d\rho = \int_{-iT}^{iT} \varphi(i\rho) \mathscr{I}(\rho,t)\,d\rho + \pi i \cdot \text{Res}(\varphi(i\rho)\mathscr{I}(\rho,t), \rho=0).
\end{equation}
The factor of $\pi$ is precisely the angle of the arc. We get
\begin{align}
    \frac{1}{i} \int_{0}^{T} \int_{\delta-iT}^{\delta+iT} \varphi(i\rho) \mathscr{I}(\rho,t)\,d\rho\,dt \ &= \ \frac{1}{i} \int_{0}^{T} \int_{-iT}^{iT} \varphi(i\rho) \mathscr{I}(\rho,t)\,d\rho\,dt \nonumber \\
													  & + \ 2\pi \int_{0}^{T} \varphi(0) \log\bigg( \frac{\sqrt{M} |d| t}{2\pi} \bigg)\,dt
.\end{align}
As a result, substituting into \eqref{almost_there} then gives
\begin{align}
    P(f_d; \varphi) \ = \ &\frac{1}{2\pi^2} \bigg( \frac{1}{i}\int_{0}^{T} \int_{-iT}^{iT} \varphi(i\rho) \mathscr{I}(\rho,t)\,d\rho\,dt + 2\pi \int_{0}^{T} \varphi(0) \log\bigg( \frac{\sqrt{M} |d| t}{2\pi} \bigg)\,dt \nonumber\\
    &+ \ \int_{0}^{T} \int_{-T}^{T} \varphi(\rho) 2\log\bigg( \frac{\sqrt{M} |d| t}{2\pi} \bigg)^2 \,d\rho\,dt \bigg) + O(T^{1/2+\varepsilon}).
\end{align}
In the first double integral, we change variables $\rho \mapsto ir$ and exploit the evenness of $\varphi$, and in $I_3$, we simply change $\rho$ to $r$.
Expanding everything out, we get
\begin{align}
    P(f_d; \varphi) \ & = \ \frac{1}{2\pi^2} \int_{0}^{T} \bigg[ 2\pi \varphi(0) \log\bigg( \frac{\sqrt{M} |d| t}{2\pi} \bigg) + \int_{-T}^{T} \varphi(r) \bigg( 2 \log^2\bigg( \frac{\sqrt{M} |d| t}{2\pi} \bigg) \nonumber\\
    &\quad + \ \bigg( \frac{L_{\star}'}{L_{\star}} \bigg)'(1+ir, f_d \otimes \overline{f}_d) + \mathscr{C}(1+ir) \nonumber \\
    &\quad + \ \frac{1}{c_{f_d}^2} \bigg( \frac{\sqrt{M} |d| t}{2\pi} \bigg)^{-2ir} L(1+ir, f_d\otimes \overline{f}_d) L(1-ir, f_d \otimes \overline{f}_d) \mathscr{A}(ir) \bigg) \bigg]\,dr\,dt \nonumber \\
    &\qquad +\ O(T^{1/2+ \varepsilon}).
\end{align}
Finally, putting
\begin{align}
    I_r \ \coloneqq \ &\int_{-T}^{T} \varphi(r) \bigg( 2\log^2\bigg( \frac{\sqrt{M} |d| t}{2\pi} \bigg) + \bigg( \frac{L'_{\star}}{L_{\star}} \bigg)' (1 + ir, f_d \otimes \overline{f}_d) \\
    &+ \ \frac{1}{c_{f_d}^2} \bigg( \frac{\sqrt{M} |d| t}{2\pi} \bigg)^{-2ir} L(1+ir, f_d \otimes \overline{f}_d)L(1-ir, f_d \otimes \overline{f}_d) \mathscr{A}(ir) + \mathscr{C}(1+ir) \bigg)\,dr \nonumber,
\end{align}
we conclude that
\begin{equation}
    P(f_d; \varphi) \ = \ \frac{1}{2\pi^2} \int_{0}^{T} \bigg[ 2\pi\varphi(0) \log\bigg( \frac{\sqrt{M}|d| t}{2\pi} \bigg) + I_r \bigg]\,dt + O(T^{1/2 + \varepsilon}).
\end{equation}
This gives the desired formula for $P(f_d; \varphi)$.
\end{proof}
\subsection{Series expansion of pair correlation}\label{subsect:series_pair_correlation}
With a formula for $P(f_d; \varphi)$, we obtain a series development for large $T$ and use it to obtain the effective matrix size.
We scale the pair-correlation by substituting 
\begin{equation}
    y \ \coloneqq \ rR/\pi \quad \text{and} \quad R \ \coloneqq \ \log \bigg( \frac{\sqrt{M} |d| T}{2\pi e} \bigg).
\end{equation}
Define the \emph{rescaled test function} $g$ by $g(y) \coloneqq \varphi(r)$.
By Proposition \ref{contour_theorem} and changing variables within the integral via $r \mapsto rR/\pi= y$, we observe
\begin{equation}\label{changeofvariable}
    \sum_{0 < \gamma,\gamma' < T} g\bigg((\gamma-\gamma')\frac{R}{\pi}\bigg) \ = \ \frac{1}{2\pi^2} \int_{0}^{T} \bigg[ 2\pi g(0) \log\bigg( \frac{\sqrt{M}|d| t}{2\pi} \bigg) + I_y\bigg]\,dt + O(T^{1/2 + \varepsilon}),
\end{equation}
where
\begin{align}
    I_y \ \coloneqq \ &\frac{\pi}{R} \int_{-T(R/\pi)}^{T(R/\pi)} g(y) \bigg( 2\log^2\bigg( \frac{\sqrt{M} |d| t}{2\pi} \bigg) + \bigg( \frac{L'_{\star}}{L_{\star}} \bigg)' \bigg(1 + \frac{i\pi y}{R}, f_d \otimes \overline{f}_d\bigg) \nonumber\\
    &+ \ \frac{1}{c_{f_d}^2} \bigg( \frac{\sqrt{M} |d| t}{2\pi} \bigg)^{-2i\pi y/R} L\bigg(1+\frac{i\pi y}{R}, f_d \otimes \overline{f}_d\bigg)L\bigg(1-\frac{i\pi y}{R}, f_d \otimes \overline{f}_d\bigg) \nonumber \\
    & \quad \times \ \mathscr{A}\bigg(\frac{i\pi y}{R}\bigg) + \mathscr{C}\bigg(1+\frac{i\pi y}{R}\bigg) \bigg)\,dy.
\end{align}
We explicitly compute the outer integral over $t$. First, we note that
\begin{equation}
    \int_{0}^{T} \log\bigg( \frac{\sqrt{M}|d| t}{2\pi} \bigg) \,dt \ = \ T \bigg( \log\bigg( \frac{\sqrt{M} |d| T}{2\pi}  \bigg)- 1 \bigg) \ = \ T\log\bigg( \frac{\sqrt{M} |d| T}{2\pi e}\bigg).
\end{equation}
Substituting the above into \eqref{changeofvariable}, we obtain
\begin{equation}
    \sum_{0 < \gamma,\gamma' < T} g\bigg((\gamma-\gamma')\frac{R}{\pi}\bigg) \ = \ \frac{g(0)}{\pi}T\log\bigg( \frac{\sqrt{M} |d| T}{2\pi e}\bigg) + \frac{1}{2\pi^2} \int_{0}^{T} I_y\,dt + O(T^{1/2 + \varepsilon})
,\end{equation}
with $I_y$ as above. The integral with respect to $t$ of $I_y$ is a double integral and can be evaluated by carefully keeping track of terms with $t$-dependence and applying the Fubini theorem. The result is
\begin{align}
    \int_{0}^{T} I_y\,dt \ & = \ \frac{2\pi T(1+R^2)}{R} \int_{-TR/\pi}^{TR/\pi} g(y) \,dy \nonumber\\
    &\quad + \ \frac{\pi T}{R}\int_{-TR/\pi}^{TR/\pi}g(y) \bigg(\frac{L_{\star}'}{L_{\star}}\bigg)' \bigg(1+\frac{i\pi y}{R}, f_d \otimes \overline{f}_d\bigg)dy \nonumber \\
    &\quad + \ \frac{\pi T}{Rc_{f_d}^2} \int_{-TR/\pi}^{TR/\pi} \frac{e^{-2\pi i y(1 + 1/R)}}{1 - 2\pi i y /R} g(y) L\bigg(1+\frac{i\pi y}{R}, f_d \otimes \overline{f}_d\bigg) \nonumber\\
    & \times \ L\bigg(1-\frac{i\pi y}{R}, f_d \otimes \overline{f}_d\bigg) \mathscr{A}\bigg(\frac{i\pi y}{R}\bigg)dy \nonumber \\
    &\quad + \ \frac{\pi T}{R} \int_{-TR/\pi}^{TR/\pi} g(y) \cdot \mathscr{C}\bigg(1+\frac{i\pi y}{R}\bigg) dy .
\end{align}
We recall that $\varphi(x) \ll 1/(1+x^2)$ for real $x$, implying $g(y)\ll 1/(1+y^2)$ for real $y$. Moreover, we are interested in obtaining a series expansion for $P(f_d; \varphi)$ when $T$ is large. For large $T$, the expression for $P(f_d; \varphi)$ is well-approximated by
\begin{align}
    \sum_{0 < \gamma,\gamma' < T} g\bigg((\gamma-\gamma')\frac{R}{\pi}\bigg) \ = \ &\frac{T}{\pi}\log \bigg( \frac{\sqrt{M} |d| T}{2\pi e} \bigg) \bigg[ g(0) + \frac{1}{2R^2 T} \int_{\mathbb{R}} g(y) \bigg( 2T(1+R^2) \nonumber \\
    &+ \ T \bigg(\frac{L_{\star}'}{L_{\star}}\bigg)' \bigg(1+\frac{i\pi y}{R}, f_d \otimes \overline{f}_d\bigg) \nonumber\\
    &+ \ \frac{1}{c_{f_d}^2} \frac{Te^{-2\pi i y(1 + 1/R)}}{1 - 2\pi i y /R} L\bigg(1+\frac{i\pi y}{R}, f_d \otimes \overline{f}_d\bigg)\nonumber\\
    & \times \ L\bigg(1-\frac{i\pi y}{R}, f_d \otimes \overline{f}_d\bigg) \mathscr{A}\bigg(\frac{i\pi y}{R}\bigg) \nonumber\\
    &+ \ T \cdot \mathscr{C}\bigg(1+\frac{i\pi y}{R}\bigg)
    \bigg)\,dy \bigg
] + O(T^{\varepsilon+1/2}),
\end{align}
as any extra constants from integrating over all of $\mathbb{R}$ are absorbed into the error term. This simplifies down to
\begin{align}\label{beforeexpansion}
    \sum_{0 < \gamma,\gamma' < T} g\bigg((\gamma-\gamma')\frac{R}{\pi}\bigg) \ = \ &\frac{T}{\pi}\log \bigg( \frac{\sqrt{M} |d| T}{2\pi e} \bigg) \bigg[ g(0) \nonumber\\
										   & + \ \int_{\mathbb{R}} g(y) \bigg( \frac{1+R^2}{R^2} + \frac{1}{2R^2} \bigg(\frac{L_{\star}'}{L_{\star}}\bigg)' \bigg(1+\frac{i\pi y}{R}, f_d \otimes \overline{f}_d\bigg) \nonumber\\
    &+ \ \frac{1}{2R^2c_{f_d}^2} \frac{e^{-2\pi i y(1 + 1/R)}}{1 - 2\pi i y /R} L\bigg(1+\frac{i\pi y}{R}, f_d \otimes \overline{f}_d\bigg)\nonumber\\
    & \times \ L\bigg(1-\frac{i\pi y}{R}, f_d \otimes \overline{f}_d\bigg) \mathscr{A}\bigg(\frac{i\pi y}{R}\bigg) \nonumber\\
    &+ \ \frac{1}{2R^2} \mathscr{C}\bigg(1+\frac{i\pi y}{R}\bigg) 
   \bigg)\,dy \bigg
] + O(T^{\varepsilon+1/2}).
\end{align}
We now perform a series expansion in $1/R$ for $R$ large. Most of the computations are routine and more details are provided in Appendix \ref{appendix:pair_correlation_series}, but critical to the argument is the fact that the Ramanujan-Petersson conjecture is known for holomorphic cusp forms. This allows us to expand $\mathscr{B}$ in a series about 1. Putting
\begin{align}
    e_1 \ &\coloneqq \ \frac{1}{2} \cdot \frac{\log(M)^2}{ |\lambda_{f_d}(M)|^{-2} \cdot M - 1} ,\\
    e_2 \ &\coloneqq \ -2 + \gamma^2 + 2\gamma_1 - \frac{\mathscr{A}''(0)}{2} - \bigg(\frac{L'}{L}\bigg)'(1, \text{ad}^2\,f_d), \\
    e_3 \ &\coloneqq \ \frac{16 + \mathscr{A}'''(0)}{12},
\end{align}
we get
\begin{align}
    \sum_{0 < \gamma,\gamma' < T} g\left((\gamma-\gamma')\frac{R}{\pi}\right) \ = \ &\frac{T}{\pi}\log \bigg( \frac{\sqrt{M} |d| T}{2\pi e} \bigg) \bigg[ g(0) + \int_{\mathbb{R}} g(y) \bigg( 1 - \bigg(\frac{\sin\pi y}{\pi y}\bigg)^2 \\
    &+ \frac{e_1 - e_2 \sin^2 \pi y}{R^2} - \frac{e_3 \pi y \sin 2\pi y}{R^3} + O(R^{-4})
    \bigg)\,dy \bigg] + O(T^{\varepsilon+1/2}), \nonumber
\end{align}
where $\gamma$ is the Euler-Mascheroni constant and $\gamma_1$ is the first Stieltjes constant. Defining a function $h$ by $h(y) = 2 g(y)$ (and thus $h$ also satisfies $h(y) \ll 1/(1+y^2)$), we may simplify to get
\begin{align}\label{afterexpansion}
    P(f_d;\varphi) \ = \ &\frac{T}{2\pi}\log \bigg( \frac{\sqrt{M} |d| T}{2\pi e} \bigg) \bigg[ h(0) + \int_{\mathbb{R}} h(y) \bigg( 1 - \bigg(\frac{\sin\pi y}{\pi y}\bigg)^2 \nonumber\\
    &+ \frac{e_1 - e_2 \sin^2 \pi y}{R^2} - \frac{e_3 \pi y \sin 2\pi y}{R^3} + O(R^{-4})
    \bigg)\,dy \bigg] + O(T^{\varepsilon+1/2}).
\end{align}
This yields the asymptotic
\begin{align}\label{asymptotic_pair_corr}
    P(f_d;\varphi) \ \sim \ N(T, f_d) \bigg[ h(0) + \int_{\mathbb{R}} h(y) \bigg( 1 - \bigg(\frac{\sin\pi y}{\pi y}\bigg)^2 \bigg)\,dy \bigg],
\end{align}
where $N(T, f_d)$ is the number of zeros $\rho$ with $\Re\,\rho \in [0,1]$ and $\Im\,\rho \in (0,T]$ \cite[Chapter 24]{IK04}. It is worth noting that the expression
\begin{equation}
    1-\bigg(\frac{\sin\pi y}{\pi y}\bigg)^2
,\end{equation}
is precisely the limiting two-point correlation function predicted by Montgomery \cite{Mon73}. Therefore, assuming the Ratios Conjectures, Montgomery's conjecture holds in the general setting of $L$-functions associated to Hecke forms of level an odd prime.

Following \cite{Con05}, the scaled pair-correlation $Q_{\U(N)}(x)$ for $U(N)$ is
\begin{equation}
    Q_{\U(N)}(x) \ = \ 1 - \bigg( \frac{\sin \pi x}{\pi x} \bigg)^2 - \frac{\sin^2 \pi x}{3N^2} + O(N^{-4}).
\end{equation}
Due to the presence of the $e_1$ term in our pair-correlation expansion, we are unable to directly match coefficients to obtain the effective matrix size. To deal with this, we minimize
\begin{equation}
    \|f(N)\|_{L^2}^2 \ = \ \int_{-t}^{t} \left| \frac{e_1 - e_2 \sin^2 \pi y}{R^2} + \frac{\sin^2 \pi y}{3N^2} \right|^2 dy
\end{equation}
for each integer $t$. We are free to do this since the integrand has unit period. Expanding and calculating the integral and optimizing with respect to $N$ in the usual manner (i.e., by differentiating with respect to $N$) yields that
\begin{equation}
    N \ = \ \frac{R}{\sqrt{3e_2 - 4e_1}}
\end{equation}
minimizes $\|f(N)\|_{L^2}^2$. Therefore, by optimizing with respect to the $L^2$ norm, we arrive at the effective matrix size of a given generic form
\begin{equation}
    N_{\text{eff}} \ = \ \frac{R}{\sqrt{3e_2 - 4e_1}}.
\end{equation}
This is the effective matrix size for the $L$-function attached to the fixed twisted generic form $f_d$. Therefore, for our family $\mathcal{F}_f^{+}(X)$ of cusp forms, all we need to do is take the average value of $e_1$ and $e_2$ as $d$ varies. 
This gives us an effective matrix size 
\begin{equation}
    N_{\text{eff}} \ = \ \frac{R}{\sqrt{3\langle e_2 \rangle - 4 \langle e_1 \rangle}}
\end{equation}
for our family provided that $3\langle e_2 \rangle - 4 \langle e_1 \rangle > 0$. (This is a reasonable restriction, especially for large $M$.)

\begin{remark}
The appearance of the Rankin-Selberg convolution of the cusp form $f_d$ in Proposition \ref{prop:int-unramified} is not a surprise.
We compare to Theorem 2.5 of \cite{CS07}.
The Rankin-Selberg convolution of $\zeta(s)$ with itself is just $\zeta(s)$.
We see that Proposition \ref{prop:int-unramified} has the same shape as in the corresponding theorem for $\zeta(s)$ and constitutes a fairly generic expression for the average of a product of shifted logarithmic derivatives of automorphic $L$-functions. Such a formula should generally be controlled by the $L$-function attached to the appropriate automorphic Rankin-Selberg lift.
It is important to note that on $\GL(2)$, the Rankin-Selberg lift counts $\zeta(s)$ as a factor.
It is likely that moments of all automorphic functions also involve contributions from $\zeta(s)$ to the main term.
\end{remark}

\section{Cutoff value}
The twist of a given form $f \in S_k^{\text{new}}(M,\text{principal})$ by a quadratic character $\psi $ has level $Md^2$, weight $k$, and nebentype $\chi_{f_d} = \chi_f\psi_d^2$ \cite[Proposition 14.19]{IK04}. Therefore, for a given form, the density of zeros near the central point $1/2$ is asymptotic to $\log d$. This yields the relation $N_{\text{std}} \ \sim \ \log d$. Thus, the values at the central point are no longer discretized on a scale of $1/\sqrt{d}$; instead, they are discretized on a scale of $\exp( (1-k)N_{\text{std}}/2)$.
If we restrict to even forms, we now want
matrices in $\SpO(2N)$ satisfying
\begin{align}
\label{eq:exciseddef} |\Lambda_A(1,N)| \ \geq \ \ c \cdot\exp((1-k)N_{\text{std}}/2). \end{align}
We now optimize the cutoff value $c$ of the excision threshold.
Our point of departure is \cite{CKRS05} and \cite{CKRS06} as modified by \cite{DHKMS12}.

\subsection{Cutoff value for forms with principal nebentypus}\label{sct:cutoff_value_for_forms_with_principal_nebentypus}
We may apply Kohnen-Zagier's formula \eqref{eqn:kohnen-zagier-formula} as our family of principal forms satisfies the assumptions for the formula
\begin{align}
	L(1 /2,f_d) \ < \ \frac{|c(d)|^2\kappa_f}{|d|^{k-1/2}} \ \implies \ L(1 /2,f_d) \ = \ 0
.\end{align}
The coefficients $c(d)$ are the Fourier coefficients of a half-integral weight modular form that is obtained via generalized Shimura correspondence.
It was Waldspurger (see \cite{Wal80,Wal81}) who related the Fourier coefficients of this half-integral weight form obtained via the Shimura correspondence to the central values of $L$-functions arising from quadratic twists of an integral weight form.
Kohnen and Zagier \cite{KZ81} then used this Waldspurger formula to establish a discretization at the central point.

Restricting to $\chi_f$ principal, the arithmetic of the coefficients $c(d)$ is not well-understood.
The authors in \cite[Section 4]{CKRS06} use conjectures arising from random matrix theory for the value distribution of elliptic curve $L$-functions to make a series of conjectures on the statistics of these coefficients $c(d)$.
The more complex arithmetic of the coefficients $c(d)$ are not accessible via random matrix theory.
Hence, the authors make use of a numerically generated proportionality constant \cite[Section 5]{CKRS06}.

Though a notion of Shimura correspondence exists for $\chi_f$ non-principal, there is no analogous formula of Waldspurger type.
For this reason, it is not known whether the values $L(1 /2,f_d)$ are discretized for forms where $\chi_f$ non-principal.
We also have no results which would allow us to predict the coefficient of the main term of the frequency of vanishing.
For this reason, theoretically deriving behavior of low-lying zeros of $L$-functions attached to forms with non-principal nebentypus remains out of reach.

One way to predict the frequency of vanishing of our family $\mathcal{F}_f^{+}(X)$ of quadratic twists of a modular form with principal nebentypus would be to calculate the probability that a random variable $Y_d$ with probability density $P_f(d,x)$ assumes a value less than $\kappa_f|c(d)|^2d^{(1-k) /2}$ and then sum asymptotically over the family.
This method was pioneered by \cite{CKRS05} and \cite{CKRS06}.
However, to obtain an asymptotic for the frequency of vanishing with the correct leading coefficient, we need to have an idea of the statistical behavior of the coefficient $c(d)$ evaluated at $\mathcal{D}_f^{+}(X)$.

For this reason the authors in \cite{DHKMS12} determine this value numerically; that is, they introduce a notion of an `effective' cutoff depending on a parameter that is determined numerically, but which does not depend on $d$.
We also determine this value numerically and write $\delta_f \kappa_f d^{(1-k) /2}$ for our `effective' cutoff, where $\delta_f$ is a numerical input.
In analogy with Equation (5.12) in \cite{DHKMS12}, we write
\begin{align}
\Prob\bigg( 0 \leq Y_d \ \leq \ \frac{\delta_f \kappa_f}{|d|^{k-1/2}} \bigg) \ & \sim \ \int_0^{\delta_f\kappa_f d^{k-1/2}} a_f(-1 /2) \ h(\log d) \ x^{-1 /2} \ dx \\	
 & = \ 2 a_f(-1 /2) \ h(\log d)\ \frac{\sqrt{\delta_f \kappa_f}}{|d|^{k/2 - 1/4}}
,\end{align}
where
\begin{align}\label{eqn:barnes-G-asymptotic}
h(N) \ \sim \ 2^{-7/8} G(1/2)\pi^{-1/4}N^{3/8}
\end{align}
is the asymptotic for the moment generating function of the corresponding group for $\mathcal{F}_f^{+}(X)$ at the pole with $G$ the Barnes $G$-function.
Following \cite{CKRS05} and \cite{DHKMS12}, we conjecture that 
\begin{align}\label{eqn:conj-starred-sum}
	\bigg| \bigg\{ L_f(s,\psi_d) \in \mathcal{F}^{+}_f(X) : d \text{ prime}, \ L_f(1 /2,\psi_d) = 0 \bigg\}  \bigg| = {{\sum_{\substack{d\leq X \\ d \text{ prime}}}}^*}\Prob\bigg( 0\leq Y_d \leq \frac{\delta_f\kappa_f}{|d|^{k-1/2}} \bigg), 
\end{align}
where the starred sum in the middle line means that the sum is only over those \emph{prime} fundamental discriminants for which $\chi_d(-M)\epsilon_f = +1$ (of which, asymptotically, there are $X /4\log X$).
The convergence of the sum \eqref{eqn:conj-starred-sum} is determined by $k$. Namely, if $k< 3$, then the sum diverges, and we have the asymptotic
\begin{align}
\eqref{eqn:conj-starred-sum} \ & \sim \ \frac{1}{4\log X} \sum_{n=1}^{\lfloor X \rfloor} 2a_f(-1 /2) h (\log X) \frac{\sqrt{\delta_f\kappa_f}}{n^{k/2-1/4}}  \\
& \sim \ \frac{1}{4\log X}2a_f(-1 /2) \sqrt{\delta_f\kappa_f} 2^{-7/8} G(1/2)\pi^{-1/4}N^{3/8}  \frac{4}{5-2k} X^{(5-2k) /4},
\end{align}
where the last line comes from \ref{eqn:barnes-G-asymptotic}.
With the same starred notation, we have
\begin{align}
	\sqrt{\delta_f} \ = \ \frac{\left| \left\{ L_f(s,\psi_d) \in \mathcal{F}_f^+ : d\text{ prime}, \ L_f(1 /2,\psi_d) = 0 \right\}  \right| }{\sum_{d\text{ prime}}^{*}2a_f(-1 /2) h(\log d) \sqrt{\kappa_f}|d|^{1/2-k}}
.\end{align}
\begin{remark}
On the other hand, if $k \geq 3$, then the sum \eqref{eqn:conj-starred-sum} converges.
This would heuristically mean there is no discretization for the zeros of twisted $L$-functions attached to a form of principal nebentype with weight 6 or greater. We expect there to be little to no repulsion for principal forms with weight equal to 4; and we explain how the numerical observations for \texttt{5.4.a.a} and \texttt{7.4.a.a} in Section \ref{sec:numerical-observations} support this conclusion.
\end{remark}

\subsection{Translating the cutoff for \texorpdfstring{$\SpO(\mathrm{even})$}{SO(even)}, weight 2}

Since repulsion from the origin rapidly decreases from weight 2 to 4, we really only need to compute the cutoff for newforms with principal nebentype and weight 2. Since it suffices to only consider the elliptic curve case when numerically calculating the cutoff value, we may employ the same process as in \cite[Section 5.2]{DHKMS12}.

\section{Numerical Observations}
\label{sec:numerical-observations}
We reiterate that we did not compute the effective matrix size for our family due to the given cuspidal newform $L$-function's Euler product not being accessible in PARI/GP or SageMath.
Recall the nearest neighbour spacing statistic is the probability density for distances between consecutive zeros, or equivalently, a normalized histogram with bin size 100 of gaps between consecutive zeros. We work with the following cuspidal newforms.
\begin{center}
\begin{tabular}{c|c|c}
    LMFDB Label & Type & Group \\
    \hline
    \texttt{11.2.a.a}, \texttt{5.4.a.a},\texttt{5.8.a.a}, \texttt{7.4.a.a} & $\chi_f$ principal, even twists & $\SpO(2N)$\\
    \hline
    \texttt{11.2.a.a}, \texttt{5.4.a.a},
    \texttt{5.8.a.a}, \texttt{7.4.a.a} & $\chi_f$ principal, odd twists & $\SpO(2N+1)$\\
    \hline
    \texttt{3.7.b.a}, \texttt{7.3.b.a} & self-CM & $\USp(2N)$\\
    \hline
    \texttt{13.2.e.a}, \texttt{11.7.b.b}, \texttt{7.4.c.a}, \texttt{17.2.d.a} & generic & $\U(N)$\\
\end{tabular}
\end{center}

In the following, we consider the `first' eigenvalues of random matrices $\SpO(\text{odd})$ and $\SpO(\text{even})$ whose characteristic polynomials are evaluated at or near 1. By `first,' we mean those eigenvalues closest to 1 on the unit circle. Note the eigenvalues of random matrices from $\SpO(\text{odd})$ are all going to be zero. 

\subsection{Families with orthogonal symmetry}
\begin{center}
\begin{figure}[htpb]	
	\begin{tabular}{c c}
		Lowest zeros (even twists) & Lowest zeros (odd twists) \\
	\includegraphics[scale=0.45]{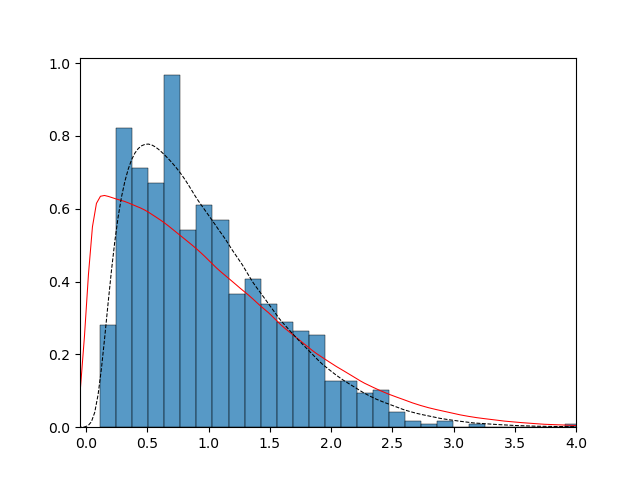} & \includegraphics[scale=0.45]{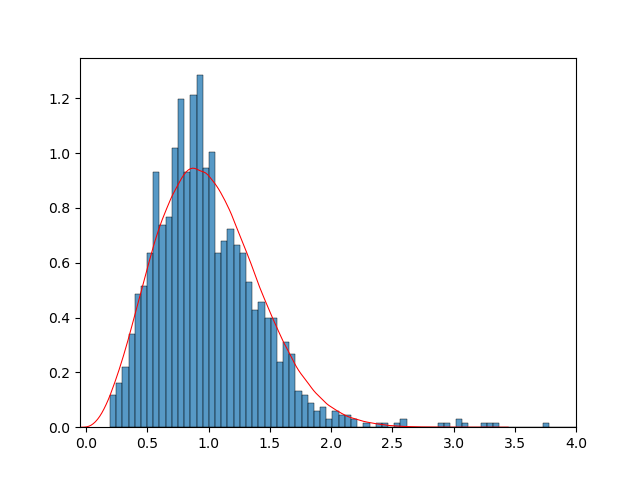}
	\end{tabular}
	\caption{The left histogram shows the distribution of lowest zeros for 1,380 even twists of \texttt{11.2.a.a} with discriminant up to 9,960; the red curve (left) shows the distribution of the first eigenvalues of 1,000,000 randomly generated $\text{SO}(18)$ matrices with characteristic polynomial evaluated at 1; the black dotted curve (left) is the same distribution but with excision. We varied the excision threshold numerically to obtain the optimal fit. The right histogram shows the distribution of second lowest zeros for 1,394 odd twists of \texttt{11.2.a.a} with discriminant up to 9,957, and the red line (right) shows the distribution of the first eigenvalues of 1,000,000 randomly generated $\text{SO}(19)$ matrices with characteristic polynomial evaluated near 1. The data have been normalized to have mean 1.}
 \label{tab:zeros-mf11w2aa}
\end{figure}	
\end{center}

We numerically computed the non-vanishing lowest-lying zeros of $\mathcal{F}_f^{+}(X)$ with $f \in S_{2}^{\text{new}}(11,\text{principal})$, which has sign $\epsilon_f=+1$. The code may be found in \cite{Yao24}. As expected, we obtain the repulsion from the origin for even and hence, the model requires the cutoff value. Recall we chose those twists with positive fundamental discriminant. If we choose twists with negative fundamental discriminants, we still recover repulsion from the origin. Our results combined with the results in \cite{DHKMS12} means that regardless of if we range over twists with negative or positive discriminants, we still recover repulsion from the origin. In addition, observe in Figure \ref{tab:zeros-mf11w2aa} the rather pronounced repulsion from the origin for non-vanishing lowest-lying zeros of even twists of \texttt{11.2.a.a}. Since the distribution of the lowest-lying zeros agrees most with that of the excised matrices (dotted line) in Figure \ref{tab:zeros-mf11w2aa}, this verifies the necessity of creating such a model.

\begin{center}
\begin{figure}[htpb]	
\begin{tabular}{c c}
Lowest non-vanishing zeros & Lowest non-vanishing zeros \\
\includegraphics[scale=0.45]{Images/nfzc_mf11w2aa_sign-1_SO21.png} & \includegraphics[scale=0.45]{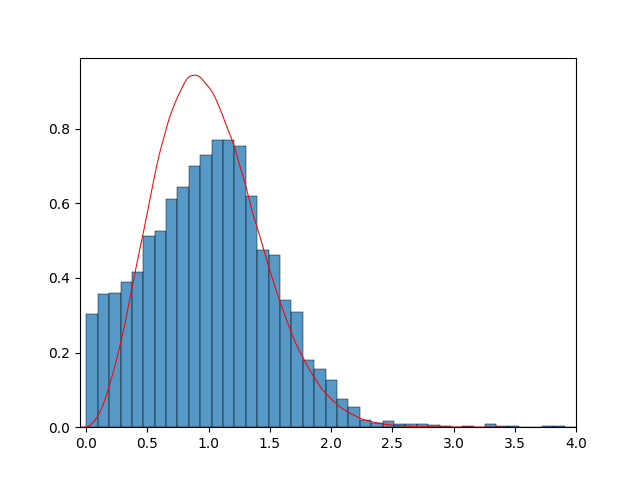}
\end{tabular}
\caption{The left histogram shows lowest non-vanishing zeros of 4,497 odd twists of \texttt{11.2.a.a} with discriminant up to 32,829, and the red curve (right) shows the distribution of the first eigenvalues of 1,000,000 randomly generated $\text{SO}(21)$ matrices with characteristic polynomial evaluated near 1. The right histogram shows lowest non-vanishing zeros of 5,463 odd twists of \texttt{7.4.a.a} with discriminant up to 41,000, and the red curve (right) shows the distribution of the first eigenvalues of 1,000,000 randomly generated $\text{SO}(21)$ matrices with characteristic polynomial evaluated near 1. The data have been normalized to have mean 1.}
\label{tab:attraction-mf11w2aa-mf7w4aa}
\end{figure}	
\end{center}

Recall that eigenvalues of random matrices from $\SpO(\text{odd})$ with characteristic polynomials evaluated at 1 do, in fact, vanish. This $\SpO(\text{odd})$ behavior agrees with the behavior of the (mostly) vanishing lowest-lying zeros of odd twists. We call to attention that those lowest-lying zeros which do not vanishing give rise to interesting behavior. In Figure \ref{tab:attraction-mf11w2aa-mf7w4aa}, the left histogram shows the (normalized) lowest non-vanishing zeros of odd twists that do not vanish. In particular, the non-vanishing lowest-lying zeros odd twists seem to `force' a particular distribution which does not look as natural as that of even twists. In fact, we see an attraction toward the origin in the right histogram of Figure \ref{tab:attraction-mf11w2aa-mf7w4aa}. On the other hand, we plotted the distribution of (non-vanishing) second-lowest lying zeros of odd twists of \texttt{11.2.a.a} in the right histogram of Figure \ref{tab:zeros-mf11w2aa}. If we only consider non-vanishing second lowest-lying zeros, then we find some agreement between zeros and eigenvalues. Future authors might consider implementing a random matrix model that incorporates the ``spiky'' behavior found at the peak.

\begin{center}
\begin{figure}[htpb]
	\begin{tabular}{c c}
	Lowest zeros (even)
 & Non-vanishing lowest zeros (odd) \\
	\includegraphics[scale=0.45]{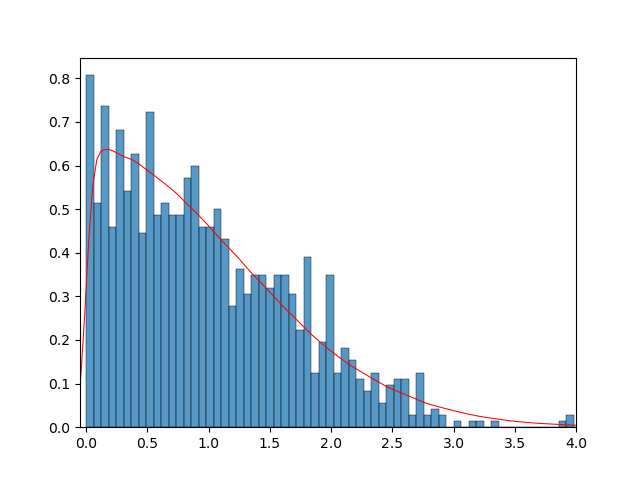}
 & \includegraphics[scale=0.45]{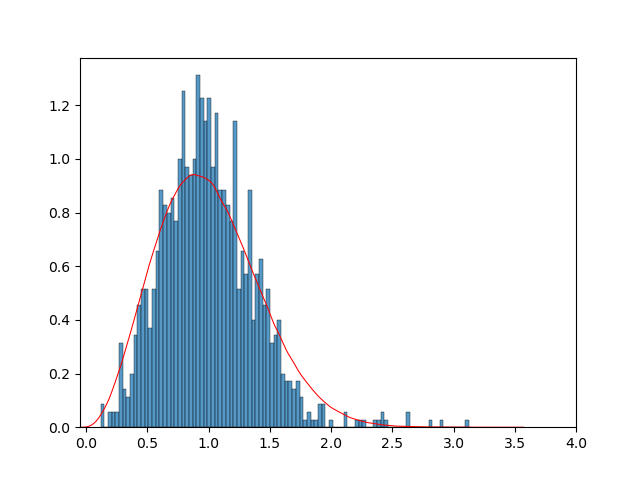}	
	\end{tabular}
	\caption{The left histogram shows non-vanishing lowest-lying zeros of 1,174 even twists of \texttt{5.8.a.a} with discriminant up to 9,228; the red curve (left) shows the distribution of the first eigenvalues of 1,000,000 randomly generated $\SpO(16)$ matrices with characteristic polynomial evaluated at 1 without excision. The right histogram shows the second lowest-lying zeros of 1,165 odd twists of \texttt{5.8.a.a} with discriminant up to 9,229; the red curve (right) shows the distribution of the first eigenvalues of 1,000,000 randomly generated $\text{SO}(17)$ matrices with characteristic polynomial evaluated near 1. The data have been normalized to have mean 1.}
 \label{tab:table-mf5w8aa}
\end{figure}	
\end{center}

We numerically computed the non-vanishing lowest-lying zeros of our family $\mathcal{F}_f^{+}(X)$ with $f \in S_{8}^{\text{new}}(5,\text{principal})$, and the results are shown in Figure \ref{tab:table-mf5w8aa}. As predicted, the distribution of the non-vanishing lowest-lying zeros of even twists matches the eigenvalues of random matrices in $\SpO(\text{even})$ with characteristic polynomials evaluated at 1. We also get agreement between the lowest-zeros of odd twists and the eigenvalues of random matrices in $\SpO(\text{odd})$ with characteristic polynomials evaluated at 1 as both vanish. As predicted, the distribution of the second (non-vanishing) lowest-lying zeros of odd twists matches the eigenvalues of random matrices in $\SpO(\text{odd})$ with characteristic polynomials evaluated near 1, respectively. There is no visible attraction toward the origin either.

\begin{center}
\begin{figure}[htpb]
	\begin{tabular}{c c}
	Lowest zeros (even)
 & Non-vanishing lowest zeros (odd) \\
 \includegraphics[scale=0.45]{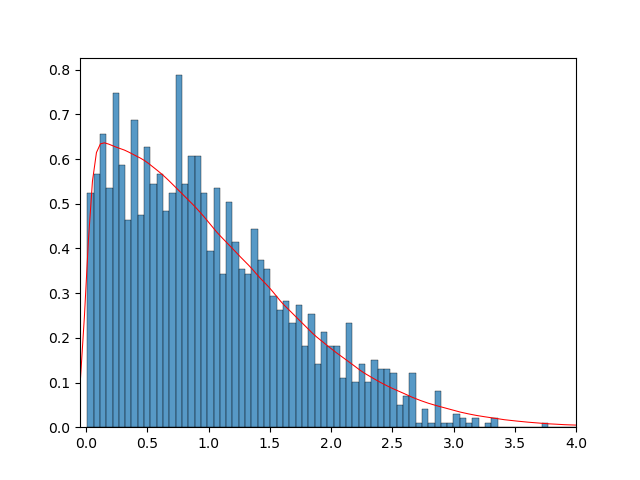}
 & \includegraphics[scale=0.45]{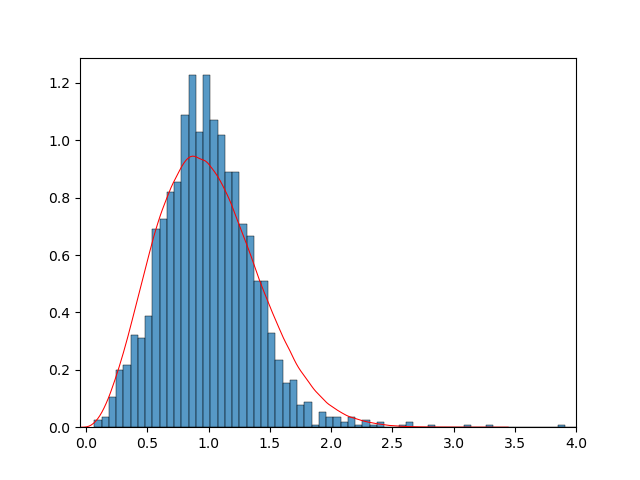}
	\end{tabular}
	\caption{The left histogram shows non-vanishing lowest-lying zeros of 1,951 even twists of \texttt{5.4.a.a} with discriminant up to 15,500; the red curve (left) shows the distribution of the first eigenvalues of 1,000,000 randomly generated $\SpO(18)$ matrices with characteristic polynomial evaluated at 1 without excision. The right histogram shows the second lowest-lying zeros of 1,951 odd twists of \texttt{5.4.a.a} with discriminant up to 15,500; the red curve (right) shows the distribution of the first eigenvalues of 1,000,000 randomly generated $\text{SO}(19)$ matrices with characteristic polynomial evaluated near 1. The data have been normalized to have mean 1.}
 \label{tab:table-mf5w4aa}
\end{figure}	
\end{center}

The question which natural arises is, do we observe a change in repulsion as we fix the level but vary the weight? Our example of a level 5, weight 4 newform with principal nebentype in Figure \ref{tab:table-mf5w4aa} suggests no such change. The authors in \cite{CSLPRRV24} respond to this question by computing the non-vanishing lowest-lying zeros for families of even twists for given newforms of fixed level 3 and weight varying between 6 and 10. In particular, their data for twists of \texttt{3.6.a.a}, \texttt{3.8.a.a}, and \texttt{3.10.a.a} support the prediction that no repulsion occurs at the origin for principal forms with weight 6 and greater. On a similar note, in addition to presenting data for \texttt{5.4.a.a}, we present data for \texttt{7.4.a.a} with two aims. First, we seek to clarify the lowest-lying zeros' behavior for principal forms of weight 4; and second, future authors may seek to determine if any changes in repulsion occurs as we fix the weight but vary the level.

\begin{center}
\begin{figure}[htpb]	
\begin{tabular}{c c}
Lowest zeros (even twists) & Lowest zeros (odd twists) \\
\includegraphics[scale=0.45]{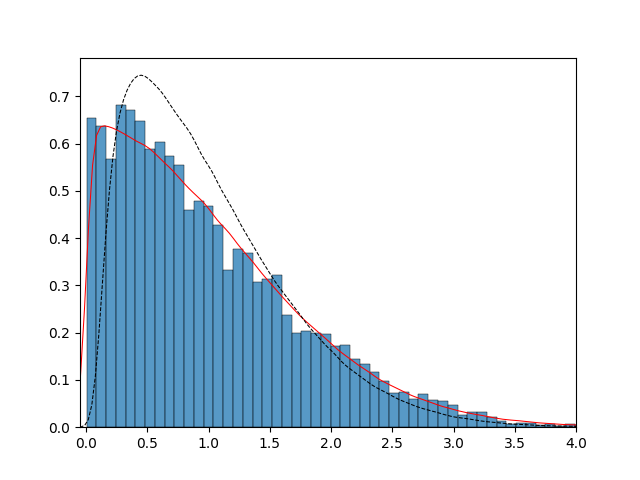} & \includegraphics[scale=0.45]{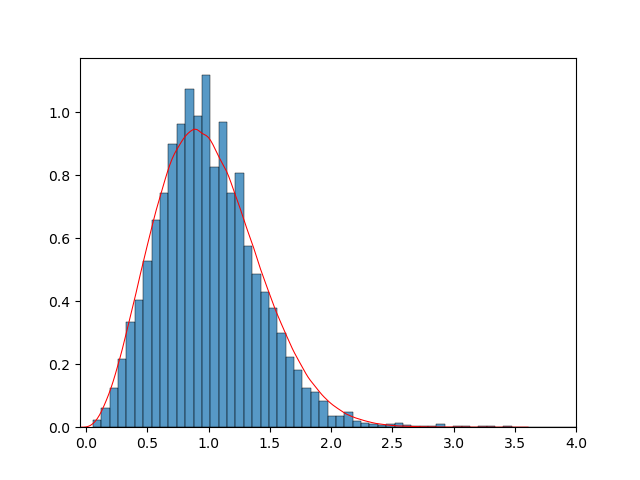}
\end{tabular}
\caption{The left histogram shows the distribution of lowest non-vanishing zeros for 5,463 odd twists of \texttt{7.4.a.a} with discriminant up to 41,000, and the red line (right) shows the distribution of the first eigenvalues of 1,000,000 randomly generated $\text{SO}(20)$ matrices with characteristic polynomial evaluated near 1. The right histogram shows second lowest-lying zeros of 5,463 odd twists of \texttt{7.4.a.a} with discriminant up to 41,000, and the red curve (right) shows the distribution of the first eigenvalues of 1,000,000 randomly generated $\text{SO}(21)$ matrices with characteristic polynomial evaluated near 1. The data have been normalized to have mean 1.}
\label{tab:zeros-mf7w4aa}
\end{figure}	
\end{center}

The sign of $f \in S_{4}^{\text{new}}(7,\text{principal})$ is $\epsilon_f=+1$. As shown in the left histogram of Figure \ref{tab:zeros-mf7w4aa}, there is little to no discernible repulsion at the origin for the even twists. There is no need for excision for this form as shown by the disagreement between the excised distribution (dotted line) in Figure \ref{tab:zeros-mf7w4aa}; that is, the non-excised random matrix model describes the distribution well. This is expected given the heuristic proposed in \cite{DHKMS12}. In the right histogram of Figure \ref{tab:zeros-mf7w4aa}, the second lowest-lying zeros of odd twists of \texttt{7.4.a.a} agree with the matrix model.

We remark a similar problem to the one arising for odd twists of \texttt{11.2.a.a} that concerns whether or not to model the second lowest-lying zeros or the lowest-lying non-vanishing zeros. On the right of Figure \ref{tab:attraction-mf11w2aa-mf7w4aa}, the red curve deviates from the data. In fact, we see the same attraction toward the origin for odd twists of \texttt{7.4.a.a} as that for odd twists of \texttt{11.2.a.a} on the left. Rather than implement a cutoff value, one might develop a new model for lowest non-vanishing zeros that accounts for this attraction by introducing a value that incorporates more first eigenvalues near the origin for families with high discriminant.

\subsection{Families with symplectic symmetry} 
\begin{center}
\begin{figure}[htpb]	
	\begin{tabular}{c c c}
		Lowest zeros ($\Delta=+1$) & Lowest zeros ($\Delta=-1$)\\
	\includegraphics[scale=0.45]{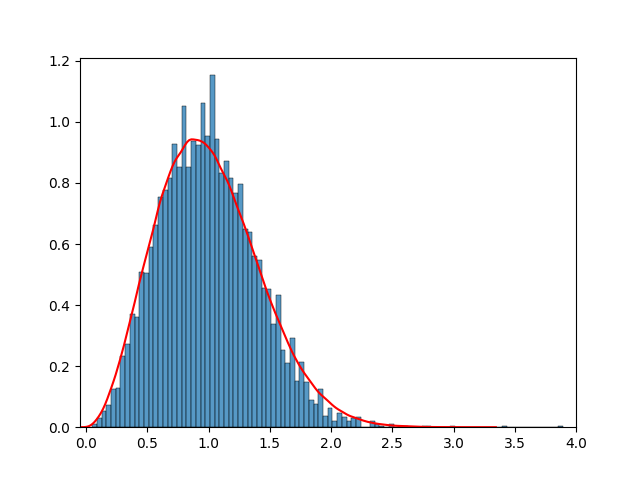} & \includegraphics[scale=0.45]{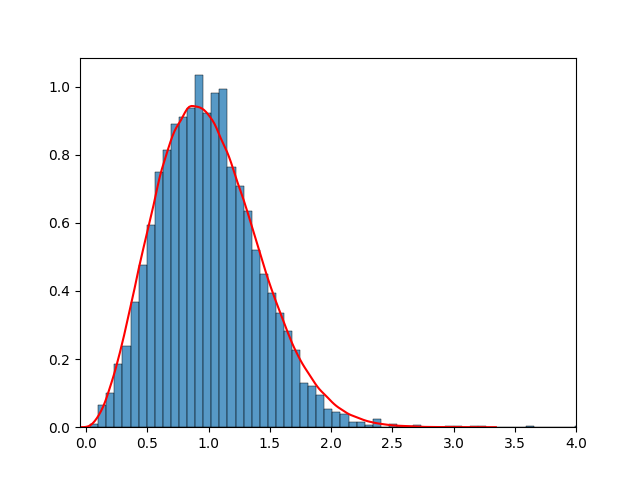} 
	\end{tabular}
	\caption{The left histogram shows the distribution of lowest-lying zeros for 5,458 twists of \texttt{3.7.b.a} with choice $\Delta =+1$ and discriminant up to 47,881, and the red line (left) shows the distribution of first eigenvalues of 1,000,000 randomly generated $\USp(20)$ matrices with characteristic polynomial evaluated at 1. The right histogram shows the distribution of lowest zeros for 5,726 twists of \texttt{3.7.b.a} with choice $\Delta =-1$ and discriminant up to 50,237, and the red line (right) shows the distribution of first eigenvalues of 1,000,000 randomly generated $\USp(20)$ matrices with characteristic polynomial evaluated at 1. The data have been normalized to have mean 1.}
 \label{tab:zeros-mf3w7ba}
\end{figure}	
\end{center}

\begin{center}
\begin{figure}[htpb]	
	\begin{tabular}{c c c}
		Lowest zeros ($\Delta=+1$) & Lowest zeros ($\Delta=-1$)\\
	\includegraphics[scale=0.45]{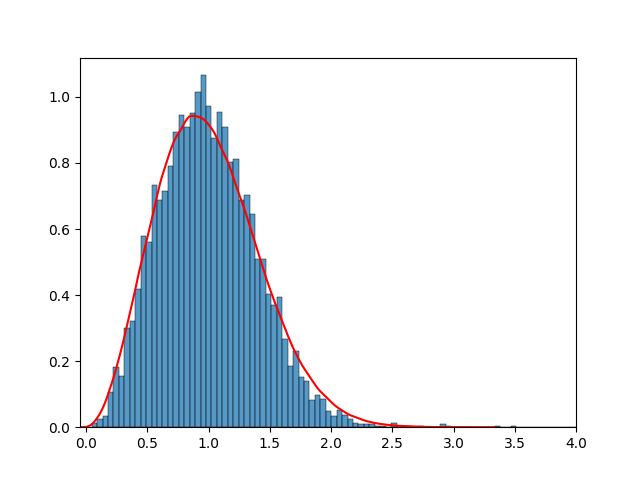} & \includegraphics[scale=0.45]{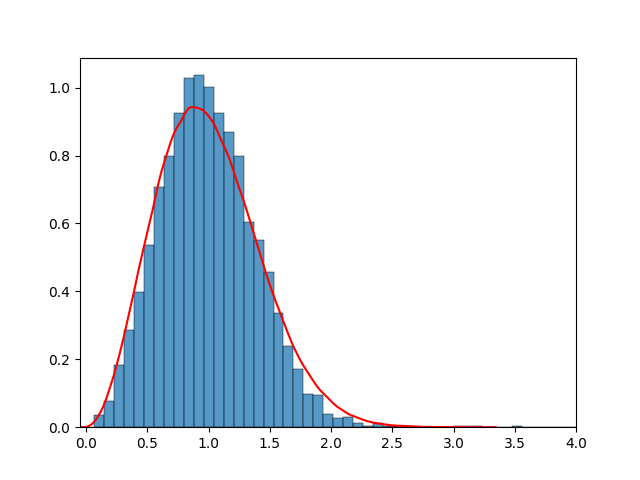} 
	\end{tabular}
	\caption{The left histogram shows the distribution of lowest-lying zeros for 5,467 twists of \texttt{7.3.b.a} with choice $\Delta =+1$ and discriminant up to 41,197, and the red line (left) shows the distribution of first eigenvalues of 1,000,000 randomly generated $\USp(20)$ matrices with characteristic polynomial evaluated at 1. The right histogram shows the distribution of lowest zeros for 5,486 twists of \texttt{7.3.b.a} with choice $\Delta =-1$ and discriminant up to 41,196, and the red line (right) shows the distribution of first eigenvalues of 1,000,000 randomly generated $\USp(20)$ matrices with characteristic polynomial evaluated at 1. The data have been normalized to have mean 1.}
 \label{tab:zeros-mf7w3ba}
\end{figure}	
\end{center}

As shown in Figure \ref{tab:zeros-mf3w7ba}, the lowest-lying zeros of $\mathcal{F}_f^{+}(X)$ for $f \in S_{7}^{\text{new}}(3,\text{self-CM})$ with sign $\epsilon_f=+1$ follow the predicted symplectic symmetry. The theory and the numerical results align as predicted. In Figure \ref{tab:zeros-mf7w3ba}, we also verify this behavior holds for \texttt{7.3.b.a}.

\subsection{Families with unitary symmetry}
\label{sec:families-with-unitary-symmetry}
\begin{center}
\begin{figure}[htpb]	
\begin{tabular}{c c}
Lowest zeros (twists) \\
\includegraphics[scale=0.55]{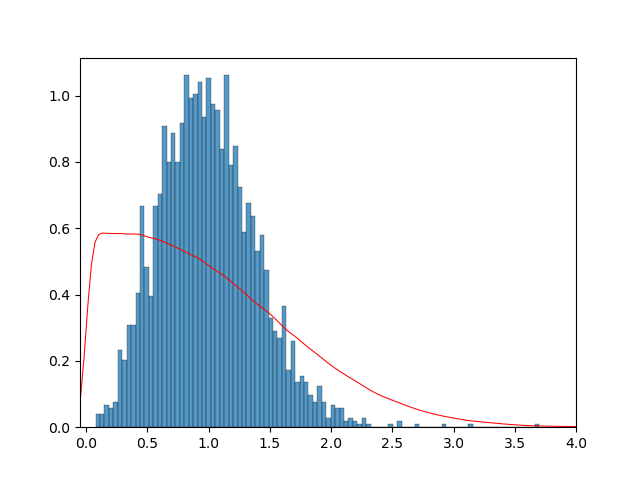} 
\end{tabular}
\caption{The histogram (blue) shows the distribution of the 2,860 lowest-lying zeros of the family of the generic form \texttt{11.7.b.b} with discriminant up to 10,277; the red line shows the distribution of the eigenvalues of 1,000,000 random matrices from $\U(9)$ whose characteristic polynomials are evaluated at 1. The data have been normalized to have mean 1.}
\label{tab:zeros-mf11w7bb}
\end{figure}	
\end{center}

As shown in Figure \ref{tab:zeros-mf11w7bb}, the lowest-lying zeros of our family associated to the generic form \texttt{11.7.b.b} does not follow the predicted unitary distribution. The distribution of the low-lying zeros seems to match the distribution of the first eigenvalues of 1,000,000 numerically generated symplectic matrices. This means we recovered self-CM behavior from a generic form. Note the form \texttt{11.7.b.a} is, in fact, self-CM. Hence, a form with predicted unitary symmetry can have a different predicted symmetry under certain conditions—which are yet to be determined. The deviating behavior showcased in Figure \ref{tab:zeros-mf11w7bb} may be explained by the fact that the one-level density for the unitary ensemble showcases no oscillatory behavior as it equals 1. Hence, there is no possibility of extracting any arithmetic nuance. In particular, the data suggests certain generic forms that have unitary symmetry would restrict to have symmetry of another ensemble.

\begin{center}
\begin{figure}[htpb]	
	\begin{tabular}{c c c}
		Lowest zeros ($\Delta=+1$) & Lowest zeros ($\Delta=-1$)\\
	\includegraphics[scale=0.45]{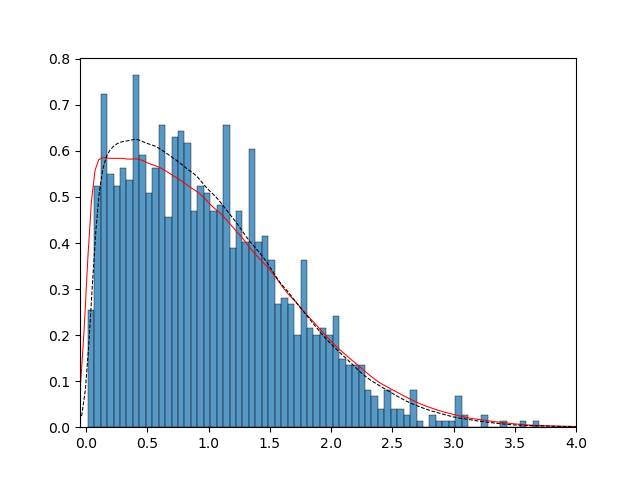} & \includegraphics[scale=0.45]{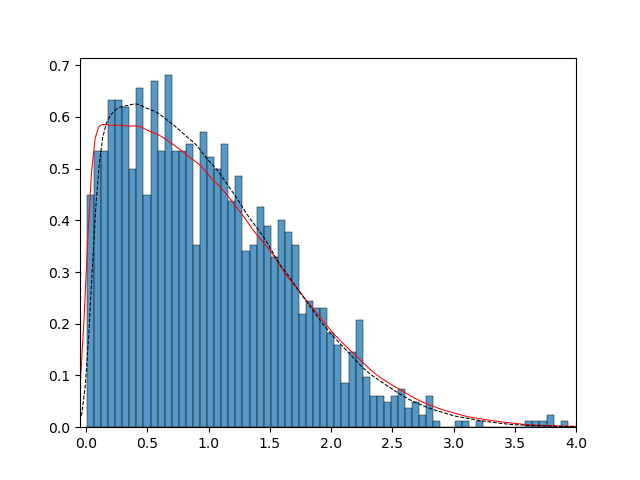} 
	\end{tabular}
	\caption{The left histogram shows the distribution of lowest-lying zeros for 1,048 twists of \texttt{13.2.e.a} with choice $\Delta =+1$ and discriminant up to 7,433; the red line (left) shows the distribution of first eigenvalues of 1,000,000 randomly generated $\U(9)$ matrices with characteristic polynomial evaluated at 1; the black dotted curve (left) is the same distribution but with excision of cutoff value 1/16. The right histogram shows the distribution of lowest zeros for 1,049 twists of \texttt{13.2.e.a} with choice $\Delta =-1$ with discriminant up to 7,429; the red line (right) shows the distribution of first eigenvalues of 1,000,000 randomly generated $\U(9)$ matrices with characteristic polynomial evaluated at 1; the black dotted curve (right) is the same distribution but with excision of cutoff value 1/16. The data have been normalized to have mean 1.}
 \label{tab:even-odd-mf13w2ea}
\end{figure}	
\end{center}

\begin{center}
\begin{figure}[htpb]	
\begin{tabular}{c c c}
Lowest zeros (twists) & Lowest zeros (twists) \\
\includegraphics[scale=0.45]{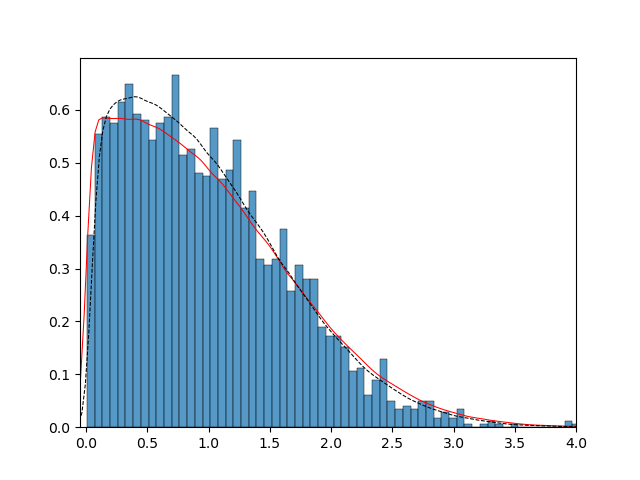} & \includegraphics[scale=0.45]{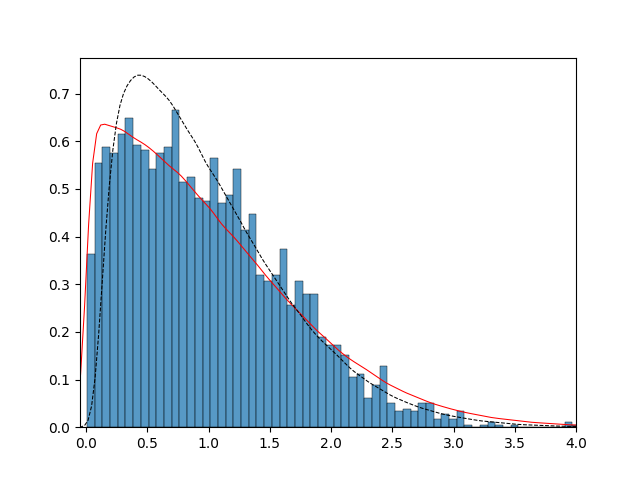} 
\end{tabular}
\caption{The left histogram shows the distribution of lowest-lying zeros for 2,097 twists of \texttt{13.2.e.a} with discriminant up to 7,500; the red curve (left) shows the distribution of first eigenvalues of 1,000,000 randomly generated matrices from $\U(9)$ with characteristic polynomial evaluated at 1; the black dotted curve (left) is the same distribution but with excision of cutoff value 1/16. We varied the excision threshold numerically to obtain the optimal fit. The right histogram shows the same twists as the left histogram but has the red line (right) showing the distribution of first eigenvalues of 1,000,000 randomly generated matrices from $\SpO(18)$; the black dotted curve (right) is the same distribution but with excision of cutoff value 1/64. The data have been normalized to have mean 1.}
 \label{tab:zeros-mf13w2ea-U-SO}
\end{figure}	
\end{center}

The next generic form we considered was \texttt{13.2.e.a}, and the histograms are presented in Figure \ref{tab:zeros-mf13w2ea-U-SO}. We took even and odd twists by setting $\psi_d(M)$ equal to either $+1$ or $-1$ to see if we recovered $\SpO(\text{even})$ and $\SpO(\text{odd})$ symmetry, respectively. The right histogram in Figure \ref{tab:zeros-mf13w2ea-U-SO} shows how well $\SpO(\text{even})$ random matrices model the twists of $\texttt{13.2.e.a}$. In Figure \ref{tab:even-odd-mf13w2ea}, we did not recover $\SpO(\text{odd})$ symmetry, which is expected given that generic forms should not be influenced by parity of the twist. 
In Figure \ref{tab:zeros-mf13w2ea-U-SO}, we see how well the first eigenvalues of 1,000,000 numerically generated of both unitary and special orthogonal (even) matrices, normalized to have mean 1, model the lowest-lying zeros of our family. This indistinguishability between the distributions presents difficulty for numerically determining the predicted ensemble for a generic form. 

We also see repulsion at the origin for \texttt{13.2.e.a}. The black dotted curve (left) is the distribution of random matrices from $\U(18)$ but with excision that was varied numerically to find the optimal fit. This might indicate the excision present in the model should be extended to all weight 2 forms regardless of principality of the nebentype. Note when we vary the excision threshold for $\SpO(18)$, a sharp incline appears. The sharp incline given by the black dotted curve (right) in Figure \ref{tab:zeros-mf13w2ea-U-SO} does not model the data well whereas the black dotted curve (left) still models the data well. Perceiving this difference gives us an approach to distinguish between forms with corresponding matrix group that is either unitary or special orthogonal.

\begin{center}
\begin{figure}[htpb]	
\begin{tabular}{c c}
Lowest zeros (twists) & Lowest zeros (twists) \\
\includegraphics[scale=0.45]{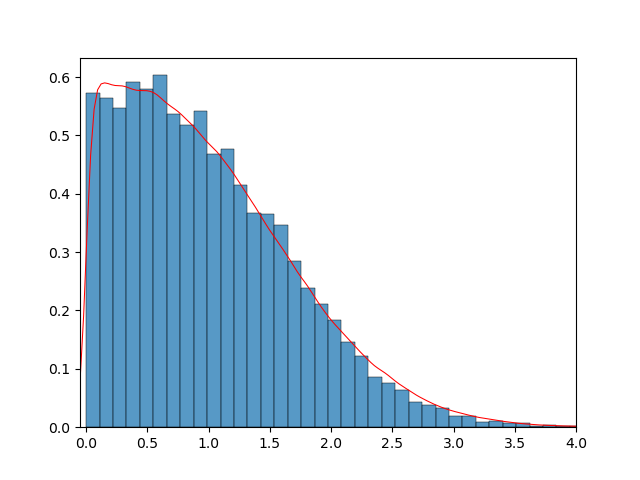} & \includegraphics[scale=0.45]{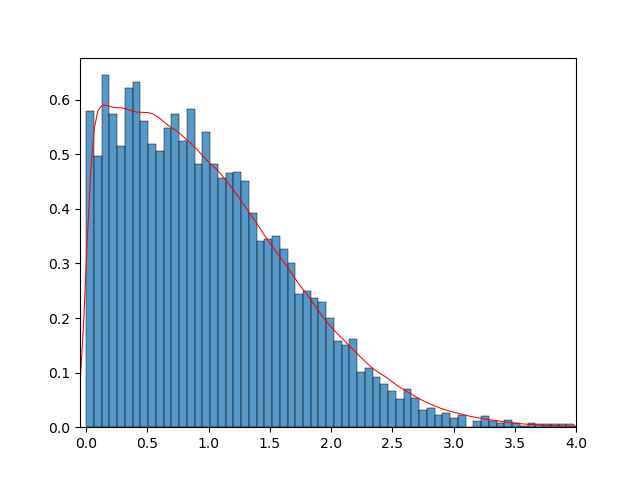}
\end{tabular}
\caption{The histogram (blue, left) shows the distribution of the 10,334 lowest-lying zeros of the family of the generic form \texttt{7.4.c.a} with discriminant up to 38,869; the red line (left) shows the distribution of the eigenvalues of 1,000,000 random matrices from the unitary ensemble $\U(10)$ whose characteristic polynomials are evaluated at 1. The histogram (blue, right) shows the distribution of the 11,889 lowest-lying zeros of the family of the generic form \texttt{17.2.d.a} with discriminant up to 25,105; the red line (left) shows the distribution of the eigenvalues of 1,000,000 random matrices from the unitary ensemble $\U(10)$ whose characteristic polynomials are evaluated at 1. The data have been normalized to have mean 1.}
\label{tab:zeros-mf7w4ca-mf17w2da}
\end{figure}	
\end{center}

In Figure \ref{tab:zeros-mf7w4ca-mf17w2da}, we provide the distribution for the lowest-lying zeros of the family for \texttt{7.4.c.a} to investigate if we see any repulsion for a weight 4 generic form. We do not. In the same figure, we also provide the distribution for the lowest-lying zeros of the family for \texttt{17.2.d.a} and we see no repulsion even though it is a (generic) form of weight 2. The question then is, what else controls discretization?

\appendix
\newpage
\section{Counting fundamental discriminants}\label{appendix:counting_fundamental_discriminants}
\begin{lemma}\label{lem:funcount}
  Let $\mathcal{D}^{+}_{f}(X)$ denote the set of even fundamental discriminants $d$ satisfying
  $0\leq d\leq X$, and fix an odd prime $M$ with $\Delta\in\{\pm1\}$. Then
\begin{equation}\label{eq:funasymp1}
  \sum_{\substack{d\in\mathcal D(X)\\(d/M) = \Delta}}1
  \ =\ \frac3{\pi^2}\frac M{2(M+1)}X+O(X^{1/2}),\end{equation}
  where $(\cdot/M)$ denotes the Legendre symbol.
  Moreover, fix $0< U<M$ an integer. Then
\begin{equation}\label{eq:funasymp2}\sum_{\substack{d\in\mathcal D(X)\\d\equiv U(M)}}1
\ =\ \frac3{\pi^2}\frac M{M^2-1}X+O(X^{1/2}).\end{equation}
\end{lemma}
\begin{proof}
We consider \eqref{eq:funasymp1} first. When $\Delta = \pm 1$, this is Lemma A.1 of~\cite{HMM11}.
For~\eqref{eq:funasymp2}, we follow a similar approach to~\cite[Lemma~B.1]{Mil08}. We have
\begin{equation}\sum_{\substack{d\in\mathcal D(X)\\d\equiv U(M)}}1\ =\ \sum_{\substack{0<d\leq X\\d\equiv1(4)\\d\equiv U(M)}}\mu(d)^2
+\sum_{\substack{0<d\leq X\\4\mid d\\d/4\equiv2(4)\\d\equiv U(M)}}\mu(d/4)^2
+\sum_{\substack{0<d\leq X\\4\mid d\\d/4\equiv3(4)\\d\equiv U(M)}}\mu(d/4)^2.\end{equation}
The Chinese remainder theorem gives
\begin{equation}\label{eq:chinese}
\sum_{\substack{d\in\mathcal D(X)\\d\equiv U\pmod {M}}}1\ =\ \sum_{\substack{0<d\leq X\\d\equiv A\pmod{ 4M}}}\mu(d)^2
+\sum_{\substack{0<d\leq X/4\\d\equiv A'\pmod{ 4M}}}\mu(d)^2
+\sum_{\substack{0<d\leq X/4\\d\equiv A''\pmod {4M}}}\mu(d)^2,\end{equation}
where $(A,4M)=1=(A'',4M)$ and $(A',4M)=2$. Recall
\begin{equation}\mu(n)^2=\sum_{m^2\mid n}\mu(m)\label{eq: square to mobius potato}\end{equation} and for all $m\in \mathbb{Z}$, we have $(m^2,4M)\neq 2$, but $A'\equiv2(4)$. For $d\equiv A(4M)$, \eqref{eq: square to mobius potato} gives 
\begin{equation}
\sum_{\substack{0<d\leq X\\d\equiv A(4M)}}\mu(d)^2\ =\ \sum_{\substack{0<d\leq X\\d\equiv A(4M)}}\sum_{m^2\mid d}\mu(m).
\end{equation}
Exchanging the order of summation and simplifying yields
\begin{align}
    &\sum_{\substack{0<d\leq X\\d\equiv A(4M)}}\mu(d)^2 \ = \ \sum_{m\leq X^{1/2}}\mu(m)\sum_{\substack{d'\leq X/m^2\\d'm^2\equiv A(4M)}}1 = \sum_{\substack{m\leq X^{1/2}(m^2,4M)\mid A}}\mu(m)
    \Big(\frac X{m^2}\frac{(m^2,4M)}{4M}+O(1)\Big) \nonumber \\
    &\quad =\ \frac{X}{4M}\sum_{\substack{m\leq X^{1/2}\\(m^2,4M)\mid A}}
    \frac{\mu(m)(m^2,4M)}{m^2}+O(X^{1/2}) \ =\ \frac X{4M}\sum_{\substack{m\leq X^{1/2}\\(m^2,4M)=1}}\frac{\mu(m)}{m^2}
    +O(X^{1/2}) \nonumber \\
    &\quad =\ \frac X{4M}\frac{(1-\frac14)^{-1}(1-\frac1{M^2})^{-1}}
    {\zeta(2)}+O(X^{1/2}) \ =\ \frac{2}{\pi^2}\frac{M}{M^2-1}X+O(X^{1/2}).
\end{align}

An analogous analysis of $d\equiv A'(4M)$ and $d\equiv A''(4M)$ yields the final result upon substituting them into~\eqref{eq:chinese}.
\end{proof}

\section{Calculating the mean density over fundamental discriminants} \label{app:mean-density-over-d}
\begin{lemma}
    \begin{equation}
    \sum_{d\in \mathcal{D}_f^{+}(X)} \bigg(\frac{\sqrt{M} |d|}{2\pi}\bigg)^{-2\pi i \tau/R}\ =\ |\mathcal{D}_f^{+}(X)| e^{-2\pi i \tau - 2\pi i \tau / R} \bigg(1-\frac{2\pi i 
    \tau}{R}\bigg)^{-1} + O(X^{1/2})
.\end{equation}
\end{lemma}

\begin{proof}
Recall $R = \log\paren{\frac{\sqrt{M}X}{2\pi e}}$. Then
\begin{align}
   \sum_{d\in \mathcal{D}_f^{+}(X)}\bigg( \frac{\sqrt{M} |d|}{2\pi} \bigg) ^{-2\pi i \tau /R}\ &=\ \sum_{d\in \mathcal{D}_f^{+}(X)} \exp\bigg(-\frac{2\pi i \tau}{R}\log(\sqrt{M} / 2\pi)\bigg)d^{-2\pi i \tau / R} \nonumber \\
   &=\ \exp\bigg(-2\pi i \tau + \frac{2\pi i \tau \log X}{R} - \frac{2\pi i \tau}{R}\bigg)\sum_{d\in \mathcal{D}_f^{+}(X)} d^{-2\pi i \tau/R} \label{eq: total sum to sub stuff into}
.\end{align}

We find an expression for $\sum_{d\in \mathcal{D}_f^{+}(X)} d^{-2\pi i \tau/R}$ using summation by parts. Put 
\[
    a_n\ =\ \begin{cases}
        1 & \qquad n\in \mathcal{D}_f^{+}(X) \\
        0 & \qquad \text{otherwise}
    \end{cases}
\]
and $\phi(u) = u^{-2\pi i \tau/R}$ to get
\begin{equation}
    \sum_{d\in \mathcal{D}_f^{+}(X)} d^{-2\pi i \tau/R}\ =\ |\mathcal{D}_f^{+}(X)| X^{-2\pi i \tau/R} - \int_1^X |\mathcal{D}_f(u)| u^{-2\pi i \tau/R} \frac{-2\pi i \tau}{R} d^*u\label{eq: sum d powers},
\end{equation}
where $d^*u = du/u$. We consider $f = \overline{f}$. Substituting the expression from Lemma \ref{lem:funcount} into \eqref{eq: sum d powers}, we get
\begin{equation}
    \sum_{d\in \mathcal{D}_f^{+}(X)} d^{-2\pi i \tau/R}= |\mathcal{D}_f^{+}(X)| X^{-2\pi i \tau/R} - \int_1^X \bigg(\frac{3 M u}{2\pi^2(M+1)} + O(u^{1/2})\bigg) u^{-2\pi i \tau/R} \frac{-2\pi i \tau}{R} d^*u. \label{eq: sum d integral}
\end{equation}
We have the following bound on the error term:
\begin{equation}\label{eq: error_term_cabbage}
    \left|\int_1^X O(u^{1/2}) u^{-2\pi i \tau/R} \frac{-2\pi i \tau}{R} d^*u \right|\ \leq\  C \int_1^X \left|u^{1/2}  u^{-2\pi i \tau/R} \frac{-2\pi i \tau}{R}\frac{1}{u} \right|du \ =\ O(X^{1/2})
\end{equation}
for some constant $C > 0$.
Moreover, for the other part of the integral in \eqref{eq: sum d integral}, we get
\begin{align}
    -\int_1^X \frac{3M u}{2\pi^2 (M+1)} u^{-2\pi i \tau/R} \frac{-2\pi i \tau}{R}d^*u\ &=\ \frac{2\pi i \tau}{R} \frac{3 M}{2\pi^2(M+1)}\int_1^X u^{-2\pi i \tau/R}du \nonumber\\
    &=\ \frac{2\pi i \tau}{R} \frac{3 M X}{2\pi^2(M+1)} \frac{X^{- 2\pi i \tau/R}}{1 - \frac{2\pi i \tau}{R}} + O(1).
\end{align}
Hence, using Lemma \ref{lem:funcount}, we obtain
\begin{equation}
    -\int_1^X \frac{3M u}{2\pi^2(M+1)} u^{-2\pi i \tau/R} \frac{-2\pi i \tau}{R}d^*u\ =\ \frac{2\pi i \tau}{R} (|\mathcal{D}_f^{+}(X)| + O(X^{1/2})) \frac{X^{-2\pi i \tau/R}}{1 - \frac{2\pi i \tau}{R}} + O(1).\quad\quad\quad
\end{equation}
Again, $\tau\in \mathbb{R}$ gives
\begin{equation}\label{eq: integral main spinach}
    -\int_1^X \frac{3M u}{2\pi^2(M+1)} u^{-2\pi i \tau/R}  \frac{-2\pi i \tau}{R}d^*u\ =\ \frac{2\pi i \tau}{R} |\mathcal{D}_f^{+}(X)|  \frac{X^{-2\pi i \tau/R}}{1 - \frac{2\pi i \tau}{R}} + O(X^{1/2}).
\end{equation}
Finally, substituting both \eqref{eq: integral main spinach} and \eqref{eq: error_term_cabbage} into \eqref{eq: sum d integral} and then substituting into \eqref{eq: total sum to sub stuff into} to get
\begin{align}
    \sum_{d\in \mathcal{D}_f^{+}(X)} \bigg(\frac{\sqrt{M} |d|}{2\pi}\bigg)^{-2\pi i \tau/R}\ =&\ \exp\bigg(-2\pi i \tau + \frac{2\pi i \tau \log X}{R} - \frac{2 \pi i \tau}{R}\bigg) \nonumber\\
    &\times\ \bigg(|\mathcal{D}_f^{+}(X)| \bigg(1-\frac{2\pi i \tau}{R}\bigg)^{-1} \exp\bigg(-\frac{2\pi i \tau}{R}\log X\bigg) + O(X^{1/2})\bigg).
\end{align}
Finally, since $\tau\in \mathbb{R}$, we arrive at 
\begin{equation}
    \sum_{d\in \mathcal{D}_f^{+}(X)} \bigg(\frac{\sqrt{M} |d|}{2\pi}\bigg)^{-2\pi i \tau/R}\ =\ |\mathcal{D}_f^{+}(X)| e^{-2\pi i \tau - 2\pi i \tau / R} \bigg(1-\frac{2\pi i 
    \tau}{R}\bigg)^{-1} + O(X^{1/2}).
\end{equation}
The case $f\neq \overline{f}$ is the same.
\end{proof}
\section{Calculating mean density over log of positive fundamental discriminants}\label{app:mean-density-over-log-d}

\begin{lemma}
Let $\mathcal D(X)$ denote the set of even fundamental discriminants $d$ satisfying
  $d\leq X$. Then,
    \begin{equation}\label{eq: beginning eggs and toast}
        \sum_{d\in \mathcal{D}_f^{+}(X)} \log \bigg( \frac{\sqrt{M} |d|}{2\pi} \bigg) = |\mathcal{D}_f^{+}(X)|\bigg(\log\bigg(\frac{\sqrt{M} X}{2\pi}\bigg) - 1\bigg) + O(X^{1/2}).
    \end{equation}
\end{lemma}

\begin{proof}
We perform summation by parts on $\sum_{d\in \mathcal{D}_f^{+}(X)}\log d$ with
\[
    a_n \ =\ \begin{cases}
        1 & \qquad n\in \mathcal{D}_f^{+}(X) \\
        0 & \qquad \text{otherwise}
    \end{cases}
\]
and $\phi(u) = \log u$ to get
\begin{equation}
     \sum_{d\in \mathcal{D}_f^{+}(X)} \log d \ =\ |\mathcal{D}_f^{+}(X)| \log X - \int_1^X |\mathcal{D}_f(u)| d^*u,
     \label{eq: sum of log d equation}
\end{equation}
where $d^*u = du/u$. We consider $f= \overline{f}$. Then by Lemma \ref{lem:funcount}, the equation \eqref{eq: sum of log d equation} becomes 
\begin{equation}
    \sum_{d\in \mathcal{D}_f^{+}(X)} \log d \ =\ |\mathcal{D}_f^{+}(X)| \log X - \int_1^X \bigg(\frac{3 M u}{2\pi^2(M+1)} + O(u^{1/2})\bigg) \, d^*u \label{eq: log d potato}.
\end{equation}
For the error term, we get
\begin{align}
    \left|\int_1^X O(u^{1/2}) \, d^*u\right| \ &\leq\ C\int_1^X  |u^{-1/2}|\, du \ =\ O(X^{1/2}) \label{eq: log error}.
\end{align}
The other part of the integral is 
\begin{align}
    -\int_1^X \frac{3 Mu}{2\pi^2(M+1)}d^*u\ &=\ -\int_1^X \frac{3M}{2\pi^2(M+1)}du
    \ =\ -|\mathcal{D}_f^{+}(X)| + O(X^{1/2})\label{eq: log other part of integral}.
\end{align}
Then, combining \eqref{eq: log error} and \eqref{eq: log other part of integral} turns \eqref{eq: log d potato} into 
\begin{equation}
    \sum_{d\in \mathcal{D}_f^{+}(X)} \log d =\ |\mathcal{D}_f^{+}(X)| \log X - |\mathcal{D}_f^{+}(X)| + O(X^{1/2}) \label{eq: almost done pepper salad}.
\end{equation}
Substituting \eqref{eq: almost done pepper salad} into \eqref{eq: beginning eggs and toast} gives
\begin{align}
    \sum_{d\in \mathcal{D}_f^{+}(X)} \log \bigg( \frac{\sqrt{M} |d|}{2\pi} \bigg)\ &=\ |\mathcal{D}_f^{+}(X)| \log\paren{\frac{\sqrt{M}}{2\pi}} + |\mathcal{D}_f^{+}(X)| \log X - |\mathcal{D}_f^{+}(X)| + O(X^{1/2}) \nonumber\\
    &=\ |\mathcal{D}_f^{+}(X)| \paren{\log\paren{\frac{\sqrt{M} X}{2\pi}} - 1} + O(X^{1/2}).
\end{align}
The case $f\neq \overline{f}$ is analogous.
\end{proof}

\section{Series expansions for pair-correlation}\label{appendix:pair_correlation_series}
We provide some details regarding the expansion of \eqref{beforeexpansion} in series for large $R$. We expand the integrand in \eqref{beforeexpansion} in powers of $1/R$; that is, we expand
\begin{align}
    1 + \frac{1}{R^2} &+ \frac{1}{2R^2} \bigg(\frac{L_{\star}'}{L_{\star}}\bigg)' \bigg(1+\frac{i\pi y}{R}, f_d \otimes \overline{f}_d\bigg) \nonumber \\
    &+ \frac{1}{2R^2c_{f_d}^2} \frac{e^{-2\pi i y(1 + 1/R)}}{1 - 2\pi i y /R} L\bigg(1+\frac{i\pi y}{R}, f_d \otimes \overline{f}_d\bigg)L\bigg(1-\frac{i\pi y}{R}, f_d \otimes \overline{f}_d\bigg) \mathscr{A}\bigg(\frac{i\pi y}{R}\bigg)  \nonumber\\
    &+ \frac{1}{2R^2} \mathscr{C}\bigg(1+\frac{i\pi y}{R}\bigg)
\end{align}
in $1/R$. First, we replace the $L$-functions arising from the Rankin-Selberg convolution of $f_d$ with its dual via the factorization
\begin{equation}
    L\left(1\pm\frac{i\pi y}{R}, f_d \otimes \overline{f}_d\right) \ = \ L\left(1\pm\frac{i\pi y}{R}, \text{ad}^2\,f_d\right) \zeta\left(1\pm\frac{i\pi y}{R}\right).
\end{equation}
Then, using the Laurent expansions of the Riemann zeta function about its pole at $s=1$, we observe that
\begin{align}
    \zeta\left(1+\frac{i\pi y}{R}\right) \zeta\left(1-\frac{i\pi y}{R}\right) \ &= \ \frac{R^2}{\pi^2 y^2} + \gamma^2 + 2\gamma_1 + O(R^{-2}).
\end{align}
Then, since the adjoint square $L$-function is entire, we may write its Taylor expansions about $s=1$; this yields
\begin{align}
    L\left(1+\frac{i\pi y}{R}, \text{ad}^2\,f_d\right) \ = \ &L\left(1, \text{ad}^2\,f_d\right) + L'\left(1, \text{ad}^2\,f_d\right) \frac{i\pi y}{R} - \frac{1}{2} L''\left(1, \text{ad}^2\,f_d\right) \frac{\pi^2 y^2}{R^2} \nonumber \\
    &- \frac{1}{6}L'''\left(1, \text{ad}^2\,f_d\right) \frac{i\pi^3y^3}{R^3} + O(R^{-4}) \\
    L\left(1-\frac{i\pi y}{R}, \text{ad}^2\,f_d\right) \ = \ &L\left(1, \text{ad}^2\,f_d\right) - L'\left(1, \text{ad}^2\,f_d\right) \frac{i\pi y}{R} - \frac{1}{2} L''\left(1, \text{ad}^2\,f_d\right) \frac{\pi^2 y^2}{R^2} \nonumber \\
    &+ \frac{1}{6}L'''\left(1, \text{ad}^2\,f_d\right) \frac{i\pi^3y^3}{R^3} + O(R^{-4}).
\end{align}
The product is then
\begin{align}
    L(1, \text{ad}^2\,f_d)^2 + \frac{\pi^2 y^2}{R^2} \bigg( L'(1,\text{ad}^2\,f_d)^2 - L(1, \text{ad}^2\,f_d) L''(1, \text{ad}^2\,f_d) \bigg) + O(R^{-4}).
\end{align}
We also write down the series expansion
\begin{equation}
    \frac{e^{-2\pi i y(1+1/R)}}{1 - 2\pi i y/R} \ = \ e^{-2\pi i y} \frac{e^{-2\pi i y/R}}{1-2\pi i y/R} \ = \ e^{-2\pi i y} \left(1 - \frac{2\pi^2 y^2}{R^2} - \frac{8i \pi^3 y^3}{3R^3} + O(R^{-4})\right).
\end{equation}
Then, the analyticity of $\mathscr{A}(\cdot)$ about the origin allows for the expansion
\begin{align}
    \mathscr{A}\left(\frac{i\pi y}{R}\right) \ &= \ \mathscr{A}(0) + \mathscr{A}'(0) \frac{i\pi y}{R} - \frac{\mathscr{A}''(0)}{2} \frac{\pi^2 y^2}{R^2} - \frac{\mathscr{A}'''(0)}{6} \frac{i\pi^3 y^3}{R^3} + O(R^{-4}).
\end{align}
Then, using the definition of the ramified $L$-function, we are able to derive the series expansion
\begin{align}
    \bigg(\frac{L_{\star}'}{L_{\star}}\bigg)'\bigg(1+\frac{i\pi y}{R}, f_d \otimes \overline{f}_d\bigg) \ & = \ -\frac{R^2}{\pi^2y^2} - (\gamma^2 + 2\gamma_1) + \bigg(\frac{L'}{L}\bigg)'(1, \text{ad}^2\,f_d) \nonumber \\
    & \qquad - \sum_{p\mid N} \bigg(\frac{L_p'}{L_p}\bigg)'(1) + O(R^{-2}).
\end{align}
Here, each $L_p$ is the appropriate local Archimedean factor (c.f. \cite{IK04}).
Finally, making use of the identity
\begin{equation}
    -\mathscr{B}(1) \ = \ \frac{\mathscr{A}''(0)}{2} + \bigg(\frac{L_{\text{ram}}'}{L_{\text{ram}}}\bigg)'(1, f_d \otimes \overline{f}_d)
\end{equation}
yields the result upon combining together the aforementioned facts.

Throughout the argument, the evenness of $g(y)$ in conjunction with the integral being taken over all of $\mathbb{R}$ allows us to eliminate terms in the series expansion that are odd in the variable $y$. In this manner, we obtain the quoted series expansion in Section \ref{subsect:series_pair_correlation}.

\newpage
\noindent


\end{document}